\pdfoutput=1
\documentclass[11pt,a4paper]{amsart}
\usepackage[T1]{fontenc}
\usepackage[utf8]{inputenc}
\usepackage{lmodern,microtype}
\usepackage[margin=26mm]{geometry}
\usepackage{amsmath,amssymb,amsthm,mathtools}
\usepackage{booktabs,array,longtable,enumitem,needspace}
\usepackage{xcolor}
\usepackage{combelow}
\usepackage[colorlinks=true,linkcolor=blue!40!black,citecolor=blue!40!black,urlcolor=blue!40!black]{hyperref}
\hypersetup{pdftitle={Geometric obstructions to lim Cohen-Macaulay sequences},pdfauthor={Cristian Anghel and Filip Chindea},pdfkeywords={lim Cohen-Macaulay sequences, local cohomology, Koszul homology, signature, Harder-Narasimhan polygons, mixed characteristic}}
\setlist{itemsep=3pt,topsep=5pt}
\newtheorem{theorem}{Theorem}[section]
\newtheorem{proposition}[theorem]{Proposition}
\newtheorem{lemma}[theorem]{Lemma}
\newtheorem{corollary}[theorem]{Corollary}
\theoremstyle{definition}

\theoremstyle{remark}
\newtheorem{remark}[theorem]{Remark}
\newtheorem{example}[theorem]{Example}
\newcommand{\CC}{\mathbf C}

\newcommand{\ZZ}{\mathbf Z}
\newcommand{\PP}{\mathbf P}
\newcommand{\OO}{\mathcal O}
\newcommand{\HH}{\mathrm H}

\newcommand{\mm}{\mathfrak m}
\newcommand{\nn}{\mathfrak n}
\newcommand{\rk}{\operatorname{rk}}
\newcommand{\length}{\operatorname{\ell}}
\newcommand{\Hom}{\operatorname{Hom}}
\newcommand{\RHom}{\operatorname{RHom}}
\newcommand{\Ext}{\operatorname{Ext}}
\newcommand{\Tor}{\operatorname{Tor}}
\newcommand{\depth}{\operatorname{depth}}
\newcommand{\Spec}{\operatorname{Spec}}
\newcommand{\Proj}{\operatorname{Proj}}
\newcommand{\coker}{\operatorname{coker}}

\title[Obstructions to lim Cohen--Macaulay sequences]{Geometric obstructions to lim Cohen--Macaulay sequences}
\author{Cristian Anghel}
\address{Simion Stoilow Institute of Mathematics of the Romanian Academy (IMAR), 21 Calea Grivi\cb{t}ei, 010702 Bucharest, Romania}
\email{cristian.anghel@imar.ro}
\author{Filip Chindea}
\address{Simion Stoilow Institute of Mathematics of the Romanian Academy (IMAR), 21 Calea Grivi\cb{t}ei, 010702 Bucharest, Romania}
\email{fchindea@imar.ro}
\date{30 November 2026}
\subjclass[2020]{Primary 13C14, 14J60; Secondary 13D45, 13D02, 13H15, 14F06, 14J29}
\keywords{lim Cohen--Macaulay sequences, local cohomology, Koszul homology, signature, Harder--Narasimhan polygons, mixed characteristic}
\begin{document}
\begin{abstract}
Lim Cohen--Macaulay sequences, introduced by Bhatt, Hochster and Ma,
are asymptotic substitutes for small
Cohen--Macaulay modules. Let $Y$ be a smooth connected complex
projective surface and $H$ an ample globally generated divisor with
$H^2<K_Y^2-8\chi(\OO_Y)$. We prove that the section ring of $(Y,H)$
admits no graded lim Cohen--Macaulay sequence, and that under an
additional upper bound on $K_Y\cdot H$ in terms of $H^2$ and
$K_Y^2-8\chi(\OO_Y)$, with no proportionality between $K_Y$ and $H$
assumed, its completion at the vertex admits no lim
Cohen--Macaulay sequence of arbitrary finite modules. The proof
combines lower bounds, linear in the rank, for the intermediate
cohomology of sheaves on $Y$ and for the local cohomology of modules
over this completion, with a duality argument that turns such
bounds into lower bounds for Koszul homology. Blown-up
Hirzebruch--Kummer covers of the plane branched along the Hesse, dual
Hesse and complete-quadrangle arrangements, and products of the dual
Hesse surfaces with projective spaces, give complete normal local
domains without lim Cohen--Macaulay sequences in every dimension at
least three, with residue field $\CC$ and an isolated singularity. For
every prime $p$ there are also mixed-characteristic examples, with
residue field $\overline{\mathbf F}_p$, in every dimension at least
four. In particular, complete local domains with algebraically closed
residue field need not admit lim Cohen--Macaulay sequences, in equal
characteristic zero or in mixed characteristic; this disproves the
existence conjecture emphasized by Hochster, in its
mixed-characteristic and characteristic-free forms. The weakly
lim Cohen--Macaulay condition is studied in a companion paper.
\end{abstract}
\maketitle
\tableofcontents

\section{Introduction}\label{sec:introduction}

Lim Cohen--Macaulay sequences were introduced by Bhatt, Hochster and
Ma \cite{BHM} (see also \cite[Section~5]{Ho}) as an asymptotic
substitute for small Cohen--Macaulay modules.
For a local ring $(B,\mathfrak m)$, let $\nu_B(M)$ denote the minimal
number of generators of a finite module. If $\mathbf x$ is a system of
parameters, set
\[
 \kappa_{\mathbf x}(M)=\max_{i\ge1}\length \HH_i(\mathbf x;M),
 \qquad
 \chi_1(\mathbf x;M)=\sum_{i\ge1}(-1)^{i-1}\length \HH_i(\mathbf x;M).
\]
A sequence $(M_j)$ of nonzero finite $B$-modules of dimension
$\dim B$ is \emph{lim Cohen--Macaulay} if
$\kappa_{\mathbf x}(M_j)=o(\nu_B(M_j))$ as $j\to\infty$, for one,
equivalently every, system of parameters \cite[Theorem~4.6]{BHM}.
For graded rings and modules we use a graded variant, defined here, with
homogeneous systems of parameters, minimal numbers of homogeneous
generators, and dimensions over the degree-zero field in place of
lengths; after localization at the irrelevant ideal it is the
condition above (Section~\ref{ex:graded-subsection}). The \emph{weakly lim
Cohen--Macaulay} condition replaces $\kappa_{\mathbf x}$ by $\chi_1$
\cite{Ma,BHM}; when the distinction matters, we call the first
condition \emph{ordinary}. The two conditions differ because the
positive Koszul homology groups can cancel in their alternating
sum.

Bhatt, Hochster and Ma show that lim Cohen--Macaulay sequences exist
over $F$-finite local rings of positive characteristic
\cite[Theorem~5.4]{BHM}, and that their existence in mixed
characteristic would imply Serre's positivity conjecture for
intersection multiplicities over all regular local rings
\cite[Theorem~6.1 and Corollary~6.3]{BHM}. Complete local domains
without small Cohen--Macaulay modules are now known \cite{OAI,An3}, so
it is natural to ask whether the asymptotic substitute always exists.
Hochster emphasizes the conjecture that complete local domains have
lim Cohen--Macaulay sequences \cite[Section~1]{Ho}, and states it for
complete local domains of mixed characteristic with algebraically
closed residue field as \cite[Conjectures~5.2 and~13.2]{Ho}. The form
for every complete local domain with algebraically closed residue
field is stated as \cite[Conjecture~1.1]{IS}, where it is attributed
to Bhatt, Hochster and Ma; see also \cite{HS}.
Theorem~\ref{main:dimensions} shows that the asymptotic substitute
need not exist, and so disproves each of these forms. It gives complete
normal local domains with residue field $\CC$ and an isolated
singularity in every dimension at least three, and
mixed-characteristic domains with residue field
$\overline{\mathbf F}_p$ in every dimension at least four, that admit
no lim Cohen--Macaulay sequence. Assuming \cite[Main claim~1.1]{Ch},
nonexistence had already been asserted for a class of
three-dimensional examples \cite[Corollary~3.2]{J};
Section~\ref{intro:related} compares it with ours. Hochster's weak
forms \cite[Conjectures~10.1 and~13.1]{Ho}, which concern weakly lim
Cohen--Macaulay sequences, remain open.

Nonexistence does not follow formally from the absence of Ulrich or
arithmetically Cohen--Macaulay bundles \cite{An1,An2}. The absence of
modules with vanishing intermediate cohomology does not exclude
sequences whose positive Koszul homology, normalized by the number of
generators, tends to zero. Moreover, a module with large local
cohomology can have negligible positive Koszul homology
(Example~\ref{transfer:example}), so a lower bound for local cohomology
must be transferred to Koszul homology through a replacement
module. We prove geometric
inequalities for intermediate cohomology and a quantitative dualizing
replacement that together give the obstruction.

Unless stated otherwise, $Y$ is a smooth connected complex
projective surface and $H$ is ample and globally generated. Write
\[
 d=H^2,\qquad \sigma=K_Y^2-8\chi(\OO_Y)
 =\tfrac13\bigl(K_Y^2-2c_2(Y)\bigr),\qquad
 \lambda=\frac{K_Y\cdot H}{d}.
\]
By the Hirzebruch signature theorem, $\sigma$ is the signature of
$Y$. Let
\[
 R=R(Y,H)=\bigoplus_{t\ge0}\HH^0(Y,\OO_Y(tH)),\qquad
 \mathfrak m=R_+,\qquad \widehat R=\widehat{R_{\mathfrak m}}.
\]
We use the full section ring; generation in degree one is not assumed.
It is a finitely generated normal domain of dimension three, and
$\widehat R$ is a complete normal local domain with residue field
$\CC$, regular outside its closed point (Proposition~\ref{loc:model};
see also \cite[Lemma~1.2]{An2} and \cite[Proposition~6.1]{An3}).

All slopes are Mumford--Takemoto slopes, and Harder--Narasimhan
factors are torsion-free and slope-semistable \cite[Theorem~1.6.7]{HL}.
On $Y$ the normalized slope of a sheaf $E$ of rank $\rho>0$ is
$\bar\mu(E)=c_1(E)\cdot H/(\rho d)$, so
$\bar\mu(E(tH))=\bar\mu(E)+t$ and $\bar\mu(\omega_Y)=\lambda$.
The symbol $\equiv$ denotes numerical equivalence.

\subsection{The geometric inequalities}
For a coherent sheaf $E$, put
$h^1_*(E)=\sum_{t\in\ZZ}h^1(Y,E(tH))$, allowing the value $\infty$.

\begin{theorem}[Graded cohomology bound]\label{main:graded}
For every coherent sheaf $E$ of positive rank on $Y$,
\[
 h^1_*(E)\ge\frac{\sigma-d}{8}\rk E.
\]
The inequality is strict if $d<\sigma$. More generally, for every integer
$q\ge1$ one may replace the constant by $q(\sigma-q^2d)/8$.
Equivalently, every finite graded $R$-module $M$ satisfies
$\dim_{\CC}\HH^2_{\mathfrak m}(M)\ge(\sigma-d)\rk_R M/8$, and the
same holds with the constant $q(\sigma-q^2d)/8$.
\end{theorem}

The proof compares Harder--Narasimhan polygons after a finite map
$Y\to\mathbf P^2$. An elementary estimate for $h^1_*$ on the plane
replaces Horrocks' splitting criterion, which \cite{An2} uses for
arithmetically Cohen--Macaulay bundles. Coherent sheaves reduce to
vector bundles (Corollary~\ref{surf:coherent-inequality}(b)). For vector bundles,
\cite[Corollary~8.3(a)]{AC}, applied as in
Corollary~\ref{surf:veronese}, gives the larger constants
$q^3d\,L(\sigma/(q^2d))\ge q(\sigma-q^2d)/6$ for $q^2d<\sigma$
(Remark~\ref{surf:stronger-bundle-bound}); our proof is independent
of it and also gives the correction terms of
Theorem~\ref{surf:refined-inequality}.

\begin{theorem}[Completed-local cohomology bound]\label{main:local}
Let $q\ge1$ be an integer and $M$ a finite $\widehat R$-module.
\begin{enumerate}[label=\textup{(\alph*)}]
\item If $K_Y\cdot H<5qd$, then
\[
 \length \HH^2_{\mathfrak m}(M)\ge
 \frac{q(\sigma-q^2d)}8\rk_{\widehat R}M.
\]
\item If $K_Y\cdot H=5qd$, then
\[
 \length \HH^2_{\mathfrak m}(M)\ge
 q\left(\frac{\sigma-q^2d}{8}-2\right)\rk_{\widehat R}M.
\]
\end{enumerate}
The length is allowed to be infinite; it is finite for reflexive modules.
The same inequalities hold over $R_{\mathfrak m}$. We write
$c_q=q(\sigma-q^2d)/8$ if $K_Y\cdot H<5qd$, and
$c_q=q((\sigma-q^2d)/8-2)$ if $K_Y\cdot H=5qd$.
\end{theorem}

For $K_Y\equiv4H$, $q=1$ and reflexive $M$, the first estimate is
the inequality proposed in \cite[Main claim~1.1]{Ch}, and Chen's
constant $15H^2/8-\chi(\OO_Y)$ equals $(\sigma-d)/8$ when
$K_Y\equiv4H$. Section~\ref{loc:section} gives a complete proof of
Theorem~\ref{main:local}, so part (a) contains that claim. We follow Chen's
elementary modification and persistence argument, replace the
proportionality $K_Y\equiv4H$ by the slope condition
$K_Y\cdot H\le5qd$, and add the endpoint and Veronese refinements.
Theorem~\ref{main:graded} is proved in Section~\ref{surf:section}
(Corollaries~\ref{surf:coherent-inequality},
\ref{surf:graded-inequality} and~\ref{surf:veronese}), and
Theorem~\ref{main:local} is Theorem~\ref{loc:inequality}.

\subsection{Obstructions to ordinary lim Cohen--Macaulay sequences}
\begin{theorem}[Graded and local criteria]\label{main:criteria}
If $d<\sigma$, then $R(Y,H)$ admits no lim Cohen--Macaulay sequence
of finite graded modules. If there is an integer $q\ge1$ such that
\begin{equation}\label{main:local-condition}
 \bigl(q^2d<\sigma\ \text{and}\ K_Y\cdot H<5qd\bigr)
 \quad\text{or}\quad
 \bigl(q^2d+16<\sigma\ \text{and}\ K_Y\cdot H\le5qd\bigr),
\end{equation}
then neither $\widehat R$ nor $R_{\mathfrak m}$ admits a lim
Cohen--Macaulay sequence of arbitrary finite modules.
In each case the obstruction is quantitative. Let $\mathbf x$ be
three sections of $\OO_Y(H)$ without a common zero in the graded
case, and any system of parameters of $\widehat R$ in the local
case. There is a constant $\eta_{\mathbf x}>0$ such that every
nonzero finite module $M$, graded in the graded case and over
$\widehat R$ in the local case, satisfies
\[
 \kappa_{\mathbf x}(M)\ge\eta_{\mathbf x}\nu(M),
\]
with $\nu=\nu_R$, the minimal number of homogeneous generators, in
the graded case and $\nu=\nu_{\widehat R}$ in the local case.
\end{theorem}

The constant $\eta_{\mathbf x}$ is explicit. In the graded case,
$\eta_{\mathbf x}=c/(dhC_{\mathbf x}+c)$ with $c=(\sigma-d)/8$,
$h=h^1_*(\OO_Y)$ and $C_{\mathbf x}=1+a_0+a_1+a_0^2$, where
$a_0=\dim_{\CC}R/(\mathbf x)R$ and $a_1=a_0-d$; the constants in
the local case, for an arbitrary system of parameters, are given in
Section~\ref{ex:local-subsection}. No bound on generator
degrees or regularity is imposed in the graded assertion.

For each of the Hesse, dual Hesse and complete-quadrangle
arrangements, of $s=12$, $9$ and $6$ lines, embed $\PP^2$ as a plane
in $\PP^{s-1}$ by the defining linear forms. For $n\ge2$ let $X_n$ be
its preimage under the $n$-th power map, a Hirzebruch--Kummer cover of
the plane branched along the arrangement. Let $Y_n$ be the blow-up of
$X_n$ along its reduced singular locus, $A_n$ the pullback of
$\OO_{X_n}(1)$ and $E_n$ the reduced exceptional divisor
(Section~\ref{ex:families}). The \emph{standard}
polarizations are $4A_n-E_n$, $6A_n-E_n$ and $3A_n-E_n$,
respectively.

For the Hesse pairs $(Y_n,4A_n-E_n)$ the local criterion
\eqref{main:local-condition} holds for every $n\ge3$. It holds for the dual Hesse pairs $(Y_n,6A_n-E_n)$ for
$n\ge5$, $n\ne7$, and for the complete-quadrangle pairs
$(Y_n,3A_n-E_n)$ for $n\ge6$, $n\ne8$ (Corollary~\ref{ex:standard-cones}).
In the quadrangle family $K_{Y_n}=(n-2)(3A_n-E_n)$; for $n=6$
this is $K_Y=4H$, the example of \cite[Corollary~4.3]{An2},
\cite[Section~7]{Ch} and \cite[Section~3]{J}. The two
omitted standard polarizations satisfy $5<\lambda\le6$, a range in
which no integer $q$ satisfies \eqref{main:local-condition}
(Lemma~\ref{ex:range}). On those same surfaces, the different
polarizations $7A_7-E_7$ and $4A_8-E_8$ do satisfy it
(Corollary~\ref{ex:other-polarizations}). The two omitted standard cones remain
open for arbitrary modules (Remark~\ref{ex:omitted}).

\begin{theorem}[Dimension and characteristic]
\label{main:dimensions}
The following hold.
\begin{enumerate}[label=\textup{(\alph*)}]
\item For every $D\ge3$ there is a complete normal local domain of
dimension $D$ that has residue field $\CC$ and an isolated singularity
and admits no lim Cohen--Macaulay sequence. One may take the completed
section ring of
\[
 Y_{32m}\times\mathbf P^{m-2},\qquad
 \operatorname{pr}_1^*(6A_{32m}-E_{32m})+
 \operatorname{pr}_2^*\OO_{\mathbf P^{m-2}}(1),\qquad m=D-1,
\]
where $Y_{32m}$ is the dual Hesse surface $Y_n$ with $n=32m$.
\item For every prime $p$ and every $D\ge4$ there is a complete
normal local domain of mixed characteristic $(0,p)$ and dimension $D$,
with residue field $\overline{\mathbf F}_p$, that admits no lim
Cohen--Macaulay sequence.
\end{enumerate}
\end{theorem}

For $D=3$ one may also take the completed cones of
Corollaries~\ref{ex:standard-cones} and~\ref{ex:other-polarizations}
that satisfy \eqref{main:local-condition}. The complex examples use a
higher-dimensional replacement and restriction to surface sections.
The mixed-characteristic examples use an integral model and transfer
from the vertex of its generic fibre; their
singular locus is positive-dimensional (Remark~\ref{arith:scope} and
the paragraph after Proposition~\ref{arith:mixed-variables}).
Theorem~\ref{main:criteria} is proved in Section~\ref{ex:section}.
Theorem~\ref{main:dimensions}(a) is Corollary~\ref{hd:dual-hesse},
and part (b) is proved after
Proposition~\ref{arith:mixed-variables}.

\subsection{Strategy of proof}
Theorems~\ref{main:graded}--\ref{main:criteria} rest on two
cohomological inequalities and one transfer; Theorem~\ref{main:dimensions}
adds two further steps.
\begin{enumerate}[label=\textup{(\arabic*)}]
\item \emph{The graded inequality} (Section~\ref{surf:section}). Push
a torsion-free sheaf $E$ forward along the finite map
$f\colon Y\to\PP^2$ of degree $d$ given by three sections of
$\OO_Y(H)$. Riemann--Roch on $Y$, Bogomolov's inequality, the Hodge index
theorem and a comparison of Harder--Narasimhan polygons show that, when
$d<\sigma$, the Harder--Narasimhan factors of $f_*E$ are far from
equality in Bogomolov's inequality. An elementary estimate for $h^1_*$
on the plane turns this excess into $h^1_*(E)\ge(\sigma-d)\rk E/8$.
\item \emph{The local inequality} (Section~\ref{loc:section}).
Following Chen \cite{Ch}, extend a reflexive module over
$\widehat R$ to a reflexive sheaf on the regular threefold
$W\to\Spec\widehat R$ obtained from the total space of $\OO_Y(-H)$,
whose closed fibre is $Y$ (Proposition~\ref{loc:model}), and modify it
along $Y$ until the normalized Harder--Narasimhan slopes
of its restriction $V$ lie in the window
$[(\lambda-1)/2,(\lambda+1)/2)$. Then
$\chi(V)\le-(\sigma-d)\rk V/8$, and this negativity persists on the
punctured spectrum up to two tail terms
(Lemma~\ref{loc:persistence}). By Lemma~\ref{loc:tail-bound}, if
$K_Y\cdot H<5d$, then $(\lambda+1)/2<3$ and the tails are bounded by
$h^0(V)$ and $h^2(V)$; at $K_Y\cdot H=5d$ the second tail costs at
most $2\rk V$ more. Veronese
subrings give the parameter $q$.
\item \emph{The transfer} (Section~\ref{transfer:section}). A
dualizing construction replaces an arbitrary module $M$ by a reflexive
module $M'$ of at least the same rank, with
$\length \HH^2_{\mathfrak m}(M')\le hC_{\mathbf x}\kappa_{\mathbf x}(M)$.
Combined with step~(1) in the graded case or step~(2) in the local
case, this gives $\kappa_{\mathbf x}(M)\ge\eta_{\mathbf x}\nu(M)$.
\item \emph{Higher dimensions} (Section~\ref{hd:section}). The
replacement extends to every dimension, and restriction to a general
surface section reduces Theorem~\ref{main:dimensions}(a) to step~(2).
\item \emph{Mixed characteristic} (Section~\ref{arith:section}).
Complete an integral model of the section ring at the vertex of its
special fibre. Flat base change to the completed cone of its generic
fibre, and then to $\CC$, carries a hypothetical lim Cohen--Macaulay
sequence over this ring to the complex cone, and formal variables give every dimension at least
four.
\end{enumerate}

\subsection{The weak condition}
The ordinary criterion does not settle the weak condition. The
companion paper \cite{AnII} studies it on the same cones. It proves an
Euler-characteristic gap for modules of depth at least two and, on
completed cones satisfying \eqref{main:local-condition}, a bound on the
ranks of torsion-free modules $M$ for each fixed value of
$\chi_1(\mathbf x;M)$ smaller than $\rk_A\widehat R$, where
$A=\CC[[\mathbf x]]$. Under $d<\sigma$ alone, without the slope
condition, it also proves that every nonzero finite graded
torsion-free $R$-module $M$ satisfies
$\rk_RM\le C\,\chi_1(\mathbf x;M)^2$, where $\mathbf x$ consists of
three sections of $\OO_Y(H)$ without a common zero and the constant
$C$ is independent of $M$. The existence of weakly lim
Cohen--Macaulay sequences on these cones remains open.

\subsection{Related work and organization}\label{intro:related}
The asymptotic framework is due to Bhatt, Hochster and Ma
\cite{BHM,Ho}. For existence in mixed characteristic for rings admitting perfectoid
towers, and the characteristic-free form of the existence conjecture,
see \cite{IS}; for the consequence that systems of parameters are
$Q$-sequences in local rings with a weakly lim Cohen--Macaulay
sequence, and an extension of Lech's formula, see \cite{HS}. The surfaces, and the absence on them of
Ulrich bundles and of arithmetically Cohen--Macaulay bundles
(equivalently, of graded maximal Cohen--Macaulay modules over their
section rings), come from \cite{An1,An2}; the absence of small
Cohen--Macaulay modules over their completed cones is discussed in
\cite{An3}. For $2d<\sigma$ the absence of small Cohen--Macaulay
modules is \cite[Theorem~1.2]{OAI}; see \cite[Theorem~10.3]{An3} for
every dimension.

Assuming \cite[Main claim~1.1]{Ch},
the case $K_Y\equiv4H$, $q=1$, $M$ reflexive of
Theorem~\ref{main:local}(a), the Hexagon record \cite{J} (contributor
A.~Jiang) asserts that $\widehat R$ admits no lim Cohen--Macaulay
sequence for every pair with $K_Y\equiv4H$ and $d<\sigma$
\cite[Corollary~3.2]{J}. For $n=4k+2$ the standard quadrangle
polarization $H=3A_n-E_n$ satisfies $K_{Y_n}\equiv4kH$ and
$k^2H^2<\sigma$, so this assertion also covers the standard quadrangle cones with
$n\equiv2\pmod4$, through the Veronese subring of degree $k$ and
Remark~\ref{loc:restriction-scalars}. For no other cone of
Corollaries~\ref{ex:standard-cones} and~\ref{ex:other-polarizations}
is $K_Y\equiv4kH$ with $k$ a positive integer
(Remark~\ref{ex:hesse-nonsubcanonical}). The record also gives a replacement for complete
local domains of any dimension with finite-length lower local
cohomology \cite[Definition~2.1, Proposition~2.3 and
Lemmas~2.4--2.5]{J}, with constants depending on those of
\cite[Theorem~4.6]{BHM}; a uniform Koszul bound
\cite[Corollary~2.11]{J}; the reduction to reflexive modules over
three-dimensional normal local domains \cite[Remark~2.10]{J}; and a
mixed-characteristic transfer with residue field $\mathbf F_p$
\cite[Remark~3.3]{J}. The record states that its proof was produced
by an AI system and not checked by its human contributor. We claim no
priority for the statements of \cite{J} just listed; the versions
used here are proved in full below. Relative to \cite{Ch} and \cite{J}, the
following are new.
\begin{enumerate}[label=\textup{(\roman*)}]
\item An elementary proof, independent of \cite{AC}, of the graded
cohomology bound,
Theorem~\ref{main:graded}, with the Harder--Narasimhan correction terms
of Theorem~\ref{surf:refined-inequality}; the pushforward argument is
that of \cite[Lemma~5.5 and Proposition~9.1]{An3}. The numerical bound
itself also follows from the stronger estimate of
\cite[Corollary~8.3(a)]{AC} (Remark~\ref{surf:stronger-bundle-bound}).
\item The local inequality, Theorem~\ref{main:local}, without
proportionality, with the endpoint and the Veronese parameter; it
combines the method of \cite{Ch} with
\cite[Sections~6 and~7 and Remark~9.2]{An3}.
\item The replacement with explicit constants: $hC_{\mathbf x}$ in
dimension three (Proposition~\ref{transfer:replacement}), together
with its graded form (Appendix~\ref{app:graded}), and the explicit
Koszul gaps $\eta_{\mathbf x}$ of Theorem~\ref{main:criteria}.
\item The cones of Corollaries~\ref{ex:standard-cones}
and~\ref{ex:other-polarizations} other than the quadrangle cones with
$n\equiv2\pmod4$, which \cite{J} covers conditionally on \cite{Ch},
and the isolated examples in every dimension at least three,
Theorem~\ref{main:dimensions}(a).
\item The arithmetic transfer with algebraically closed residue field
and the mixed-characteristic examples in every dimension at least
four, Theorem~\ref{main:dimensions}(b).
\end{enumerate}

Section~\ref{surf:section} proves Theorem~\ref{main:graded} and
Section~\ref{loc:section} proves Theorem~\ref{main:local}.
Section~\ref{transfer:section} gives the replacement and the Koszul
gap; Section~\ref{ex:section} proves Theorem~\ref{main:criteria} and
treats the three families. Sections~\ref{hd:section}
and~\ref{arith:section} prove Theorem~\ref{main:dimensions}, and
Section~\ref{sec:questions} lists open questions.
Appendix~\ref{app:graded} proves the graded forms of the results of
Section~\ref{transfer:section}.

\section{The surface inequality and graded modules}\label{surf:section}

This section proves Theorem~\ref{main:graded}. We keep the notation
\(d\), \(\sigma\), \(\lambda\), \(\bar\mu\) and \(h^1_*\) of the
Introduction, and write \(E(t)=E\otimes\OO_Y(tH)\), so that
\(h^1_*(E)=\sum_{t\in\mathbf Z}h^1(Y,E(t))\in[0,\infty]\). Slopes on
\(Y\) are taken with respect to \(H\), and slopes on \(\mathbf P^2\)
with respect to a line. For a torsion-free sheaf \(G\) of positive
rank \(\rho\), put
\[
 \Delta(G)=2\rho c_2(G)-(\rho-1)c_1(G)^2,
 \qquad \delta(G)=\frac{\Delta(G)}{\rho^2}.
\]
We identify degree-four rational cohomology with \(\mathbf Q\).

\subsection{Riemann--Roch and the basic identity}
\label{surf:basic-subsection}

For a coherent sheaf \(F\) of rank \(N>0\) on \(\mathbf P^2\), set
\[
 b(F)=\frac{c_1(F)}N,
 \qquad
 v(F)=\frac{2\operatorname{ch}_2(F)}N-b(F)^2.
\]
Here \(c_1(F)\in\mathbf Z\) is the degree of \(F\). Riemann--Roch takes the form
\begin{equation}\label{surf:plane-rr}
 \frac{\chi(F(t))}{N}
 =\frac12\left(t+b(F)+\frac32\right)^2
   +\frac{v(F)}2-\frac18
 \qquad(t\in\mathbf Z).
\end{equation}
Since \(2\operatorname{ch}_2=c_1^2-2c_2\), for torsion-free \(F\)
one has
\begin{equation}\label{surf:variance-discriminant}
 v(F)=-\delta(F).
\end{equation}

Let
\(0=F_0\subset F_1\subset\cdots\subset F_k=F\)
be the Harder--Narasimhan filtration of a torsion-free sheaf \(F\).
Write \(Q_j=F_j/F_{j-1}\), \(n_j=\rk Q_j\) and
\(\mu_j=b(Q_j)\). Define its slope variance by
\[
 v_{\mathrm{HN}}(F)
 =\frac1N\sum_{j=1}^kn_j\bigl(\mu_j-b(F)\bigr)^2
 =\frac1N\sum_{j=1}^kn_j\mu_j^2-b(F)^2.
\]
Additivity of rank, \(c_1\), and \(\operatorname{ch}_2\) gives
\[
 \frac{2\operatorname{ch}_2(F)}N
 =\frac1N\sum_jn_j\bigl(v(Q_j)+\mu_j^2\bigr).
\]
Consequently
\begin{equation}\label{surf:hn-discriminants}
 v_{\mathrm{HN}}(F)-v(F)
 =\frac1N\sum_{j=1}^kn_j\delta(Q_j)\ge0,
\end{equation}
where the inequality is Bogomolov's inequality~\cite[Theorem 3.4.1]{HL}.

Since \(|H|\) is base-point free, three general sections
\(x_0,x_1,x_2\) of \(\OO_Y(H)\) have no common zero. They define
\[
 f:Y\longrightarrow\mathbf P^2,
 \qquad f^*\OO_{\mathbf P^2}(1)=\OO_Y(H).
\]
A positive-dimensional fibre would contain a curve on which \(H\)
has degree zero, contrary to ampleness. Thus \(f\) is finite and
surjective, of degree \(d\). Its direct-image functor is exact, and
the projection formula gives
\begin{equation}\label{surf:finite-cohomology}
 \HH^i(Y,G(t))=\HH^i\bigl(\mathbf P^2,(f_*G)(t)\bigr)
 \qquad(i\ge0,\ t\in\mathbf Z)
\end{equation}
for every coherent sheaf \(G\). If \(G\) is torsion-free of rank
\(\rho\), then \(f_*G\) is torsion-free of rank \(d\rho\).

For such a sheaf \(G\) on \(Y\), put
\[
 \xi=\frac{c_1(G)}\rho,\qquad
 \beta=\xi-\frac12K_Y,\qquad
 \alpha(G)=\beta-\frac{\beta\cdot H}{d}H
 \quad\text{in }\operatorname{NS}(Y)_{\mathbf Q}.
\]
Then \(\alpha(G)\cdot H=0\), so the Hodge index theorem gives
\(\alpha(G)^2\le0\), and
\begin{equation}\label{surf:orthogonal-part}
 \beta\cdot H=d\left(\bar\mu(G)-\frac\lambda2\right),
 \qquad
 \beta^2=\alpha(G)^2
       +d\left(\bar\mu(G)-\frac\lambda2\right)^2.
\end{equation}
Riemann--Roch on \(Y\) gives
\begin{equation}\label{surf:surface-rr}
 \frac{\chi(G(t))}\rho
 =\frac d2t^2+(\beta\cdot H)t
   +\frac12\beta^2-\frac\sigma8-\frac{\delta(G)}2.
\end{equation}
Indeed,
\(\operatorname{ch}_2(G)/\rho=(\xi^2-\delta(G))/2\), and
\[
 \frac12\xi\cdot(\xi-K_Y)+\chi(\OO_Y)
 =\frac12\beta^2-\frac\sigma8.
\]

\begin{proposition}[The pushforward identity]
\label{surf:pushforward-identity}
Let \(E\) be a torsion-free sheaf of rank \(r>0\) on \(Y\), with
Harder--Narasimhan factors \(G_1,\ldots,G_s\) of ranks \(r_i\).
Put
\[
 w_i=\frac{r_i}r,\qquad
 b_i=b(f_*G_i),\qquad
 \overline b=\sum_iw_ib_i,\qquad
 v_b(E)=\sum_iw_i(b_i-\overline b)^2.
\]
Then
\begin{equation}\label{surf:pushforward-slopes}
 b_i=\bar\mu(G_i)-\frac\lambda2-\frac32,
 \qquad b_1>\cdots>b_s,
\end{equation}
and
\begin{equation}\label{surf:pushforward-variance}
 v(f_*E)=\sum_iw_i v(f_*G_i)+v_b(E),
 \qquad
 v(f_*G_i)=\frac14-\frac\sigma{4d}
          -\frac{\delta(G_i)-\alpha(G_i)^2}{d}.
\end{equation}
Moreover, \(\delta(G_i)\ge0\) and \(\alpha(G_i)^2\le0\).
\end{proposition}
\begin{proof}
We recall the argument of~\cite[Lemma 5.5 and Proposition 9.1]{An3}.
The filtration of \(f_*E\) induced by that of \(E\) has quotients
\(f_*G_i\). Additivity of rank, \(c_1\), and
\(\operatorname{ch}_2\) gives the first identity in
\eqref{surf:pushforward-variance}, just as in
\eqref{surf:hn-discriminants}.

For an arbitrary torsion-free \(G\) of rank \(\rho\), compare
\eqref{surf:surface-rr} with \eqref{surf:plane-rr}, using
\(\chi(G(t))=\chi((f_*G)(t))\) and \(\rk(f_*G)=d\rho\).
The coefficients of \(t\) give
\[
 b(f_*G)+\frac32=\frac{\beta\cdot H}{d},
\]
which proves \eqref{surf:pushforward-slopes}. The constant terms give
\[
 \frac{v(f_*G)}2-\frac18+\frac{(\beta\cdot H)^2}{2d^2}
 =\frac1d\left(\frac12\beta^2-\frac\sigma8
                         -\frac{\delta(G)}2\right).
\]
Substitute
\(\beta^2=\alpha(G)^2+(\beta\cdot H)^2/d\) to obtain the second
identity in \eqref{surf:pushforward-variance}. Bogomolov's inequality
applies to each semistable \(G_i\); the Hodge index theorem gives
the assertion about \(\alpha(G_i)\). Finally, the normalized slopes
of the \(G_i\) decrease strictly, so the \(b_i\) do as well.
\end{proof}

\subsection{Polygons and intermediate cohomology on the plane}
\label{surf:plane-subsection}

For a torsion-free sheaf \(F\) of rank \(N\) on \(\mathbf P^2\),
with Harder--Narasimhan data \((n_j,\mu_j)\) as above, its
Harder--Narasimhan polygon is the continuous function
\(p_F:[0,N]\to\mathbf R\) with \(p_F(0)=0\), affine of slope
\(\mu_j\) on
\([n_1+\cdots+n_{j-1},n_1+\cdots+n_j]\). It is concave and satisfies
\(p_F(N)=\deg F\).

\begin{lemma}[Subsheaves and the polygon]\label{surf:polygon}
For every subsheaf \(F'\subset F\) of rank \(\rho\),
\(\deg F'\le p_F(\rho)\).
\end{lemma}
\begin{proof}
This is the usual polygon comparison; cf.\ \cite{Sha}.
Set \(F'_j=F'\cap F_j\). The quotient \(F'_j/F'_{j-1}\) is a
subsheaf of \(Q_j\), of some rank \(0\le\rho_j\le n_j\). If
\(\rho_j=0\), it is zero because \(Q_j\) is torsion-free; otherwise
its degree is at most \(\mu_j\rho_j\), by semistability. Hence
\[
 \deg F'\le\sum_j\mu_j\rho_j,
 \qquad 0\le\rho_j\le n_j,
 \quad\sum_j\rho_j=\rho.
\]
Since \(\mu_1>\cdots>\mu_k\), the maximum of
\(\sum_j\mu_j y_j\) under these constraints is attained by filling
\(y_1,y_2,\ldots\) in order. Its value is \(p_F(\rho)\).
\end{proof}

\begin{lemma}[Comparison of slope variances]
\label{surf:variance-comparison}
Suppose
\(0=F^0\subset F^1\subset\cdots\subset F^s=F\)
is a filtration by subsheaves whose successive quotients have
positive ranks \(m_i\) and slopes \(b_1\ge\cdots\ge b_s\).
Then
\[
 \sum_{i=1}^s\frac{m_i}{N}\bigl(b_i-b(F)\bigr)^2
 \le v_{\mathrm{HN}}(F).
\]
\end{lemma}
\begin{proof}
Of the filtration, only the ordering of the slopes \(b_i\) is used;
its quotients need not be semistable.
Put \(\rho_i=\rk F^i\), and let \(\psi:[0,N]\to\mathbf R\) be the
continuous function with \(\psi(0)=0\) that is affine of slope
\(b_i\) on \([\rho_{i-1},\rho_i]\).
Then \(\psi(\rho_i)=\deg F^i\) and
\(\psi(N)=\deg F=p_F(N)\). By Lemma~\ref{surf:polygon},
\(\psi(\rho_i)\le p_F(\rho_i)\). On each interval
\([\rho_{i-1},\rho_i]\), concavity places \(p_F\) above its chord,
which in turn lies above \(\psi\). Thus
\(\varphi=p_F-\psi\ge0\), with \(\varphi(0)=\varphi(N)=0\).

Both functions are piecewise affine, and
\[
 \int_0^N(p_F')^2-\int_0^N(\psi')^2
 =\int_0^N(p_F'-\psi')^2+2\int_0^N\psi'\varphi'.
\]
Abel summation gives
\[
 \int_0^N\psi'\varphi'
 =\sum_{i=1}^s b_i\bigl(\varphi(\rho_i)-\varphi(\rho_{i-1})\bigr)
 =\sum_{i=1}^{s-1}(b_i-b_{i+1})\varphi(\rho_i)\ge0.
\]
It follows that
\[
 \sum_jn_j\mu_j^2=\int_0^N(p_F')^2
 \ge\int_0^N(\psi')^2=\sum_im_ib_i^2.
\]
Divide by \(N\) and subtract the square of the common mean
\(b(F)\). For a split bundle \(F=\bigoplus_j\OO(a_j)\), this
comparison is~\cite[Lemma 2.4]{An2}.
\end{proof}

\begin{proposition}[Cohomology in slope windows]
\label{surf:plane-h1}
Let \(F\) be torsion-free of rank \(N\) on \(\mathbf P^2\), with
Harder--Narasimhan factors \(Q_j\) of ranks \(n_j\) and slopes
\(\mu_j\). For \(t\in\mathbf Z\), put
\(J_t=\{j:-t-2\le\mu_j\le-t-1\}\). Then
\begin{equation}\label{surf:window-h1}
 \begin{split}
 h^1(\mathbf P^2,F(t))
 &\ge\sum_{j\in J_t}n_j
 \left(\frac{\delta(Q_j)}2+\frac18
       -\frac12\left(t+\mu_j+\frac32\right)^2\right)\\
 &\ge\frac12\sum_{j\in J_t}n_j\delta(Q_j).
 \end{split}
\end{equation}
Consequently
\begin{equation}\label{surf:total-h1}
 \sum_{t\in\mathbf Z}h^1(\mathbf P^2,F(t))
 \ge\frac12\sum_jn_j\delta(Q_j)
 =\frac N2\bigl(v_{\mathrm{HN}}(F)-v(F)\bigr).
\end{equation}
Equality in \eqref{surf:total-h1} is possible only when both sides
are zero.
\end{proposition}
\begin{proof}
Fix \(t\). If \(J_t\) is nonempty, its indices are consecutive,
say \(j_-,\ldots,j_+\). Put
\[
 F'=F_{j_--1},\qquad
 F''=F/F_{j_+},\qquad
 F^{\mathrm{mid}}=F_{j_+}/F_{j_--1}.
\]
These sheaves are torsion-free or zero, with the corresponding
segments of the Harder--Narasimhan filtration.

First, \(\HH^2(F'(t))=0\). If \(F'\ne0\), all its Harder--Narasimhan slopes
after twisting exceed \(-1\), since its factors lie outside \(J_t\)
on the side of larger slopes. A nonzero homomorphism
\(F'(t)\to\OO(-3)\) would have image a torsion-free rank-one
quotient of degree at most \(-3\), contradicting the lower bound on
the slopes of torsion-free quotients of \(F'(t)\). Serre duality
therefore gives the asserted vanishing. Similarly,
\(\HH^0(F''(t))=0\): its maximal slope after twisting is less than
\(-2\), so it cannot contain a copy of \(\OO\).

The inclusion \(F'\subset F\) gives
\[
 h^0(F(t))\ge h^0(F'(t))
 =\chi(F'(t))+h^1(F'(t))\ge\chi(F'(t)).
\]
The map \(\HH^2(F(t))\to \HH^2(F''(t))\) is surjective, since
third cohomology vanishes on a surface, and hence
\[
 h^2(F(t))\ge h^2(F''(t))
 =\chi(F''(t))+h^1(F''(t))\ge\chi(F''(t)).
\]
Additivity of Euler characteristic now yields
\[
 h^1(F(t))\ge-\chi(F^{\mathrm{mid}}(t))
 =-\sum_{j\in J_t}\chi(Q_j(t)).
\]
By \eqref{surf:plane-rr} and \eqref{surf:variance-discriminant},
\[
 -\chi(Q_j(t))
 =n_j\left(\frac{\delta(Q_j)}2+\frac18
       -\frac12\left(t+\mu_j+\frac32\right)^2\right).
\]
For \(j\in J_t\), the expression inside the last square has
absolute value at most \(1/2\). This proves
\eqref{surf:window-h1}; when \(J_t\) is empty, it is immediate.

Every index \(j\) belongs to at least one \(J_t\), since
\([-\mu_j-2,-\mu_j-1]\) is a closed interval of length one.
All summands in \eqref{surf:window-h1} are nonnegative by
Bogomolov's inequality. Summing over \(t\) proves the inequality in
\eqref{surf:total-h1}, and its last equality is
\eqref{surf:hn-discriminants}.

To analyze equality, retain the discarded terms. Let
\(\theta_j=\mu_j-\lfloor\mu_j\rfloor\). If \(\mu_j\notin\mathbf Z\),
then \(j\) occurs only for \(t=-\lfloor\mu_j\rfloor-2\), and
\[
 \frac18-\frac12\left(t+\mu_j+\frac32\right)^2
 =\frac12\theta_j(1-\theta_j).
\]
If \(\mu_j\in\mathbf Z\), it occurs for both
\(t=-\mu_j-2\) and \(t=-\mu_j-1\), each time contributing
\(n_j\delta(Q_j)/2\). Thus
\begin{equation}\label{surf:fractional-slopes}
 \sum_t h^1(\mathbf P^2,F(t))
 \ge\frac12\sum_j n_j\bigl(\delta(Q_j)+\theta_j(1-\theta_j)\bigr)
    +\frac12\sum_{j:\,\mu_j\in\mathbf Z}n_j\delta(Q_j).
\end{equation}
Equality in \eqref{surf:total-h1} forces every \(\theta_j\) to be
zero. Then the last sum in \eqref{surf:fractional-slopes} forces
\(\sum_jn_j\delta(Q_j)=0\), proving the assertion.
\end{proof}

\subsection{The surface inequality, including its correction terms}
\label{surf:inequality-subsection}

\begin{theorem}[The refined surface inequality]
\label{surf:refined-inequality}
Let \(E\) be a torsion-free sheaf of rank \(r>0\) on \(Y\), with
Harder--Narasimhan factors \(G_i\) of ranks \(r_i\). Then
\begin{equation}\label{surf:refined-bound}
 h^1_*(E)\ge
 \frac{\sigma-d}{8}r
 +\frac12\sum_i r_i\bigl(\delta(G_i)-\alpha(G_i)^2\bigr)
 +\frac{dr}{2}\bigl(v_{\mathrm{HN}}(f_*E)-v_b(E)\bigr).
\end{equation}
The two correction terms are nonnegative. Equality is possible
only when both sides of \eqref{surf:refined-bound} are zero.
\end{theorem}
\begin{proof}
Set \(F=f_*E\), of rank \(N=dr\). By
\eqref{surf:finite-cohomology} and Proposition~\ref{surf:plane-h1},
\begin{equation}\label{surf:pushforward-h1}
 h^1_*(E)\ge\frac N2\bigl(v_{\mathrm{HN}}(F)-v(F)\bigr),
\end{equation}
with equality only when both sides are zero. Proposition
\ref{surf:pushforward-identity} gives
\[
 v(F)=\frac14-\frac\sigma{4d}
      -\frac1d\sum_iw_i\bigl(\delta(G_i)-\alpha(G_i)^2\bigr)
      +v_b(E).
\]
Therefore
\begin{equation}\label{surf:variance-decomposition}
 v_{\mathrm{HN}}(F)-v(F)
 =\frac{\sigma-d}{4d}
  +\frac1d\sum_iw_i\bigl(\delta(G_i)-\alpha(G_i)^2\bigr)
  +v_{\mathrm{HN}}(F)-v_b(E).
\end{equation}
Multiplication by \(N/2=dr/2\) identifies its right-hand side
with that of \eqref{surf:refined-bound}. The first correction term
is nonnegative by Proposition~\ref{surf:pushforward-identity}.

Finally, the images of the Harder--Narasimhan filtration of \(E\)
under the exact functor \(f_*\) filter \(F\) with quotients
\(f_*G_i\), of ranks \(dr_i\) and strictly decreasing slopes
\(b_i\). Their weighted mean is \(b(F)\), so Lemma
\ref{surf:variance-comparison} gives
\(v_b(E)\le v_{\mathrm{HN}}(F)\). This proves all assertions.
\end{proof}

\begin{corollary}[The elementary surface bound]
\label{surf:coherent-inequality}
Let \(E\) be a coherent sheaf of rank \(r>0\) on \(Y\), with torsion
subsheaf \(T\).
\begin{enumerate}[label=\textup{(\alph*)}]
\item One has
\begin{equation}\label{surf:elementary-bound}
 h^1_*(E)\ge\frac{\sigma-d}{8}r,
\end{equation}
and this inequality is strict if \(d<\sigma\).
\item If \(h^1_*(E)<\infty\), then \(T\) has finite length, \(E/T\) is
locally free, and \(h^1_*(E)=h^1_*(E/T)\).
\end{enumerate}
\end{corollary}
\begin{proof}
Set \(E'=E/T\), a torsion-free sheaf of rank \(r\). Since the support
of \(T\) has dimension at most one, \(\HH^2(Y,T(t))=0\), and
\(\HH^1(Y,E(t))\to \HH^1(Y,E'(t))\) is surjective for every \(t\).

(a) Thus \(h^1_*(E)\ge h^1_*(E')\), and
Theorem~\ref{surf:refined-inequality}, after dropping its two
nonnegative correction terms, gives \eqref{surf:elementary-bound}.
If \(d<\sigma\), the retained bound is positive. Equality in
\eqref{surf:elementary-bound} would force equality in
Theorem~\ref{surf:refined-inequality} for \(E'\), with both sides
positive, which the equality assertion of that theorem excludes.

(b) Suppose first that \(T\) has one-dimensional support. For
\(t\ll0\) one has \(\HH^0(E'(t))=0\), so \(\HH^1(T(t))\) injects into
\(\HH^1(E(t))\). Moreover
\(h^1(T(t))\ge-\chi(T(t))=-(c_1(T)\cdot H)t-\chi(T)\) with
\(c_1(T)\cdot H>0\), so \(h^1_*(E)=\infty\). If \(T\) has finite
length, then \(\HH^1(E(t))\simeq \HH^1(E'(t))\) for every \(t\). If, in
addition, \(E'\) is not locally free, then the finite-length sheaf
\(E'^{\vee\vee}/E'\) is nonzero, and \(\HH^0(E'^{\vee\vee}/E')\)
injects into \(\HH^1(E'(t))\) for \(t\ll0\), so again
\(h^1_*(E)=\infty\). Hence \(h^1_*(E)<\infty\) forces \(T\) to have
finite length and \(E'\) to be locally free, and then
\(h^1_*(E)=h^1_*(E')\).
\end{proof}

By (a) and (b), a nonzero torsion-free sheaf \(E\) with
\(h^1_*(E)=0\) is locally free and satisfies \(\sigma\le d\); this
recovers~\cite[Theorem 1.1]{An2}. If also \(\sigma=d\), both
correction terms in \eqref{surf:refined-bound} vanish for \(E\).

The correction terms carry information even when \(\sigma<d\). On
\(Y=\mathbf P^1\times\mathbf P^1\) with \(H\) of type \((1,1)\), one
has \(d=2\) and \(\sigma=0\). For \(E=\OO(a,b)\) and \(g=|a-b|\), one
has \(\delta(E)=0\) and \(\alpha(E)^2=-g^2/2\), so
\eqref{surf:refined-bound} without its last term gives
\(h^1_*(E)\ge(g^2-1)/4\), whereas the K\"unneth formula gives
\(h^1_*(E)=(g^3-g)/6\).

\subsection{Graded modules and Veronese amplification}
\label{surf:graded-subsection}

Let \(R=R(Y,H)\) and \(\mm=R_+\) be as in the Introduction, and put
\(S=\CC[x_0,x_1,x_2]\). The sections defining \(f\) induce a map
\(S\to R\). It is injective because \(f\) is dominant. It is finite
because \(R\) is the module of sections in nonnegative degrees of the
coherent sheaf \(f_*\OO_Y\) on \(\mathbf P^2\); such a module is
finitely generated over \(S\), since its graded pieces are
finite-dimensional and, by Castelnuovo--Mumford regularity,
multiplication by \(S_1\) is surjective in large degrees. Hence
\(R/S_+R\) is finite-dimensional and \(\sqrt{S_+R}=\mm\). We have
\(\Proj R=Y\) and \(\widetilde{R(t)}=\OO_Y(tH)\), without assuming
that \(R\) is generated in degree one~\cite[Lemma 1.2]{An2}: the open
sets \(D_+(x_i)\) cover \(\Proj R\) because \(\sqrt{S_+R}=\mm\), each
of them is the affine open set \(\{x_i\ne0\}\) of \(Y\), and
\(R(t)_{(x_i)}=x_i^tR_{(x_i)}\) is the module of sections of
\(\OO_Y(tH)\) over it.

If \(M\) is a finitely generated graded \(R\)-module and
\(E=\widetilde M\) on \(Y\), its associated sheaf as an
\(S\)-module is \(f_*E\). The local-cohomology comparison gives
\begin{equation}\label{surf:graded-cohomology}
 \HH^2_{\mm}(M)=\HH^2_{S_+}(M)
 \simeq\bigoplus_{t\in\mathbf Z}
             \HH^1\bigl(\mathbf P^2,(f_*E)(t)\bigr)
 =\bigoplus_{t\in\mathbf Z}\HH^1(Y,E(t)).
\end{equation}
The first equality uses independence of local cohomology from the
base ring and from the radical of the support ideal; the middle
isomorphism is the usual graded local-cohomology correspondence~\cite[Theorem A4.1]{Eis}; the last equality follows from finiteness
of \(f\). These identifications are natural in \(M\); hence they are
\(R\)-linear when \(R_s\) acts on the last direct sum by
multiplication with sections of \(\OO_Y(sH)\).

\begin{corollary}[The graded inequality]\label{surf:graded-inequality}
For every finitely generated graded \(R\)-module \(M\),
\[
 \dim_{\CC}\HH^2_{\mm}(M)\ge\frac{\sigma-d}{8}\rk_R M.
\]
\end{corollary}
\begin{proof}
The rank of \(f_*E\) is \(\rk_S M=d\rk_R M\), so
\(\rk E=\rk_R M\). If this rank is zero, the assertion is
immediate. Otherwise it follows from
Corollary~\ref{surf:coherent-inequality} and
\eqref{surf:graded-cohomology}.
\end{proof}

Conversely, every coherent sheaf \(E\) on \(Y\) is the sheaf
associated to a finitely generated graded \(R\)-module, for instance
\(M=\bigoplus_{t\ge0}\HH^0(Y,E(t))\), which is finitely generated for the
same reason as \(R\).
Equation~\eqref{surf:graded-cohomology} then identifies
\(\dim_{\CC}\HH^2_{\mm}(M)\) with \(h^1_*(E)\). Thus the coherent
and graded formulations are equivalent.

\begin{corollary}[Veronese amplification]\label{surf:veronese}
For every coherent sheaf \(E\) of rank \(r>0\) on \(Y\) and every
integer \(q\ge1\),
\[
 h^1_*(E)\ge\frac{q(\sigma-q^2d)}8r,
\]
and the inequality is strict if \(q^2d<\sigma\). Likewise, every
finitely generated graded \(R\)-module \(M\) satisfies
\(\dim_{\CC}\HH^2_{\mm}(M)\ge q(\sigma-q^2d)\rk_RM/8\).
\end{corollary}
\begin{proof}
The divisor \(qH\) is ample and globally generated, of square
\(q^2d\). Apply Corollary~\ref{surf:coherent-inequality}(a) for the
polarization \(qH\) to
\[
 E_q=\bigoplus_{i=0}^{q-1}E(iH),\qquad \rk E_q=qr.
\]
Every integer has a unique expression \(qt+i\),
\(0\le i<q\), and therefore
\[
 \sum_{t\in\mathbf Z}h^1(Y,E_q(tqH))
 =\sum_{t\in\mathbf Z}\sum_{i=0}^{q-1}
          h^1(Y,E((qt+i)H))
  =h^1_*(E).
\]
The resulting inequality is the first claim, including strictness.
For a graded module \(M\) of positive rank, apply it to
\(\widetilde M\), which has rank \(\rk_RM\) (proof of
Corollary~\ref{surf:graded-inequality}), and use
\eqref{surf:graded-cohomology}; for rank zero there is nothing to
prove.
\end{proof}

\begin{remark}[Comparison with the stronger bundle bound]
\label{surf:stronger-bundle-bound}
When \(\sigma\ge3d\), the choice
\(q=\lfloor\sqrt{\sigma/(3d)}\rfloor\) makes the Veronese bound of
order \(\sigma^{3/2}d^{-1/2}r\). For vector bundles, and hence, by
Corollary~\ref{surf:coherent-inequality}(b), for every coherent sheaf
of positive rank, \cite[Corollary 8.3(a)]{AC}
gives the stronger inequality
\[
 h^1_*(E)\ge rd\,L(\sigma/d)\qquad(\sigma\ge d),
\]
where \(L\) is the piecewise linear function, linear between
successive squares, with
\[
 L(g^2)=\frac{g^3-g}{12}\qquad(g\in\mathbf Z,\ g\ge1).
\]
The invariant denoted by \(\tau\) in that reference is our
\(\sigma\). In particular,
\(h^1_*(E)\ge(\sigma-d)r/6\). Applied to
\(\bigoplus_{i=0}^{q-1}E(iH)\) with the polarization \(qH\), as in the
proof of Corollary~\ref{surf:veronese}, the same result gives
\[
 h^1_*(E)\ge q^3d\,L\bigl(\sigma/(q^2d)\bigr)r
 \ge\frac{q(\sigma-q^2d)}6r
 \qquad(q\ge1,\ q^2d\le\sigma).
\]
We shall not use the stronger estimate. For graded modules, the graded constant
\((\sigma-d)/8\) used later may be replaced by \(dL(\sigma/d)\)
(Remark~\ref{ex:stronger-bundle-constant}); the local constants are
not affected.
\end{remark}

\section{The local inequality}\label{loc:section}

Let $Y$, $H$, $d$, $\sigma$, $\lambda$, $R=R(Y,H)$, $\mm=R_+$ and the
normalized slope $\bar\mu$ be as in Section~\ref{surf:section}, and
write $\widehat R$ for the completion of $R_{\mm}$. For a torsion-free
sheaf $V$ on $Y$, we write
$\bar\mu_{\max}(V)$ and $\bar\mu_{\min}(V)$ for the normalized slopes
of the first and last Harder--Narasimhan factors of $V$; every nonzero
subsheaf of $V$ has slope at most $\bar\mu_{\max}(V)$, and every
torsion-free quotient of positive rank has slope at least
$\bar\mu_{\min}(V)$. For the polarization $H$ itself we follow
\cite{Ch} and \cite[Sections~6 and~7 and Remark~9.2]{An3}; passing to
Veronese subrings then gives the following form of the local
inequality, which is Theorem~\ref{main:local}.

\begin{theorem}[Local cohomology inequality]\label{loc:inequality}
Let $q\ge1$ be an integer with $K_Y\cdot H\le5qd$, and put
\[
 \varepsilon_q=
 \begin{cases}
 0,&K_Y\cdot H<5qd,\\
 1,&K_Y\cdot H=5qd.
 \end{cases}
\]
Every finitely generated module $M$ over $\widehat R$, or over
$R_{\mm}$, satisfies
\begin{equation}\label{loc:main-bound}
 \length\HH^2_{\mm}(M)
 \ge q\left(\frac{\sigma-q^2d}{8}-2\varepsilon_q\right)\rk M,
\end{equation}
where $\rk M$ is the rank over $\widehat R$, respectively $R_{\mm}$.
The left side is allowed to be infinite. It is finite when $M$ is
reflexive.
\end{theorem}

Section~\ref{loc:model-subsection} constructs a regular model $W$ of
$\Spec\widehat R$ with closed fibre $Y$ and extends reflexive modules
to reflexive sheaves on $W$. Section~\ref{loc:slopes} modifies such an
extension along $Y$ until the Harder--Narasimhan slopes of its
restriction to $Y$ lie in a window of length one,
Section~\ref{loc:tails} bounds the two tails, and
Section~\ref{loc:persistence-subsection} shows that the resulting
negative Euler characteristic persists on the punctured spectrum.
Section~\ref{loc:qone} combines these steps for $q=1$, and
Section~\ref{loc:veronese-subsection} passes to Veronese subrings.

\subsection{The regular model}\label{loc:model-subsection}

\begin{proposition}[{The contraction model; cf.\ \cite[Proposition~2.1]{Ch},
\cite[Section~2]{OAI}, \cite[Proposition~6.1]{An3}}]\label{loc:model}
The ring $\widehat R$ is a complete normal local domain of dimension
three, with residue field $\CC$, regular outside its closed point. There is a proper morphism
$\pi:W\to\Spec\widehat R$ from a regular integral three-dimensional
scheme, obtained by base change from the contraction of the zero
section of the total space of $\OO_Y(-H)$, with the following
properties:
\begin{enumerate}
\item The closed fibre is the zero section $Y\subset W$, an effective
Cartier divisor, and $\OO_W(-Y)|_Y\simeq\OO_Y(H)$.
\item The morphism restricts to an isomorphism
$U:=W\setminus Y\simeq\Spec\widehat R\setminus\{\mm\}$.
\item The sheaf $\OO_W(-Y)$ is $\pi$-ample.
\end{enumerate}
Moreover, three sections $x_0,x_1,x_2$ of $\OO_Y(H)$ without a common
zero define a finite local inclusion
\[
 A=\CC[[x_0,x_1,x_2]]\subset\widehat R,
 \qquad \rk_A\widehat R=d,
\]
with the same residue field.
\end{proposition}
\begin{proof}
We include the proof, since $R$ is not assumed to be generated in
degree one. The three sections define a finite morphism $f:Y\to\mathbf P^2$
and a graded inclusion $S=\CC[x_0,x_1,x_2]\to R$. The $S$-module
$R$ consists of the sections in nonnegative degrees of the coherent
sheaf $f_*\OO_Y$, so it is finite of rank $d$, and
$(x_0,x_1,x_2)R$ is $\mm$-primary.

Let $W_0$ be the total space of $\OO_Y(-H)$. Since
$\OO_Y(-H)=f^*\OO_{\mathbf P^2}(-1)$, its morphism to the total space
of $\OO_{\mathbf P^2}(-1)$ induced by $f$ is the base change of the
finite morphism $f$, hence finite.
The latter space is the blow-up of $\mathbf A^3$ at the origin;
therefore $W_0$ is proper over $\Spec S$. The morphism
$\pi_0:W_0\to\Spec R$ defined by
\[
 \Gamma(W_0,\OO_{W_0})
 =\bigoplus_{t\ge0}\HH^0(Y,\OO_Y(tH))=R
\]
is proper as well: it is a morphism over $\Spec S$ to the separated
$S$-scheme $\Spec R$.

For a nonzero $x\in R_1$, the open set $Y_x=\{x\ne0\}$ is affine because a
multiple of $H$ is very ample, and
\[
 R[x^{-1}]=\Gamma(Y_x,\OO_Y)[x,x^{-1}].
\]
These charts for $x=x_0,x_1,x_2$ cover the complement of the vertex
and identify it with the complement of the zero section in $W_0$.
In particular that complement is regular. Negative twists of
$\OO_Y$ by $H$ have no sections, so the ring of global functions on
$W_0\setminus Y$ is
$\bigoplus_{t\in\mathbf Z}\HH^0(Y,\OO_Y(tH))=R$. Global functions
on an integral normal scheme form an integrally closed ring; hence
$R$ is normal. It is finitely generated because it is finite over $S$.

Locally on $Y$, an element of $R_t$ is the $t$th power of the fibre
coordinate multiplied by a local function. The sections $x_0,x_1,x_2$
are the fibre coordinate multiplied by functions generating the unit
ideal. Thus $\mm\OO_{W_0}=\OO_{W_0}(-Y)$: the closed fibre is the
zero section scheme-theoretically, with conormal bundle $\OO_Y(H)$.
The line bundle $\OO_{W_0}(-Y)$ is the pullback of the relatively
ample sheaf of the blow-up and is relatively ample over $\Spec R$.

Now set $W=W_0\times_{\Spec R}\Spec\widehat R$. The local ring
$R_{\mm}$ is excellent, so its completion homomorphism is regular and
$\widehat R$ is normal \cite[Tags 0AH2 and 0C23]{St}. Composed with
the localization $R\to R_{\mm}$, it gives a regular homomorphism
$R\to\widehat R$. Applying \cite[Tag 07C1]{St} on an affine open cover
of the finite-type $R$-scheme $W_0$, we see that the projection
$W\to W_0$ is regular, hence flat with regular fibres; since $W_0$ is
regular, so is $W$ \cite[Tag 031E]{St}. Properness, the Cartier fibre, its conormal bundle, relative
ampleness and the punctured-spectrum isomorphism are preserved.
The latter isomorphism shows that $\widehat R$ is regular outside
its closed point. Every connected component of $W$ is proper over
the local scheme $\Spec\widehat R$ and therefore meets the connected
closed fibre $Y$. Hence $W$ is connected; a connected regular scheme
is integral. Every closed point of $W$ lies on $Y$, because $\pi$ is
proper over the local scheme $\Spec\widehat R$; as $Y$ is a Cartier
divisor of dimension two, $W$ has dimension three. Finally
$\widehat R=R\otimes_S A$ is finite of rank $d$ over $A$, and the
inclusion is injective by flatness of $S\to A$.
\end{proof}

Let $\iota:Y\hookrightarrow W$ and $j:U\hookrightarrow W$ be the
inclusions. The morphism $j$ is affine because $Y$ is Cartier. For a
coherent sheaf $\mathcal F$ put
$\mathcal F(mY)=\mathcal F\otimes\OO_W(mY)$. If $\mathcal F$ is
torsion-free, a local equation of $Y$ is regular on it. With
$V=\mathcal F|_Y=\iota^*\mathcal F$, this gives exact sequences
\begin{equation}\label{loc:restriction-sequence}
 0\longrightarrow\mathcal F(-(m+1)Y)
 \longrightarrow\mathcal F(-mY)
 \longrightarrow\iota_*(V(mH))\longrightarrow0
 \qquad(m\in\mathbf Z).
\end{equation}
In particular $\mathcal F(mY)|_Y=V(-mH)$.

\begin{lemma}[Extension and restriction]\label{loc:extension}
The following hold.
\begin{enumerate}
\item For every finitely generated $\widehat R$-module $M$,
$\HH^2_{\mm}(M)\simeq\HH^1(U,\widetilde M|_U)$.
\item If $M$ is reflexive, there is a reflexive coherent sheaf
$\mathcal F$ on $W$ with $\mathcal F|_U\simeq\widetilde M|_U$.
\item If $\mathcal F$ is reflexive on $W$, its restriction
$V=\iota^*\mathcal F=\mathcal F\otimes\OO_Y$, with no double dual
taken, is torsion-free of rank $\rk\mathcal F$.
\item If $M$ is reflexive, $\HH^2_{\mm}(M)$ has finite length.
\end{enumerate}
Items (2)--(4) are \cite[Lemmas~2.4, 2.3 and~2.5]{Ch}; compare
\cite[Lemma~6.2 and the following paragraph]{An3}. Item (1) is
standard.
\end{lemma}
\begin{proof}
For the first assertion, the local-cohomology sequence has the segment
\[
 \HH^1(\Spec\widehat R,\widetilde M)
 \longrightarrow\HH^1(U,\widetilde M|_U)
 \longrightarrow\HH^2_{\mm}(M)
 \longrightarrow\HH^2(\Spec\widehat R,\widetilde M),
\]
whose outer terms vanish because $\Spec\widehat R$ is affine. For the
second, $\pi^*\widetilde M$ is a coherent sheaf on $W$ extending
$\widetilde M|_U$. A reflexive module over a regular local ring of
dimension at most two is free, so $\widetilde M|_U$ is locally free.
The sheaf $\mathcal F=(\pi^*\widetilde M)^{\vee\vee}$ therefore
has the required restriction.

For the third assertion, a reflexive sheaf on the regular scheme $W$
is locally free at points of codimension at most two. Thus $V$ is
locally free at the generic and codimension-one points of $Y$.
At a closed point of $Y$, one has $\depth\mathcal F\ge2$, and a
local equation of $Y$ is $\mathcal F$-regular; hence $\depth V\ge1$.
No nongeneric point of $Y$ is associated to $V$, which is consequently
torsion-free of the stated rank.

For the last assertion, use the subring $A$ of
Proposition~\ref{loc:model}. For a nonmaximal prime $\mathfrak p$ of
$A$, every prime $\mathfrak q$ of $\widehat R$ over it has height
$\dim A_{\mathfrak p}$, by going-down and incomparability \cite[Tags 00H8 and 00GT]{St}. The ring
$\widehat R_{\mathfrak q}$ is regular and $M_{\mathfrak q}$ is free.
It follows, by checking a parameter sequence at the maximal ideals of
the finite semilocal $A_{\mathfrak p}$-algebra
$\widehat R\otimes_AA_{\mathfrak p}$, that
$M_{\mathfrak p}$ is maximal Cohen--Macaulay over the regular ring
$A_{\mathfrak p}$, hence free. Thus $\Ext_A^1(M,A)$ has finite
length. Its Matlis dual is
$\HH^2_{\nn}(M)=\HH^2_{\mm}(M)$, where $\nn=(x_0,x_1,x_2)A$, by local duality
\cite[Theorem~3.5.8]{BrHe}, and Matlis duality preserves finite length.
\end{proof}

\subsection{Normalization of slopes}\label{loc:slopes}

\begin{lemma}[A prescribed interval of slopes]\label{loc:normalization}
Let $\mathcal F$ be a nonzero reflexive coherent sheaf on $W$ and let
$u\in\mathbf R$. There is a reflexive coherent sheaf $\mathcal F'$
with $\mathcal F'|_U\simeq\mathcal F|_U$ such that every
Harder--Narasimhan factor $G$ of $\mathcal F'|_Y$ satisfies
$u\le\bar\mu(G)<u+1$.
\end{lemma}
\begin{proof}
The argument is the elementary modification of
\cite[Lemma~3.1 and Proposition~3.3]{Ch}, \cite[Lemma~3.1]{OAI} and
\cite[Lemmas~6.3 and~6.4]{An3}. Those sources use the window
$[3/2,5/2)$, for $K_Y\equiv4H$, or an unspecified window of length
one; the argument is the same for $[u,u+1)$.
Write $V=\mathcal F|_Y$. Since $\mathcal F(mY)|_Y=V(-mH)$, replacing
$\mathcal F$ by $\mathcal F(mY)$ for a suitable integer $m$, we may
assume $\bar\mu_{\max}(V)<u+1$.

If $\bar\mu_{\min}(V)<u$, let $V\twoheadrightarrow Q$ be its last
Harder--Narasimhan quotient and let $V_1$ be the kernel. The sheaf
$Q$ is torsion-free and semistable of slope $\bar\mu_{\min}(V)$.
Set
\[
 \mathcal F_1=\ker(\mathcal F\longrightarrow\iota_*V
                                  \longrightarrow\iota_*Q).
\]
Then $\mathcal F(-Y)\subset\mathcal F_1\subset\mathcal F$, so the
new sheaf is torsion-free and has the same restriction to $U$.
At the generic point of $Y$, the sheaf $\mathcal F_1$ is
torsion-free over a discrete valuation ring, hence free. At a point of
$Y$ of codimension two or three in $W$, the sheaf $\mathcal F$ has
depth at least two and $\iota_*Q$ has depth at least one. The depth
lemma gives $\depth\mathcal F_1\ge2$. Away from $Y$ nothing changes,
so $\mathcal F_1$ is reflexive \cite[Proposition~1.4.1]{BrHe}.
Its restriction is torsion-free by Lemma~\ref{loc:extension}(3).

Restriction to $Y$, using the Cartier resolution and
$\Tor^{\OO_W}_1(\mathcal F,\OO_Y)=0$, gives
\[
 0\longrightarrow Q(H)\longrightarrow\mathcal F_1|_Y
 \longrightarrow V_1\longrightarrow0,
\]
because $\Tor^{\OO_W}_1(\iota_*Q,\OO_Y)=Q(H)$. Hence
\[
 \bar\mu_{\max}(\mathcal F_1|_Y)
 \le\max\{\bar\mu(Q)+1,\bar\mu_{\max}(V_1)\}<u+1.
\]
The $H$-degree increases by $d\rk Q\ge d$:
\[
 c_1(\mathcal F_1|_Y)\cdot H
 =c_1(V)\cdot H+d\rk Q.
\]
Repeating this operation leaves the rank fixed, while the integer
$c_1(V)\cdot H$ increases and remains less than
$d\rk(V)(u+1)$. The process therefore terminates, at which point
$\bar\mu_{\min}(V)\ge u$.
\end{proof}

\begin{lemma}[Euler characteristic in the central interval]
\label{loc:euler}
Let $V$ be torsion-free of rank $r>0$, and suppose every
Harder--Narasimhan factor $G$ satisfies
$|\bar\mu(G)-\lambda/2|\le1/2$. Then
\[
 \chi(V)\le\frac{d-\sigma}{8}r.
\]
\end{lemma}
\begin{proof}
This is the Riemann--Roch estimate of \cite[Remark~9.2]{An3}; for
$K_Y\equiv4H$ compare \cite[Proposition~4.1]{Ch}. For a
Harder--Narasimhan factor $G$ of $V$, of rank $\rho$,
\eqref{surf:orthogonal-part} and \eqref{surf:surface-rr} with $t=0$
give
\[
 \frac{\chi(G)}\rho
 =-\frac\sigma8+\frac{\alpha(G)^2}2
   +\frac d2\left(\bar\mu(G)-\frac\lambda2\right)^2
   -\frac{\delta(G)}2
 \le-\frac\sigma8+\frac d2
                  \left(\bar\mu(G)-\frac\lambda2\right)^2
 \le\frac{d-\sigma}{8},
\]
since $\alpha(G)^2\le0$ and $\delta(G)\ge0$
(Proposition~\ref{surf:pushforward-identity}).
Summing over the Harder--Narasimhan filtration proves the assertion.
\end{proof}

Lemma~\ref{loc:euler} needs no proportionality between $K_Y$ and
$H$: the Hodge index theorem absorbs the component of
$\xi-K_Y/2$ orthogonal to $H$. Only the width and center of the slope interval
enter the estimate.

By~\eqref{surf:pushforward-slopes}, the window
$(\lambda-1)/2\le\bar\mu(G)<(\lambda+1)/2$ used in
Section~\ref{loc:qone} is the window $-2\le b(f_*G)<-1$ for the
average slope of the pushforward to $\mathbf P^2$, inside the slope
range $-2\le\mu_j\le-1$ that defines $J_0$ in
Proposition~\ref{surf:plane-h1}. This concerns the average slope only:
the Harder--Narasimhan factors of $f_*G$ need not have slopes in
$[-2,-1]$. For instance, if $K_Y\equiv4H$ and $G=\OO_Y(2H)$, then
$\bar\mu(G)=2$ and $b(f_*G)=-3/2$, but $\OO_{\mathbf P^2}(2)\subset
f_*G$, so $\mu_{\max}(f_*G)\ge2$. The proof of Lemma~\ref{loc:euler}
takes place on $Y$ and does not use $f_*G$.

\subsection{The two tails}\label{loc:tails}

\begin{lemma}[Two-section convexity]\label{loc:two-sections}
Let $G$ be a torsion-free sheaf on $Y$ and set
$g_m=h^0(Y,G(-mH))$. Then $2g_{m+1}-g_{m+2}\le g_m$ for every
$m\in\mathbf Z$. Consequently $g_1+g_2\le g_0+g_3$, and in
particular $g_1+g_2\le g_0$ if $g_3=0$.
\end{lemma}
\begin{proof}
Choose two general sections $z,z'$ of $\OO_Y(H)$. Their common zero
scheme is finite and avoids the finite non-locally-free locus of $G$.
The Koszul complex
\[
 0\longrightarrow G(-2H)
 \xrightarrow{(z',-z)}G(-H)^{\oplus2}
 \xrightarrow{(z,z')}G
\]
is exact in its first two positions: away from the common zeros it
is split exact, and at a common zero $G$ is free and $z,z'$ is a
regular sequence. Twisting by $-mH$ and taking global sections, which is left exact,
gives $2g_{m+1}-g_{m+2}\le g_m$. For $m=0,1$ we have
$g_0\ge2g_1-g_2$ and $g_1-2g_2+g_3\ge0$, whence
\[
 g_0+g_3\ge g_1+g_2+(g_1-2g_2+g_3)\ge g_1+g_2.
\]
For $g_3=0$ this is \cite[Lemma~4.2]{Ch}.
\end{proof}

\begin{lemma}[Sections at nonpositive slope]\label{loc:nonpositive}
If every positive-rank subsheaf of a torsion-free sheaf $G$ on $Y$
has slope at most zero, then $h^0(Y,G)\le\rk G$.
\end{lemma}
\begin{proof}
Induct on the rank. There is nothing to prove if $G$ has no section.
Otherwise a nonzero section gives $\OO_Y\hookrightarrow G$. Let
$\mathcal L$ be the saturation of its image. Its double dual is
$\OO_Y(D)$, where $D$ is the effective divisor of the induced section.
The slope assumption gives $D\cdot H\le0$, so ampleness forces $D=0$.
The inclusions
$\OO_Y\subset\mathcal L\subset\mathcal L^{\vee\vee}\simeq\OO_Y$
therefore give $\mathcal L=\OO_Y$, and $G'=G/\OO_Y$ is torsion-free.
Every positive-rank subsheaf of $G'$ still has slope at most zero:
the inverse image of a positive-degree subsheaf would have the same
first Chern class and rank larger by one, contrary to the hypothesis
on $G$. By induction and connectedness of $Y$,
$h^0(G)\le h^0(\OO_Y)+h^0(G')\le1+\rk G'=\rk G$.
\end{proof}

For a coherent sheaf $V$ define
\[
 T_-(V)=\sum_{m\ge1}h^0(Y,V(-mH)),\qquad
 T_+(V)=\sum_{m\ge1}h^2(Y,V(mH)).
\]

\begin{lemma}[Tail bounds, including the endpoint]\label{loc:tail-bound}
Let $V$ be a torsion-free sheaf of rank $r$ on $Y$. Suppose
$\lambda\le5$ and every Harder--Narasimhan factor $G$ of $V$
satisfies
\[
 \frac{\lambda-1}{2}\le\bar\mu(G)<\frac{\lambda+1}{2}.
\]
Then both tails are finite and
\[
 T_-(V)\le h^0(V),\qquad T_+(V)\le h^2(V)+2\varepsilon r,
 \quad
 \varepsilon=\begin{cases}0,&\lambda<5,\\1,&\lambda=5.\end{cases}
\]
For $K_Y\equiv4H$, compare \cite[Proposition~4.3]{Ch}.
\end{lemma}
\begin{proof}
For $m\ge3$, a nonzero section of $V(-mH)$ would give a slope-zero
subsheaf, whereas
$\bar\mu_{\max}(V(-mH))=\bar\mu_{\max}(V)-m\le\bar\mu_{\max}(V)-3
<(\lambda+1)/2-3\le0$. Thus $h^0(V(-mH))=0$ for $m\ge3$, and
Lemma~\ref{loc:two-sections} gives $T_-(V)\le h^0(V)$.

Put $V^D=\Hom(V,\omega_Y)$, a locally free sheaf of rank $r$, and
$g_m=h^0(V^D(-mH))$. Serre duality
\cite[Theorem~III.7.6]{Ha} identifies $h^2(V(mH))$ with $g_m$.
Every positive-rank subsheaf $\mathcal S\subset V^D$ satisfies
\begin{equation}\label{loc:dual-slope}
 \bar\mu(\mathcal S)\le\lambda-\bar\mu_{\min}(V)
                         \le\frac{\lambda+1}{2}.
\end{equation}
Indeed, we may assume $\mathcal S$ saturated, since saturation does not
decrease $\bar\mu(\mathcal S)$, and dualize $0\to\mathcal S\to V^D\to V^D/\mathcal S\to0$ with
respect to $\omega_Y$. Since $\Hom(V^D,\omega_Y)=V^{\vee\vee}$ and
$V^D/\mathcal S$ is torsion-free, the cokernel of
$V^{\vee\vee}\to\Hom(\mathcal S,\omega_Y)$ lies in
$\mathcal Ext^1(V^D/\mathcal S,\omega_Y)$, which has finite support. The image
$Q$ of the composite from $V$ is therefore a torsion-free quotient of
rank $\rho=\rk\mathcal S$ with $c_1(Q)=\rho K_Y-c_1(\mathcal S)$. Its slope is at least
$\bar\mu_{\min}(V)$, proving \eqref{loc:dual-slope}.

Thus $g_m=0$ if $m>(\lambda+1)/2$. If $\lambda<5$, this holds for
$m\ge3$, and Lemma~\ref{loc:two-sections} gives
$T_+(V)=g_1+g_2\le g_0=h^2(V)$. At $\lambda=5$ it holds for
$m\ge4$, and the same lemma gives
\[
 T_+(V)=g_1+g_2+g_3\le g_0+2g_3.
\]
Every positive-rank subsheaf of $V^D(-3H)$ has slope at most zero
by \eqref{loc:dual-slope}, so Lemma~\ref{loc:nonpositive} gives
$g_3\le r$. This proves the endpoint bound with the loss $2r$.
\end{proof}

The restriction on $\lambda$ is used only in this tail argument:
strictly below the endpoint, $(\lambda+1)/2<3$; at the endpoint,
the remaining term is controlled by $h^0(V^D(-3H))\le r$.

\subsection{Persistence}\label{loc:persistence-subsection}

\begin{lemma}[Persistence of negative Euler characteristic]
\label{loc:persistence}
Let $\mathcal F$ be torsion-free on $W$ and $V=\mathcal F|_Y$.
If both tails of $V$ are finite, then
\[
 \dim_{\CC}\HH^1(U,\mathcal F|_U)
 \ge h^1(V)-T_+(V)-T_-(V).
\]
\end{lemma}
\begin{proof}
We follow the persistence argument of
\cite[Proposition~5.1 and Corollary~5.4]{Ch} and
\cite[Proposition~7.3]{An3}, and give the proof for torsion-free
$\mathcal F$. For $i\ge1$ and every integer $m$,
$\HH^i(W,\mathcal F(-mY))$ is a finite $\widehat R$-module by
properness and is supported at the closed point. It therefore has
finite length, equal to its complex dimension. Relative Serre
vanishing \cite[Tag 02O1]{St} gives its vanishing for $m\gg0$;
the threshold may depend on $\mathcal F$.

The sequences \eqref{loc:restriction-sequence} give
\[
 \length\HH^2(W,\mathcal F(-mY))
 \le\length\HH^2(W,\mathcal F(-(m+1)Y))+h^2(V(mH)).
\]
Descending induction for $m\ge1$ yields
$\length\HH^2(W,\mathcal F(-Y))\le T_+(V)$. The sequence for
$m=0$ then gives
\[
 \length\HH^1(W,\mathcal F)
 \ge h^1(V)-\length\HH^2(W,\mathcal F(-Y))
 \ge h^1(V)-T_+(V).
\]
For $m\ge1$, the quotient $\mathcal F(mY)/\mathcal F$ has a
filtration with factors $\iota_*(V(-kH))$, $1\le k\le m$.
The kernel of
$\HH^1(W,\mathcal F)\to\HH^1(W,\mathcal F(mY))$ is a quotient
of the global sections of this quotient sheaf, so its dimension is
at most $\sum_{k=1}^m h^0(V(-kH))\le T_-(V)$.

Since $j$ is affine and $\mathcal F$ is torsion-free,
$j_*(\mathcal F|_U)=\varinjlim_m\mathcal F(mY)$. Cohomology
commutes with these filtered colimits on the Noetherian scheme $W$
\cite[Tag 01FF]{St}. Therefore (see also \cite[Tag 09MR]{St})
\[
 \HH^1(U,\mathcal F|_U)
 =\HH^1(W,j_*(\mathcal F|_U))
 =\varinjlim_m\HH^1(W,\mathcal F(mY)).
\]
The kernel of the map from $\HH^1(W,\mathcal F)$ to this colimit is the
increasing union of the preceding kernels, hence still has dimension
at most $T_-(V)$. Subtracting it proves the assertion.
\end{proof}

\subsection{The inequality for the original polarization}
\label{loc:qone}

Suppose $\lambda\le5$ and let $\varepsilon$ be as in
Lemma~\ref{loc:tail-bound}. We first prove
\begin{equation}\label{loc:qone-bound}
 \length\HH^2_{\mm}(M)
 \ge\left(\frac{\sigma-d}{8}-2\varepsilon\right)\rk M
\end{equation}
for every finitely generated $\widehat R$-module $M$.
If $M$ is reflexive of rank $r>0$, Lemmas~\ref{loc:extension} and
\ref{loc:normalization}, with $u=(\lambda-1)/2$, give a reflexive
extension $\mathcal F$ whose torsion-free restriction $V$ has rank
$r$ and all Harder--Narasimhan slopes in
$[(\lambda-1)/2,(\lambda+1)/2)$. Lemmas~\ref{loc:extension}(1),
\ref{loc:persistence}, \ref{loc:tail-bound} and~\ref{loc:euler} then
give
\begin{align*}
 \length\HH^2_{\mm}(M)
 &=\dim_{\CC}\HH^1(U,\mathcal F|_U)\\
 &\ge h^1(V)-T_+(V)-T_-(V)\\
 &\ge h^1(V)-h^2(V)-h^0(V)-2\varepsilon r\\
 &=-\chi(V)-2\varepsilon r
 \ge\left(\frac{\sigma-d}{8}-2\varepsilon\right)r.
\end{align*}

For a general $M$ we argue as in \cite[Remark~2.10]{J}. Let
$M_{\mathrm{tor}}$ be the torsion submodule. Since $\widehat R$ is
normal, the injection $M/M_{\mathrm{tor}}\to M^{**}$ is an
isomorphism at primes of height at most one, so its cokernel has dimension at most
one. Also $\dim M_{\mathrm{tor}}\le2$. Consequently both maps
\[
 \HH^2_{\mm}(M)\longrightarrow
 \HH^2_{\mm}(M/M_{\mathrm{tor}})\longrightarrow
 \HH^2_{\mm}(M^{**})
\]
are surjective: their cokernels embed in
$\HH^3_{\mm}(M_{\mathrm{tor}})=0$ and
$\HH^2_{\mm}(M^{**}/(M/M_{\mathrm{tor}}))=0$, respectively.
The reflexive hull has the same rank, proving \eqref{loc:qone-bound}
for all modules, including rank zero.

If $M$ is finite over $R_{\mm}$, flat base change and the fact that
local cohomology is $\mm$-torsion give
$\HH^2_{\mm}(M)\simeq\HH^2_{\mm}(M\otimes\widehat R)$, while
completion preserves rank. This proves the uncompleted assertion
for $q=1$.

\subsection{Veronese subrings}\label{loc:veronese-subsection}

To obtain the constants $q(\sigma-q^2d)/8$, we pass to Veronese
subrings. For an integer $q\ge1$, let $R^{(q)}=\bigoplus_{t\ge0}R_{qt}$, with $R_{qt}$
placed in degree $t$. This is the section ring of $(Y,qH)$; write
$\mm^{(q)}$ for its irrelevant ideal and $\widehat{R^{(q)}}$ for
its completion at that ideal.

\begin{lemma}[Veronese subrings]\label{loc:veronese}
The ring $R$ is finite and torsion-free of rank $q$ over $R^{(q)}$,
and $\sqrt{\mm^{(q)}R}=\mm$. The induced inclusion
$\widehat{R^{(q)}}\to\widehat R$ is finite and local, and every
finite $\widehat R$-module $M$ satisfies
\[
 \rk_{\widehat{R^{(q)}}}M=q\rk_{\widehat R}M,
 \qquad \HH^i_{\mm^{(q)}}(M)=\HH^i_{\mm}(M)\quad(i\ge0),
\]
with the same lengths over the two completed rings.
\end{lemma}
\begin{proof}
Every homogeneous $z\in R$ has $z^q\in R^{(q)}$, so the finite-type
complex algebra $R$ is integral, hence finite, over $R^{(q)}$, and
the asserted radical equality follows. As an $R^{(q)}$-module,
\[
 R=\bigoplus_{i=0}^{q-1}R_{[i]},\qquad
 R_{[i]}=\bigoplus_tR_{qt+i}.
\]
Every summand is nonzero because $H$ is globally generated. For a
nonzero homogeneous $z\in R_{[i]}$, multiplication by $z^{q-1}$ embeds
$R_{[i]}$ in $R^{(q)}$, so each summand is torsion-free of rank one.
This proves the rank assertion before completion.

The $\mm^{(q)}R$-adic and $\mm$-adic topologies agree. Hence
\[
 \widehat R=R\otimes_{R^{(q)}}\widehat{R^{(q)}}
\]
is finite of rank $q$, and the completed inclusion is injective by
flatness of completion. Since $\widehat R$ and
$\widehat{R^{(q)}}$ are domains (Proposition~\ref{loc:model} for $H$
and for $qH$), ranks of finite $\widehat R$-modules multiply by $q$
under restriction of scalars. The ideal $\mm^{(q)}\widehat R$ is
$\mm$-primary, giving equality of local cohomology. Both residue
fields are $\CC$, so lengths, finite or infinite, agree.
\end{proof}

\begin{proof}[Proof of Theorem~\ref{loc:inequality}]
For $K_Y\cdot H\le5qd$, the divisor $qH$ is ample and globally
generated, has square $q^2d$, and satisfies
$K_Y\cdot(qH)\le5(qH)^2$, with strict inequality exactly when
$K_Y\cdot H<5qd$; the signature does not change. Regard $M$ as
a finite module over $\widehat{R^{(q)}}$. Applying
\eqref{loc:qone-bound} to $(Y,qH)$ and then Lemma~\ref{loc:veronese}
gives
\[
 \length\HH^2_{\mm}(M)
 \ge\left(\frac{\sigma-q^2d}{8}-2\varepsilon_q\right)
       \rk_{\widehat{R^{(q)}}}M
 =q\left(\frac{\sigma-q^2d}{8}-2\varepsilon_q\right)
       \rk_{\widehat R}M.
\]
The uncompleted assertion follows by completion as above.
Finiteness for reflexive $\widehat R$-modules is
Lemma~\ref{loc:extension}; a reflexive $R_{\mm}$-module $M$ has
reflexive completion $M\otimes\widehat R$ by flat base change, with
the same local cohomology.
\end{proof}

\begin{remark}[Restriction of scalars]\label{loc:restriction-scalars}
Let $B'\subset B$ be a finite local inclusion of local domains with
the same residue field. A system of parameters $\mathbf x$ of $B'$
is one of $B$, and the Koszul homology of a $B$-module on $\mathbf x$
has the same length over the two rings. Moreover
$\nu_{B'}(M)\ge\nu_B(M)$. Thus a lim Cohen--Macaulay sequence over
$B$ remains one over $B'$ after restriction of scalars; this is a
special case of \cite[Proposition~4.10(c)]{BHM}. Applied to the
completed Veronese inclusion, this transfers nonexistence from the
cone over $(Y,qH)$ to the cone over $(Y,H)$. Lemma~\ref{loc:veronese}
gives more: if every finite $\widehat{R^{(q)}}$-module $N$ satisfies
$\length\HH^2_{\mm^{(q)}}(N)\ge c\rk N$, then every finite
$\widehat R$-module $M$ satisfies
$\length\HH^2_{\mm}(M)\ge qc\rk_{\widehat R}M$.
\end{remark}

\section{From local cohomology to Koszul homology}
\label{transfer:section}

This section turns lower bounds for the second local cohomology of
reflexive modules into lower bounds for the Koszul homology of
arbitrary modules. The local cohomology of the module itself cannot
be used (Example~\ref{transfer:example}).
Proposition~\ref{transfer:replacement} replaces a finite module $M$
by a reflexive module $M'$ whose second local cohomology is bounded
by a multiple of $\kappa_{\mathbf x}(M)$, and
Proposition~\ref{transfer:gap} draws the consequence.

\subsection{The local setting}\label{transfer:setting}
Let $(B,\mm)$ be a complete normal local domain of dimension three
containing a field, and let $\Bbbk\subset B$ be a coefficient field.
For an arbitrary system of parameters $\mathbf x=(x_0,x_1,x_2)$,
the Cohen structure theorem \cite[Theorem~A.22]{BrHe} gives a finite
inclusion of a power-series ring
\[
 A=\Bbbk[[x_0,x_1,x_2]]\subset B,
 \qquad\nn=(\mathbf x)A.
\]
Since $\sqrt{\nn B}=\mm$, local cohomology modules with
supports in $\nn$ and in $\mm$ coincide, and so do depths
over $A$ and $B$. Their residue fields are both $\Bbbk$, so finite
lengths agree as well. Put
\begin{equation}\label{transfer:constants}
 e=\rk_AB,\qquad h=\length\HH^2_{\mm}(B),\qquad
 a_i=\dim_{\Bbbk}\Tor_i^A(\Bbbk,B),\qquad
 C_{\mathbf x}=1+a_0+a_1+a_0^2.
\end{equation}
Here $a_0=\length(B/(\mathbf x)B)$, and $e$ is the multiplicity of
the parameter ideal on $B$. By Lemma~\ref{transfer:base-ring}, $a_1=a_0-e$;
hence these constants depend only on $B$ and $\mathbf x$, not on the chosen
coefficient field.

For a homologically bounded-below complex $X$ with finitely
generated homology,
write $\beta_i^A(X)=\dim_{\Bbbk}\HH_i(\Bbbk\otimes_A^{\mathbf L}X)$.
For a finite module $M$, this is $\dim_{\Bbbk}\Tor_i^A(\Bbbk,M)$.
Because the Koszul complex resolves $\Bbbk$ over $A$,
\begin{equation}\label{transfer:koszul-betti}
 \length\HH_i(\mathbf x;M)=\beta_i^A(M),\qquad
 \kappa_{\mathbf x}(M)=\max_{1\le i\le3}\beta_i^A(M),\qquad
 \nu_B(M)\le\beta_0^A(M).
\end{equation}

\begin{lemma}[The ring $B$ over $A$]\label{transfer:base-ring}
One has $\operatorname{pd}_AB\le1$, so $a_i=0$ for $i\ge2$ and
$a_0-a_1=e$. Moreover $\Ext_A^1(B,A)$ has finite length $h$,
$a_1\le h$, and
\[
 C_{\mathbf x}\le1+e+2h+(e+h)^2.
\]
\end{lemma}
\begin{proof}
Normality implies $S_2$, so $\depth B\ge2$ and
Auslander--Buchsbaum gives $\operatorname{pd}_AB=3-\depth B\le1$.
The alternating rank of a free resolution yields $a_0-a_1=e$.
For a nonmaximal prime $\mathfrak p$ of $A$, the semilocal ring
$B_{\mathfrak p}$ is normal of dimension at most two, hence
Cohen--Macaulay at its maximal ideals. These all have height
$\dim A_{\mathfrak p}$ by going-down and incomparability \cite[Tags 00H8 and 00GT]{St}.
It follows that $B_{\mathfrak p}$ is maximal Cohen--Macaulay over
the regular ring $A_{\mathfrak p}$ and is free. Consequently
$\Ext_A^1(B,A)$ has finite length. Local duality
\cite[Theorem~3.5.8]{BrHe} identifies its Matlis dual with
$\HH^2_{\nn}(B)=\HH^2_{\mm}(B)$; hence $h$ is finite and
equals $\length\Ext_A^1(B,A)$.

Dualizing the minimal resolution
$0\to A^{a_1}\to A^{a_0}\to B\to0$ presents $\Ext_A^1(B,A)$
as the cokernel of a map $A^{a_0}\to A^{a_1}$ with entries in
$\nn$. Its number of generators is therefore $a_1$,
so $a_1\le h$ and $a_0=e+a_1\le e+h$. Substitution proves the
stated bound on $C_{\mathbf x}$.
\end{proof}

\subsection{Duality and reflexive replacement}\label{transfer:duality}
We index complexes homologically, with
$\HH_i(X[j])=\HH_{i-j}(X)$, and place a module in degree zero.
For a complex $X$ of $B$-modules with bounded, finitely generated
homology, set
\[
 \omega^\bullet=\RHom_A(B,A)[3],\qquad
 \mathbb D(X)=\RHom_B(X,\omega^\bullet)
       \simeq\RHom_A(X,A)[3].
\]
The second identification is adjunction, and regularity of $A$
gives $\mathbb D(\mathbb D(X))\simeq X$. We use the following consequences:
\begin{enumerate}
\item[(D1)] For a finite module $M$,
$\HH_i(\mathbb D(M))=\Ext_A^{3-i}(M,A)$, vanishing outside
$i=0,1,2,3$. In particular $\omega^\bullet$ has homology only in
degrees two and three, with
$\length\HH_2(\omega^\bullet)=h$ and
$\HH_3(\omega^\bullet)=\Hom_A(B,A)$.
\item[(D2)] One has
$\Bbbk\otimes_A^{\mathbf L}\mathbb D(M)\simeq\RHom_A(M,\Bbbk)[3]$, and hence
$\beta_i^A(\mathbb D(M))=\beta_{3-i}^A(M)$.
\item[(D3)] For a finite module $L$,
$\mathbb D(L[3])\simeq\RHom_A(L,A)$ and
$\HH^i_{\mm}(\mathbb D(L[3]))=0$ for $i\ne3$.
\end{enumerate}
In (D3), $\HH^i_{\mm}$ denotes the cohomology of $R\Gamma_{\mm}$;
in our homological indexing,
$\HH^i_{\mm}(X)=\HH_{-i}(R\Gamma_{\mm}X)$. To prove (D3), recall that
$R\Gamma_{\nn}(A)\simeq E_A(\Bbbk)[-3]$, where $E_A(\Bbbk)$ is the
injective hull of $\Bbbk$; thus it is concentrated in homological degree
$-3$.
Since $L$ has a finite free resolution over $A$,
\[
 R\Gamma_{\nn}\RHom_A(L,A)
 \simeq\RHom_A(L,R\Gamma_{\nn}A)
\]
is $\Hom_A(L,E_A(\Bbbk))[-3]$.

\begin{proposition}[{Reflexive replacement; cf.\ \cite[Definition~2.1, Proposition~2.3 and Lemmas~2.4--2.5]{J}}]
\label{transfer:replacement}
For every finite $B$-module $M$ there are a finite reflexive
$B$-module $M'$ and a homomorphism $M'\to M$ with finite-length
cokernel, such that
\begin{align}
 \rk_BM'&\ge\rk_BM,\label{transfer:replacement-rank}\\
 \length\HH^2_{\mm}(M')
 &\le h\left[\beta_1^A(M)+a_0\beta_2^A(M)
                     +(a_1+a_0^2)\beta_3^A(M)\right]
 \le hC_{\mathbf x}\kappa_{\mathbf x}(M).
 \label{transfer:replacement-bound}
\end{align}
\end{proposition}
\begin{proof}
\emph{Step 1: a minimal resolution.}
By (D1), $\mathbb D(M)$ has finitely generated homology in degrees
$0,\ldots,3$.
Choose a quasi-isomorphism $P\to \mathbb D(M)$ from a complex of finite
free $B$-modules, with $P_i=0$ for $i<0$ and all differentials having
entries in $\mm$. Put $\tau_i=\rk_BP_i$. For $i\ge0$, let
$P_{<i}$ be the subcomplex in degrees less than $i$ and let
$P_{\ge i}=P/P_{<i}$. Then $P_{\ge i}$ has no homology below $i$,
and
$\HH_i(P_{\ge i})=\coker(P_{i+1}\to P_i)$ needs exactly
$\tau_i$ generators over $B$.

\emph{Step 2: the first three ranks.}
Since $P_{\ge i}$ has no homology below degree $i$,
\[
 \HH_i(\Bbbk\otimes_A^{\mathbf L}P_{\ge i})
 =\Bbbk\otimes_A\HH_i(P_{\ge i})
 \twoheadrightarrow \Bbbk\otimes_B\HH_i(P_{\ge i}),
\]
and therefore
\begin{equation}\label{transfer:tau-connective}
 \tau_i\le\beta_i^A(P_{\ge i}).
\end{equation}
The exact sequence
$0\to P_i[i]\to P_{\ge i}\to P_{\ge i+1}\to0$ and
$\beta_j^A(P_i[i])=\tau_i a_{j-i}$, with $a_u=0$ for $u<0$, imply
\begin{equation}\label{transfer:tau-recurrence}
 \beta_j^A(P_{\ge i+1})
 \le\beta_j^A(P_{\ge i})+\tau_i a_{j-1-i}.
\end{equation}
Starting from (D2),
$\beta_j^A(P_{\ge0})=\beta_{3-j}^A(M)$, use
\eqref{transfer:tau-connective} and
\eqref{transfer:tau-recurrence} successively to obtain
\[
 \tau_0\le\beta_3^A(M),\qquad
 \tau_1\le\beta_2^A(M)+a_0\tau_0,\qquad
 \tau_2\le\beta_1^A(M)+a_1\tau_0+a_0\tau_1.
\]
In particular
\begin{equation}\label{transfer:tau-two}
 \tau_2\le\beta_1^A(M)+a_0\beta_2^A(M)
                         +(a_1+a_0^2)\beta_3^A(M).
\end{equation}

\emph{Step 3: constructing the replacement.}
Since $\HH_j(P)=0$ for $j\ge4$, the complex $P_{\ge3}$ has
homology only in degree three. Thus
$P_{\ge3}\simeq L[3]$, where $L=\coker(P_4\to P_3)$.
Dualizing $0\to P_{<3}\to P\to P_{\ge3}\to0$ gives a triangle
\begin{equation}\label{transfer:replacement-triangle}
 \mathbb D(L[3])\longrightarrow M\longrightarrow \mathbb D(P_{<3})
                         \longrightarrow \mathbb D(L[3])[1].
\end{equation}
The factors of $P_{<3}$ are $P_i[i]$ for $i=0,1,2$, whose duals
are $(\omega^\bullet)^{\tau_i}[-i]$. By (D1) these have homology
only in degrees $2-i,3-i$. It follows that $\mathbb D(P_{<2})$ has homology
only in degrees $1,2,3$. Dualizing the exact triangle
$P_{<2}\to P_{<3}\to P_2[2]\to P_{<2}[1]$ shows that $\mathbb D(P_{<3})$ has no negative homology and that
$\HH_0(\mathbb D(P_{<3}))$ is a quotient of
$\HH_2(\omega^\bullet)^{\tau_2}$. Hence
\begin{equation}\label{transfer:finite-error}
 \length\HH_0(\mathbb D(P_{<3}))\le h\tau_2.
\end{equation}
On the other hand
$\HH_j(\mathbb D(L[3]))=\Ext_A^{-j}(L,A)=0$ for $j>0$. Set
\[
 M'=\HH_0(\mathbb D(L[3]))=\Hom_A(L,A),\qquad
 Z=\HH_{-1}(\mathbb D(L[3])).
\]
The homology sequence of \eqref{transfer:replacement-triangle}
gives
\[
 M'\longrightarrow M\longrightarrow\HH_0(\mathbb D(P_{<3}))
 \longrightarrow Z\longrightarrow0,
\]
and, for $j\le-2$,
$\HH_j(\mathbb D(L[3]))\simeq\HH_{j+1}(\mathbb D(P_{<3}))=0$.
Thus the cokernel of $M'\to M$ has finite length, proving
\eqref{transfer:replacement-rank}, while $Z$ has length at most
$h\tau_2$. The only possible homology of $\mathbb D(L[3])$ is now $M'$ in
degree zero and $Z$ in degree minus one.

\emph{Step 4: local cohomology.}
There is a triangle
$M'\to \mathbb D(L[3])\to Z[-1]\to M'[1]$. Apply $R\Gamma_{\mm}$.
By (D3) its middle term has local cohomology only in degree three,
whereas the finite-length module $Z[-1]$ has local cohomology $Z$
in degree one and zero otherwise. The long exact sequence gives
\[
 \HH^0_{\mm}(M')=\HH^1_{\mm}(M')=0,
 \qquad\HH^2_{\mm}(M')\simeq Z.
\]
Together with \eqref{transfer:finite-error} and
\eqref{transfer:tau-two}, this proves
\eqref{transfer:replacement-bound}.

\emph{Step 5: reflexivity.}
An $A$-presentation $A^{s_1}\to A^{s_0}\to L\to0$ expresses
$M'=\Hom_A(L,A)$ as the kernel of a homomorphism between finite
free $A$-modules. The image of this homomorphism, a submodule of a
free module, is torsion-free, so the depth lemma gives
\[
 \depth M'_{\mathfrak p}\ge\min\{2,\dim A_{\mathfrak p}\}
 \quad\text{for every prime }\mathfrak p\subset A.
\]
The module $M'$ is $A$-torsion-free and hence $B$-torsion-free:
the minimal polynomial of a nonzero $b\in B$ over the fraction field
of $A$ has coefficients in the normal ring $A$ and a nonzero
constant term, so $by=0$ with $y\in M'$ would make that constant term
annihilate $y$.
For $\mathfrak q\subset B$ over $\mathfrak p\subset A$, one has
$\dim B_{\mathfrak q}=\dim A_{\mathfrak p}$. An $M'_{\mathfrak p}$-regular sequence in $\mathfrak pA_{\mathfrak p}$
remains $M'_{\mathfrak q}$-regular, by flatness of localization and
Nakayama's lemma, so the same $S_2$ bound holds over $B$.
Torsion-freeness and $S_2$ over the normal domain $B$ imply
reflexivity \cite[Proposition~1.4.1]{BrHe}.
\end{proof}

The numbers $\tau_i$ in this proof are the Bass numbers of $M$ over
$B$. Indeed $\RHom_B(\Bbbk,\omega^\bullet)\simeq\RHom_A(\Bbbk,A)[3]
\simeq\Bbbk$, the complex $P\otimes_B\Bbbk$ has zero differential, and
$\RHom_B(\Bbbk\otimes_B^{\mathbf L}\mathbb D(M),\omega^\bullet)\simeq
\RHom_B(\Bbbk,M)$, so $\tau_i=\dim_{\Bbbk}\Ext_B^i(\Bbbk,M)$. Likewise
$\beta_{3-i}^A(M)=\dim_{\Bbbk}\Ext_A^i(\Bbbk,M)$ are the Bass numbers
of $M$ over $A$. Thus \eqref{transfer:tau-two} bounds the second Bass
number of $M$ over $B$ by its Bass numbers over $A$, and the proof
gives the sharper form
\[
 \length\HH^2_{\mm}(M')\le h\dim_{\Bbbk}\Ext_B^2(\Bbbk,M).
\]
If $B$ is Cohen--Macaulay, then $h=0$ and $M'$ is maximal
Cohen--Macaulay or zero, as in Example~\ref{transfer:example}.

\begin{proposition}[{A uniform Koszul gap; cf.\ \cite[Corollary~2.11]{J}}]\label{transfer:gap}
Suppose that a constant $c>0$ satisfies
\begin{equation}\label{transfer:reflexive-gap}
 \length\HH^2_{\mm}(N)\ge c\rk_BN
 \quad\text{for every finite reflexive $B$-module }N.
\end{equation}
Then $h\ge c$, and every nonzero finite $B$-module $M$ satisfies
\begin{equation}\label{transfer:uniform-koszul}
 \kappa_{\mathbf x}(M)
 \ge\frac{c}{ehC_{\mathbf x}+c}\,\beta_0^A(M)
 \ge\frac{c}{ehC_{\mathbf x}+c}\,\nu_B(M).
\end{equation}
In particular, $B$ admits no lim Cohen--Macaulay sequence.
\end{proposition}
\begin{proof}
Taking $N=B$ in \eqref{transfer:reflexive-gap} gives $h\ge c$.
For $M'$ from Proposition~\ref{transfer:replacement},
\[
 c\rk_BM\le c\rk_BM'\le\length\HH^2_{\mm}(M')
                         \le hC_{\mathbf x}\kappa_{\mathbf x}(M).
\]
A minimal $A$-free resolution of $M$ has length at most three, so
\[
 \beta_0^A(M)-\beta_1^A(M)+\beta_2^A(M)-\beta_3^A(M)
 =\rk_AM=e\rk_BM.
\]
Its last differential is injective, giving
$\beta_3^A(M)\le\beta_2^A(M)$. Consequently
\[
 \nu_B(M)\le\beta_0^A(M)
 \le e\rk_BM+\beta_1^A(M)
 \le\left(\frac{ehC_{\mathbf x}}c+1\right)
                    \kappa_{\mathbf x}(M),
\]
which is \eqref{transfer:uniform-koszul}, including rank-zero modules.
The system of parameters was arbitrary, so no sequence of nonzero
finite modules can satisfy the lim Cohen--Macaulay condition.
\end{proof}

Theorem~\ref{loc:inequality} supplies \eqref{transfer:reflexive-gap}
for $B=\widehat R$ with $c=c_q$, so Proposition~\ref{transfer:gap}
gives the local assertion of Theorem~\ref{main:criteria} whenever
$c_q>0$ (Section~\ref{ex:local-subsection}).

\begin{remark}[The graded setting]\label{transfer:graded}
Lemma~\ref{transfer:base-ring} and
Propositions~\ref{transfer:replacement} and~\ref{transfer:gap} have
graded forms for a finitely generated graded normal domain $B$ of
dimension three over $\Bbbk=B_0$ and parameters
$x_0,x_1,x_2\in B_1$: the module $M'$ is graded, dimensions over
$\Bbbk$ replace lengths, and the reflexive gap is assumed only for
graded reflexive modules. They are proved in
Appendix~\ref{app:graded} (Lemma~\ref{app:base-ring} and
Propositions~\ref{app:graded-replacement} and~\ref{app:graded-gap})
and give the graded assertion of Theorem~\ref{main:criteria} in
Section~\ref{ex:graded-subsection}.
The companion paper \cite{AnII} uses
Proposition~\ref{transfer:replacement} in \cite[Proposition~2.6]{AnII}.
Its Euler-characteristic gap for modules of depth at least two, in
local and graded form \cite[Proposition~2.2 and Lemma~3.16]{AnII},
rests on a dual construction related to that of
Proposition~\ref{transfer:replacement}.
\end{remark}

\begin{example}[Why the module must be replaced]\label{transfer:example}
Let $S=\CC[x_0,x_1,x_2]$, let $m\ge1$, and let $\Omega_m$ be the
second syzygy of $S/(x_0^m,x_1^m,x_2^m)$. Its minimal resolution is
\[
 0\longrightarrow S(-3m)\longrightarrow S(-2m)^3
 \longrightarrow\Omega_m\longrightarrow0,
\]
and graded local duality, or the Koszul complex, gives
$\HH^2_{S_+}(\Omega_m)\simeq S/(x_0^m,x_1^m,x_2^m)$.
The modules $\Theta_m=S^m\oplus\Omega_m$ have rank $m+2$, and
\[
 \nu_S(\Theta_m)=m+3,\qquad
 \dim_{\CC}\HH_1(\mathbf x;\Theta_m)=1,\qquad
 \HH_2(\mathbf x;\Theta_m)=\HH_3(\mathbf x;\Theta_m)=0,
\]
whereas $\dim_{\CC}\HH^2_{S_+}(\Theta_m)=m^3$.
Thus $\kappa_{\mathbf x}(\Theta_m)/\nu_S(\Theta_m)\to0$ while
$\dim\HH^2_{S_+}(\Theta_m)/\rk\Theta_m\to\infty$.
Hence $\dim_{\CC}\HH^2_{S_+}(M)$ is bounded by no multiple of
$\kappa_{\mathbf x}(M)$, even after normalizing by rank or by
$\nu_S(M)$. Nevertheless, a gap such as \eqref{transfer:reflexive-gap}
excludes
lim Cohen--Macaulay sequences through
Proposition~\ref{transfer:replacement}, which changes the modules.
Over $S$, where $h=0$, the graded replacement $M'$ of $\Theta_m$
given by Proposition~\ref{app:graded-replacement} is free: by that
proposition, $\HH^2_{S_+}(M')=0$, so the
reflexive module $M'$ is maximal Cohen--Macaulay. After completion,
$(\Theta_m)$ is a lim Cohen--Macaulay sequence over
$\CC[[x_0,x_1,x_2]]$, although the second local cohomology of
$\Theta_m$ has length $m^3$ and its rank is only $m+2$. Essentially the same example appears in
\cite[Remark~2.9]{J}, where it shows that existence must be compared
after changing the witnesses.
\end{example}

\section{Proofs of the criteria and examples}\label{ex:section}
This section applies Section~\ref{transfer:section} to the completed
vertex $B=\widehat R$ of $R=R(Y,H)$, with maximal ideal $\mm$.
Sections~\ref{ex:graded-subsection} and~\ref{ex:local-subsection}
prove the graded and local criteria of Theorem~\ref{main:criteria},
Section~\ref{ex:range-subsection} determines when the local criterion
applies, and Section~\ref{ex:families} treats the three families.

\subsection{The graded criterion}\label{ex:graded-subsection}
\begin{proof}[Proof of the graded assertion of Theorem~\ref{main:criteria}]
Let $\mathbf x=(x_0,x_1,x_2)$ be three sections of $\OO_Y(H)$ without a
common zero, and put $S=\CC[x_0,x_1,x_2]\subset R$. Apply Proposition~\ref{app:graded-gap} to $R$ over $S$. Here
\[
 e=\rk_S R=d,\qquad h=\dim_{\CC}\HH^2_{\mm}(R)=h^1_*(\OO_Y),
 \qquad a_0=\dim_{\CC}R/(\mathbf x)R,\qquad a_1=a_0-d.
\]
Theorem~\ref{main:graded} supplies the cohomology-gap hypothesis for
all finite graded modules, with $c=(\sigma-d)/8>0$, or with any positive Veronese constant
$q(\sigma-q^2d)/8$. Consequently
$0<c\le h<\infty$, and Proposition~\ref{app:graded-gap} gives
\[
 \kappa_{\mathbf x}(M)\ge
 \frac{c}{dhC_{\mathbf x}+c}\,\beta_0^S(M)\ge
 \frac{c}{dhC_{\mathbf x}+c}\,\nu_R(M),\qquad
 C_{\mathbf x}=1+a_0+a_1+a_0^2.
\]
The replacement is built in the graded category
(Proposition~\ref{app:graded-replacement}); this is why the graded inequality of
Theorem~\ref{main:graded}, which needs no slope condition, suffices.

For any sequence of nonzero finite graded $R$-modules $M_j$, the
images of $\mathbf x$ form a system of parameters of $R_{\mm}$ and
of $B$. The modules $\HH_i(\mathbf x;M_j)$ and
$M_j/\mm M_j$ are finite-dimensional and annihilated by a power of
$\mm$, so localization and completion do not change them. Their
lengths are their dimensions over $\CC$. The lim Cohen--Macaulay
condition therefore fails for this system of parameters and hence
for every system of parameters \cite[Theorem~4.6]{BHM}.
\end{proof}

\begin{remark}\label{ex:stronger-bundle-constant}
A graded reflexive $R$-module is the section module
$\bigoplus_{t\in\mathbf Z}\HH^0(Y,E(tH))$ of a vector bundle $E$ on
$Y$. Thus the preceding proof can instead use
\cite[Corollary~8.3(a)]{AC}, giving the graded criterion with the
larger constant $c=dL(\sigma/d)\ge(\sigma-d)/6$, where $L$ is the
function of Remark~\ref{surf:stronger-bundle-bound}. The same
corollary bounds the gap moment, in the sense of \cite{AC}, of every vector bundle below by
$\sigma/d$, whereas arithmetically Cohen--Macaulay bundles have gap
moment $1$. By \cite[Remark~8.4(1)]{AC}, the region of the abstract cone of
cohomology tables defined by such a bound is not closed under entrywise
limits; the tables in the example given there are not claimed to come from bundles
on a fixed surface. The proof above uses instead the total length
$h^1_*$ and the graded replacement.
\end{remark}

\begin{corollary}[Separation of total Betti vectors]\label{ex:betti-cone}
Assume $d<\sigma$ and put $\eta_{\mathbf x}=c/(dhC_{\mathbf x}+c)$
as above. Let $\mathcal C\subset\mathbf R^4$ be the convex cone generated by the
vectors $(\beta_0^S(M),\ldots,\beta_3^S(M))$ of finite graded
$R$-modules, where $\beta_i^S(M)=\dim_{\CC}\Tor_i^S(\CC,M)$.
Every point $(b_0,b_1,b_2,b_3)$ in its closure satisfies
\[
 b_1+b_2+b_3\ge\eta_{\mathbf x}b_0.
\]
In particular $(1,0,0,0)$ does not belong to that closure.
\end{corollary}
\begin{proof}
The equality $\beta_i^S(M)=\dim_{\CC}\HH_i(\mathbf x;M)$ and the
first inequality in the preceding proof give
$\beta_1^S(M)+\beta_2^S(M)+\beta_3^S(M)
\ge\max_{i\ge1}\beta_i^S(M)\ge\eta_{\mathbf x}\beta_0^S(M)$.
This linear inequality passes to the cone and to its closure.
\end{proof}

Corollary~\ref{ex:betti-cone} concerns total Betti vectors; entrywise
limits of graded Betti tables behave differently. Over $S$, for
example, the graded Betti tables of $S/(x_0^j)$ converge entrywise to
that of $S$, since the relation moves to degree $j$, although every
total Betti vector is $(1,1,0,0)$.

\subsection{The local criterion}\label{ex:local-subsection}
\begin{proof}[Proof of the local assertion of Theorem~\ref{main:criteria}]
By Proposition~\ref{loc:model}, the completed vertex
$B$ is a normal complete local domain of dimension three with residue field $\CC$. Under
\eqref{main:local-condition}, Theorem~\ref{main:local} gives the
cohomology-gap hypothesis for all finite $B$-modules with the
constant $c=c_q>0$ of that theorem. Proposition~\ref{transfer:gap}
gives, for every system of parameters $\mathbf x$, with
$A=\CC[[\mathbf x]]$, $e=\rk_AB$, $a_0=\length(B/(\mathbf x)B)$,
$a_1=a_0-e$ and $C_{\mathbf x}=1+a_0+a_1+a_0^2$ as in
\eqref{transfer:constants}, and with
$h=\length \HH^2_{\mm}(B)=\dim_{\CC}\HH^2_{\mm}(R)=h^1_*(\OO_Y)$, which
does not depend on $\mathbf x$, the bound
\[
 \kappa_{\mathbf x}(M)\ge\eta_{\mathbf x}\nu_B(M),\qquad
 \eta_{\mathbf x}=\frac{c_q}{ehC_{\mathbf x}+c_q},
\]
and excludes ordinary lim Cohen--Macaulay sequences.

If the parameters come from three sections of $\OO_Y(H)$ without
common zero, then $A=\CC[[x_0,x_1,x_2]]$, $B=R\otimes_S A$, and
\[
 e=d,\qquad a_i=\dim_{\CC}\Tor_i^S(\CC,R).
\]
For any finite $R_{\mm}$-module $M$, completion preserves the
finite-length modules $\HH_i(\mathbf x;M)$ and the number of
generators. Thus the same normalized lower bound holds before
completion, excluding sequences there too; compare
\cite[Proposition~4.12(c)]{BHM}. In particular, neither $B$ nor
$R_{\mm}$ has a small Cohen--Macaulay module, since such a module
forms a constant lim Cohen--Macaulay sequence.
\end{proof}

\subsection{The range of the criterion}\label{ex:range-subsection}
\begin{lemma}[Signature--slope inequality]\label{ex:range}
Assume $d<\sigma$.
\begin{enumerate}[label=\textup{(\alph*)}]
\item $Y$ is of general type, $\lambda>0$, and
$\sigma\le\lambda^2d/9$; in particular $\lambda>3$.
\item An integer $q\ge1$ satisfies \eqref{main:local-condition} if
and only if either $\lambda/5<q<\sqrt{\sigma/d}$, or
$q=\lambda/5$ and $q^2d+16<\sigma$. Every such $q$ satisfies
$\lambda/5\le q<\lambda/3$.
\item Such an integer exists if and only if either
$q_0=\lfloor\lambda/5\rfloor+1$ satisfies $q_0^2d<\sigma$, or
$\lambda/5$ is an integer and $(\lambda/5)^2d+16<\sigma$.
There is none if $5<\lambda\le6$.
\end{enumerate}
\end{lemma}
\begin{proof}
Since $d$ is a positive integer, $\sigma\ge2$. The classification
of surfaces implies that such a surface is of general type:
$\sigma=1$ for $\mathbf P^2$, $\sigma\le0$ for other minimal
surfaces of Kodaira dimension at most one, and blowing up lowers
$\sigma$ by one \cite[Chapter~VI]{BHPV}. For general type,
some positive multiple of $K_Y$ is effective and nonzero, so
$K_Y\cdot H>0$; moreover, the Bogomolov--Miyaoka--Yau inequality gives
$K_Y^2\le3c_2(Y)$ \cite[Theorem~VII.4.1]{BHPV}; it is stated there
for minimal surfaces and persists under blowing up, which lowers
$K_Y^2$ by one and raises $c_2(Y)$ by one. Together with
Noether's formula this says $K_Y^2\le9\chi(\OO_Y)$, or
$\sigma\le K_Y^2/9$. Hodge index gives
$K_Y^2\le(K_Y\cdot H)^2/d=\lambda^2d$, proving (a).

The first assertion of (b) rewrites
\eqref{main:local-condition}: in its second alternative, either
$q>\lambda/5$ and the first alternative holds as well, or
$q=\lambda/5$. The upper bound $q<\lambda/3$ follows from (a).
For (c), $q_0$ is the least integer strictly larger than
$\lambda/5$, and $q^2d$ increases with $q$. If $5<\lambda\le6$,
the interval $[\lambda/5,\lambda/3)$ contains no integer.
\end{proof}
In particular, $q=1$ concerns general-type surfaces with
$3<\lambda\le5$.

\subsection{Three families}
\label{ex:families}
We recall the surfaces of \cite{An1}, following
\cite[Section~4]{An3}.
For each of the three arrangements below, let $s$ be its number of lines.
Their defining linear forms span a three-dimensional space and embed
$\mathbf P^2$ as a plane
$\Pi\subset\mathbf P^{s-1}$. For $n\ge2$, let $X_n$ be its
inverse image under $[z_i]\mapsto[z_i^n]$, let $Y_n\to X_n$ be
the blow-up of the reduced singular locus, let $A_n$ be the
pullback of $\OO_{X_n}(1)$, and let $E_n$ be the reduced
exceptional divisor.

Let $\zeta$ be a primitive cube root of unity. In coordinates
$x,y,z$ on $\mathbf P^2$, the Hesse arrangement consists of the three
coordinate lines and the nine lines
$x+\zeta^ay+\zeta^bz=0$, $0\le a,b\le2$, and has nine quadruple
and twelve double points. The dual Hesse arrangement has nine
lines and twelve triple points; the complete quadrangle has six
lines, four triple points and three double points. The first two
are defined over $\mathbf Q(\zeta)$, and the quadrangle, through four
rational points in general position, over $\mathbf Q$.
The power maps and blow-ups require no additional roots of unity,
so $Y_n,A_n,E_n$ are defined over these fields for every $n$.

In all three cases $A_n\cdot E_n=0$, and
\cite[Proposition~4.1]{An3} gives
\[
\begin{array}{rlll}
 \text{Hesse:}&A_n^2=n^9,&E_n^2=-9n^9,
 &K_{Y_n}=(9n-12)A_n-(2n-3)E_n,\\
 \text{dual Hesse:}&A_n^2=n^6,&E_n^2=-12n^6,
 &K_{Y_n}=(6n-9)A_n-(n-2)E_n,\\
 \text{quadrangle:}&A_n^2=n^3,&E_n^2=-4n^3,
 &K_{Y_n}=(3n-6)A_n-(n-2)E_n.
\end{array}
\]
The signatures are given in \cite[Proposition~4.1(e), Table~2]{An3}.
Let $H_n$ be the \emph{standard} polarization $4A_n-E_n$,
$6A_n-E_n$ or $3A_n-E_n$, respectively
(Table~\ref{ex:family-table}); the standard cones are the section
rings $R(Y_n,H_n)$ and their completed vertices. These divisors are
very ample for all $n\ge2$: the ideal sheaf $\mathcal I_{\Sigma_n}$ of
the reduced singular locus $\Sigma_n$ of $X_n$, twisted by $3$, $5$
and $2$, is globally generated by four cubic monomials, by
square-free monomials of degree five and by three quadratic
monomials, respectively \cite[Proposition~5.1, Lemmas~7.1 and~7.3]{An1}, so the
blow-up very-ampleness criterion of \cite[Lemma~4.2]{An3} applies; see
also \cite[Proposition~4.3 and the proof of Corollary~9.3]{An3}. Apart from
the very-ampleness assertion of
Corollary~\ref{ex:other-polarizations}, the results below use only
that $H_n$ is ample and globally generated. In
the quadrangle family $K_{Y_n}=(n-2)H_n$ in $\operatorname{Pic}Y_n$,
by the display above. For $n=6$ this is $K_Y=4H$, and $(Y_6,3A_6-E_6)$, with $H^2=1080$ and
$\chi(\OO_Y)=1926$, is the example of \cite[Corollary~4.3]{An2},
\cite[Section~7]{Ch} and \cite[Section~3]{J}.

\begin{table}[htbp]
\centering
\small
\begin{tabular}{c c c c c c c}
Arrangement&$H_n$&$d$&$\sigma$&$K\cdot H$&$\lambda$&$d<\sigma$\\\hline
Hesse&$4A_n-E_n$&$7$&$3n^2-11$&$18n-21$&$(18n-21)/7$&$n\ge3$\\
Dual Hesse&$6A_n-E_n$&$24$&$2n^2-13$&$24n-30$&$(4n-5)/4$&$n\ge5$\\
Quadrangle&$3A_n-E_n$&$5$&$(n^2-10)/3$&$5(n-2)$&$n-2$&$n\ge6$
\end{tabular}
\caption{The three families. In each row, $d$, $\sigma$ and
$K\cdot H$ are divided by $n^9$ (Hesse), $n^6$ (dual Hesse) or $n^3$
(quadrangle).}
\label{ex:family-table}
\end{table}

\begin{corollary}\label{ex:standard-cones}
Let $(Y,H)=(Y_n,H_n)$ be a standard pair. If it is a Hesse pair with
$n\ge3$, a dual Hesse pair with $n\ge5$, $n\ne7$, or a
complete-quadrangle pair with $n\ge6$, $n\ne8$, then neither $B$ nor
$R_{\mm}$ admits an ordinary lim Cohen--Macaulay sequence. If
$d<\sigma$, including at the two omitted indices, then $R$ admits no
lim Cohen--Macaulay sequence of graded modules.
\end{corollary}
\begin{proof}
Use Theorem~\ref{main:criteria} and Lemma~\ref{ex:range}, with
$q_0=\lfloor\lambda/5\rfloor+1\le\lambda/5+1$.

For Hesse, $n=3$ gives $\lambda=33/7$, $q_0=1$ and $7<16$;
$n=4$ gives $\lambda=51/7$, $q_0=2$ and $28<37$. If $n\ge5$,
\[
 q_0\le\frac{18n+14}{35},\qquad
 7q_0^2\le\frac{(18n+14)^2}{175}<3n^2-11,
\]
because $201n^2-504n-2121>0$ for $n\ge5$.

For dual Hesse, $n=5,6$ gives $\lambda=15/4,19/4$, $q_0=1$,
and $24<37,59$. For $n=7$, $\lambda=23/4\in(5,6]$, so no $q$
exists. For $n=8,9$, $q_0=2$ and $96<115,149$. If $n\ge10$,
\[
 q_0\le\frac{4n+15}{20},\qquad
 24q_0^2\le\frac{3(4n+15)^2}{50}<2n^2-13,
\]
because $52n^2-360n-1325>0$ for $n\ge10$.

For the quadrangle, $\lambda=n-2$ and the strict inequality is
$15q_0^2<n^2-10$. For $n=6$, $q_0=1$ and $15<26$.
For $n=7$, $\lambda=5$, and the endpoint with $q=1$ applies since
$\sigma-d=8\cdot7^3=2744>16$; the constant $(\sigma-d)/8-2$ is $341$. For $n=8$, $\lambda=6$ and no
$q$ exists. For $n=9,10,11$, $q_0=2$ and $60<71,90,111$.
For $n=12$, $\lambda=10$ and the endpoint with $q=2$ applies,
since $\sigma-4d=(74/3)12^3=42624>16$, with constant
$(\sigma-4d)/8-2=5326$ on the cone over $(Y,2H)$. If $n\ge13$,
\[
 q_0\le\frac{n+3}{5},\qquad
 15q_0^2\le\frac{3(n+3)^2}{5}<n^2-10,
\]
because $2n^2-18n-77>0$ for $n\ge13$.
The graded assertion follows directly from $d<\sigma$.
\end{proof}

For Hesse, the integer $q=\lfloor3n/5\rfloor$ also satisfies the
strict criterion for every $n\ge3$. This is checked directly for $n=3,4$; for $n\ge5$,
$35q>21n-35\ge18n-21$, that is, $q>\lambda/5$, and
$7q^2\le63n^2/25<3n^2-11$. Table~\ref{ex:hesse-constants} lists
$\lambda$, $q_0$ and the constant $c_{q_0}=q_0(\sigma-q_0^2d)/8$ on
the original cone for $3\le n\le8$. At $n=5$, the alternative $q=3$ gives
$c_3=(3/8)5^9$, whereas $q_0=2$ gives $c_2=9\cdot5^9$.
The endpoint criterion can improve on $c_{q_0}$: for $n=7$, $\lambda/5=3$, and
$q=3$ gives
\[
 3\left(\frac{\sigma-9d}{8}-2\right)
 =\frac{219}{8}7^9-6>12\cdot7^9=c_{q_0}.
\]
For graded modules the Veronese form of
Theorem~\ref{main:graded} gives $q(\sigma-q^2d)/8$ for every integer
$q\ge1$ with $q^2d<\sigma$, without a slope condition.

\begin{table}[htbp]
\centering
\small
\begin{tabular}{c c c c c c c}
$n$&$3$&$4$&$5$&$6$&$7$&$8$\\\hline
$\lambda$&$33/7$&$51/7$&$69/7$&$87/7$&$15$&$123/7$\\
$q_0$&$1$&$2$&$2$&$3$&$4$&$4$\\
$q_0^2d/n^9$&$7$&$28$&$28$&$63$&$112$&$112$\\
$\sigma/n^9$&$16$&$37$&$64$&$97$&$136$&$181$\\
$c_{q_0}/n^9$&$9/8$&$9/4$&$9$&$51/4$&$12$&$69/2$
\end{tabular}
\caption{Hesse pairs $(Y_n,4A_n-E_n)$. Every finite $B$-module $M$
satisfies $\length \HH^2_{\mm}(M)\ge c_{q_0}\rk M$.}
\label{ex:hesse-constants}
\end{table}

\begin{remark}[Proportionality of the canonical class]
\label{ex:hesse-nonsubcanonical}
The intersection form is nondegenerate on the span of $A_n,E_n$.
If $K_{Y_n}\equiv t(4A_n-E_n)$ on a Hesse surface, coefficient
comparison gives $t=(9n-12)/4=2n-3$, hence $n=0$. Thus no Hesse pair
with $n\ge3$, and no multiple of its polarization, satisfies the
proportionality hypothesis of \cite{Ch}. Likewise, on the dual Hesse
surfaces $K_{Y_n}-t(6A_n-E_n)=(6n-9-6t)A_n-(n-2-t)E_n$ is nonzero for
every $t$; for the polarizations $J$ of
Corollary~\ref{ex:other-polarizations}, $K_Y\equiv tJ$ would require
$7t=33$ and $t=5$, respectively $4t=18$ and $t=6$; and on the
quadrangle surfaces $K_{Y_n}\equiv(n-2)H_n$, which is $4kH_n$ with
$k\in\mathbf Z$ exactly when $n\equiv2\pmod4$. Hence the quadrangle
cones with $n\equiv2\pmod4$ are the only cones of
Corollaries~\ref{ex:standard-cones} and~\ref{ex:other-polarizations}
with $K_Y\equiv4kH$ for a positive integer $k$. For the Hesse pairs,
the absence of small Cohen--Macaulay modules for every $n\ge3$ is
recorded in
\cite[Corollary~9.3]{An3} as a consequence of
\cite[Theorem~1.2]{OAI}; the present result also excludes
asymptotically negligible positive Koszul homology.
\end{remark}

\begin{remark}[The two omitted standard cones]\label{ex:omitted}
For dual Hesse at $n=7$, $\lambda=23/4$ forces every integer
$q\ge\lambda/5$ to satisfy $q\ge2$; then
$q^2d\ge96\cdot7^6>85\cdot7^6=\sigma$. For the quadrangle at
$n=8$, $\lambda=6$ again requires $q\ge2$, whereas
$q^2d\ge20\cdot8^3>18\cdot8^3=\sigma$. Neither standard cone is
settled by the strict or endpoint local criterion. Graded
sequences are excluded by Theorem~\ref{main:criteria}. Both cones
satisfy $2d<\sigma$, namely $48\cdot7^6<85\cdot7^6$ and
$10\cdot8^3<18\cdot8^3$, so they have no small Cohen--Macaulay
modules \cite[Theorem~1.2]{OAI}; what remains open is the lim
Cohen--Macaulay question for arbitrary modules
(Section~\ref{sec:questions}).
\end{remark}

\begin{corollary}[Other polarizations at the omitted indices]
\label{ex:other-polarizations}
Take $J=7A_7-E_7$ on the dual Hesse surface $Y_7$, or
$J=4A_8-E_8$ on the complete-quadrangle surface $Y_8$. These
polarizations are very ample, and neither the completed vertex $B$
of $R(Y,J)$ nor $R(Y,J)_{\mm}$ admits an ordinary lim Cohen--Macaulay
sequence. Every finite $B$-module $M$ satisfies
$\length \HH^2_{\mm}(M)\ge6\cdot7^6\rk M$ in the first case, and
$\length \HH^2_{\mm}(M)\ge384\rk M$ in the second.
\end{corollary}
\begin{proof}
In each case $J=H_n+A_n$, with $H_n$ very ample and $A_n$
globally generated, so $J$ is very ample. The intersections are
\[
\begin{array}{c|ccc}
 &J^2&K\cdot J&\sigma\\\hline
 Y_7\text{ (dual Hesse)}&37\cdot7^6&171\cdot7^6&85\cdot7^6\\
 Y_8\text{ (quadrangle)}&12\cdot8^3&48\cdot8^3&18\cdot8^3
\end{array}
\]
Thus $J^2<\sigma$; the ratio $K\cdot J/J^2$ equals $171/37$,
respectively $4$, so $K\cdot J<5J^2$; and $(\sigma-J^2)/8$ equals
$6\cdot7^6$, respectively $384$.
Apply Theorems~\ref{main:local} and~\ref{main:criteria} with $q=1$.
\end{proof}

\section{Isolated singularities in every dimension}\label{hd:section}
The reflexive replacement of Proposition~\ref{transfer:replacement}
extends to arbitrary dimension (Proposition~\ref{hd:replacement}).
After restricting the replaced module to a surface section, we can
apply Theorem~\ref{main:local}. This gives complete normal local
domains with residue field $\CC$ and an isolated singularity, without
ordinary lim Cohen--Macaulay sequences, in every dimension at least three.

\subsection{Replacement in arbitrary dimension}
Let $m\ge2$, and let $(B,\mm)$ be a complete normal local domain
of dimension $m+1$ containing a field. Assume that $B_{\mathfrak q}$
is Cohen--Macaulay for every prime $\mathfrak q\ne\mm$. Choose a
coefficient field $\Bbbk\subset B$, a system of parameters
$\mathbf x=(x_0,\ldots,x_m)$, and the corresponding Noether normalization
$A=\Bbbk[[x_0,\ldots,x_m]]\subset B$, which is finite. Write $\nn$ for its maximal
ideal and put
\[
 e=\rk_A B,\qquad a_i=\dim_{\Bbbk}\Tor_i^A(\Bbbk,B),\qquad
 b=\sum_{i\ge0}a_i,\qquad
 h_B=\sum_{j=2}^{m}\length \HH^j_{\mm}(B).
\]
Depth, local cohomology and finite lengths agree over $A$ and $B$.
With $\beta_i^A$ as in Section~\ref{transfer:setting}, the Koszul
complex on $\mathbf x$ resolves $\Bbbk$ over $A$, so, as in
\eqref{transfer:koszul-betti},
\[
 \length \HH_i(\mathbf x;M)=\beta_i^A(M),\qquad
 \kappa_{\mathbf x}(M)=\max_{1\le i\le m+1}\beta_i^A(M).
\]
In particular $a_i=\length \HH_i(\mathbf x;B)$ and $e$, the multiplicity
of $(\mathbf x)B$, depend only on $B$ and $\mathbf x$. For $m=2$, the
punctured Cohen--Macaulay hypothesis holds by Serre's criterion,
$h_B=h$, and $b=a_0+a_1$ by Lemma~\ref{transfer:base-ring}.

By local duality \cite[Theorem~3.5.8]{BrHe},
$\HH^i_{\nn}(F)$ is Matlis dual to $\Ext_A^{m+1-i}(F,A)$.
Either has finite length exactly when the other does, and their
lengths then agree. We also use that if a finite module $F$ over a regular local ring
$\Lambda$ has projective dimension $t\ge1$, then
$\Ext^t_{\Lambda}(F,\Lambda)\ne0$: the dual of the last differential in a
minimal free resolution is not surjective by Nakayama's lemma.

\begin{lemma}[Punctured freeness]\label{hd:punctured-freeness}
\begin{enumerate}[label=\textup{(\alph*)}]
\item $B_{\mathfrak p}$ is free over $A_{\mathfrak p}$ for every
prime $\mathfrak p\ne\nn$. The modules $\Ext_A^i(B,A)$ for $i\ge1$
have finite length, vanish for $i\ge m$, and satisfy
$\sum_{i\ge1}\length\Ext_A^i(B,A)=h_B<\infty$.
\item If $F$ is a finite $B$-module with $\HH^j_{\mm}(F)$ of finite
length for every $j\le m$, then $F_{\mathfrak p}$ is free over
$A_{\mathfrak p}$ for every
prime $\mathfrak p\ne\nn$.
For every prime $\mathfrak q\ne\mm$, $F_{\mathfrak q}$ is zero or
maximal Cohen--Macaulay over $B_{\mathfrak q}$. If $B_{\mathfrak q}$
is regular for every prime $\mathfrak q\ne\mm$, then each
$F_{\mathfrak q}$ is free of rank $\rk_B F$.
\end{enumerate}
\end{lemma}
\begin{proof}
For (a), the maximal ideals of the semilocal ring $B_{\mathfrak p}$
correspond to primes $\mathfrak q$ of $B$ over $\mathfrak p$.
They differ from $\mm$, and have height $\dim A_{\mathfrak p}$
by going-down and incomparability \cite[Tags 00H8 and 00GT]{St}. The rings $B_{\mathfrak q}$ are
Cohen--Macaulay, so $B_{\mathfrak p}$ is maximal Cohen--Macaulay
over the regular local ring $A_{\mathfrak p}$ and is free by
Auslander--Buchsbaum. The positive Ext modules consequently have
finite length. Local duality identifies their lengths with those
of $\HH^{m+1-i}_{\mm}(B)$, which vanish when $m+1-i\le1$ because
$\depth B\ge2$.

For (b), local duality makes every positive $\Ext_A^i(F,A)$ a
finite-length module. Hence all positive Ext modules of
$F_{\mathfrak p}$ over $A_{\mathfrak p}$ vanish away from
$\nn$, forcing $F_{\mathfrak p}$ to be free by the preceding
top-Ext observation. For $\mathfrak q\ne\mm$, its contraction
$\mathfrak p$ differs from $\nn$, since $\mm$ is the only prime
over $\nn$. A system of parameters of $A_{\mathfrak p}$ is
regular on the free module $F_{\mathfrak p}$, and remains regular
on its nonzero localization $F_{\mathfrak q}$. This sequence has
length $\dim A_{\mathfrak p}=\dim B_{\mathfrak q}$. Thus
$F_{\mathfrak q}$ is zero or maximal Cohen--Macaulay. If the ring
is regular, Auslander--Buchsbaum makes it free; its rank agrees
with the generic rank since $B$ is a domain.
\end{proof}

As in Section~\ref{transfer:duality}, we index complexes
homologically and set
\[
 \omega^\bullet=\RHom_A(B,A)[m+1],\qquad
 \mathbb D(X)=\RHom_B(X,\omega^\bullet)
 \simeq\RHom_A(X,A)[m+1].
\]
Then $\mathbb D(\mathbb D(X))\simeq X$ for every complex $X$ with bounded, finitely
generated homology. The duality identities used in the
three-dimensional replacement become:
\begin{enumerate}[label=\textup{(D\arabic*$'$)}]
\item For finite $M$,
$\HH_i(\mathbb D(M))=\Ext_A^{m+1-i}(M,A)$, zero outside $0\le i\le m+1$.
Moreover $\HH_i(\omega^\bullet)=0$ for $i\le1$ and $i>m+1$,
$\sum_{i=2}^{m}\length \HH_i(\omega^\bullet)=h_B$, and
$\HH_{m+1}(\omega^\bullet)=\Hom_A(B,A)$.
\item $\beta_i^A(\mathbb D(M))=\beta_{m+1-i}^A(M)$.
\item For finite $L$, $\mathbb D(L[m+1])\simeq\RHom_A(L,A)$ and
$\HH^i_{\mm}(\mathbb D(L[m+1]))=0$ for $i\ne m+1$.
\end{enumerate}
These follow from regular local duality as (D1)--(D3) do in
Section~\ref{transfer:duality}; the assertions on
$\HH_i(\omega^\bullet)=\Ext_A^{m+1-i}(B,A)$ use
Lemma~\ref{hd:punctured-freeness}(a).

\begin{proposition}[Replacement in arbitrary dimension]
\label{hd:replacement}
For every finite $B$-module $M$ there are a finite reflexive
$B$-module $M'$ and a map $M'\to M$ with finite-length cokernel such that
$\rk_BM'\ge\rk_BM$, $\HH^0_{\mm}(M')=\HH^1_{\mm}(M')=0$, the modules
$\HH^j_{\mm}(M')$ have finite length for $j\le m$, and
\begin{equation}\label{hd:replacement-bound}
 \sum_{j=2}^m\length \HH^j_{\mm}(M')
 \le h_B(1+b)^m\sum_{i=1}^{m+1}\beta_i^A(M)
 \le(m+1)(1+b)^mh_B\kappa_{\mathbf x}(M).
\end{equation}
\end{proposition}
\begin{proof}
The module $M'$ is that of \cite[Definition~2.1]{J}, and the second
bound in \eqref{hd:replacement-bound} follows from
\cite[Lemmas~2.4 and~2.5]{J} for the parameters of $A$, for which the
constant $u$ of \cite[(6)]{J} may be taken equal to one. We give the
argument in our notation, following the proof of
Proposition~\ref{transfer:replacement}. Choose a minimal resolution
$P\to \mathbb D(M)$ by finite free $B$-modules, with $P_i=0$ for $i<0$, set
$\tau_i=\rk_B P_i$, and define $P_{<i}$ and $P_{\ge i}=P/P_{<i}$ as
there. The inequalities \eqref{transfer:tau-connective}
and~\eqref{transfer:tau-recurrence} hold with the same proofs, where
$a_s=0$ for $s<0$. Starting from
$\beta_j^A(P_{\ge0})=\beta_{m+1-j}^A(M)$, which is (D2$'$), they give,
for $0\le i\le m$,
\[
 \tau_i\le\beta_{m+1-i}^A(M)
 +\sum_{l<i}a_{i-1-l}\tau_l
 \le\beta_{m+1-i}^A(M)+b\sum_{l<i}\tau_l.
\]
Induction gives
$\sum_{l\le i}\tau_l\le(1+b)^i
\sum_{l\le i}\beta_{m+1-l}^A(M)$, and therefore
\begin{equation}\label{hd:free-ranks}
 \sum_{i=0}^m\tau_i\le(1+b)^m\sum_{i=1}^{m+1}\beta_i^A(M)
 \le(m+1)(1+b)^m\kappa_{\mathbf x}(M).
\end{equation}

Since $\HH_j(P)=0$ for $j>m+1$, the quotient
$P_{\ge m+1}$ is quasi-isomorphic to $L[m+1]$, where
$L=\coker(P_{m+2}\to P_{m+1})$. Dualizing gives a triangle
\[
 \mathbb D(L[m+1])\longrightarrow M\longrightarrow \mathbb D(P_{<m+1})
 \longrightarrow \mathbb D(L[m+1])[1].
\]
The finite filtration of $P_{<m+1}$ has factors $P_i[i]$ for
$0\le i\le m$, with
$\HH_j(\mathbb D(P_i[i]))=\HH_{j+i}(\omega^\bullet)^{\tau_i}$.
For $j\le0$ these modules have finite length, and the filtration
therefore yields
\begin{equation}\label{hd:dual-length}
 \sum_{j\le0}\length \HH_j(\mathbb D(P_{<m+1}))
 \le\sum_{i=0}^m\tau_i\sum_{l\le i}\length \HH_l(\omega^\bullet)
 \le h_B\sum_{i=0}^m\tau_i.
\end{equation}
Put $M'=\HH_0(\mathbb D(L[m+1]))=\Hom_A(L,A)$.
Since $\HH_j(\mathbb D(L[m+1]))=\Ext_A^{-j}(L,A)$, there is no homology in
positive degrees. The triangle
shows that the cokernel of $M'\to M$ is a submodule of
$\HH_0(\mathbb D(P_{<m+1}))$, that $\HH_{-1}(\mathbb D(L[m+1]))$ is a quotient of
that module, and that
\[
 \HH_j(\mathbb D(L[m+1]))\simeq \HH_{j+1}(\mathbb D(P_{<m+1}))\quad(j\le-2).
\]
Consequently the cokernel has finite length, $\rk M'\ge\rk M$,
and the sum of the lengths of the negative homology of
$\mathbb D(L[m+1])$ is at most $h_B\sum_{i=0}^m\tau_i$.

Let $T$ be the truncation of $\mathbb D(L[m+1])$ to negative homological
degrees. There
is a triangle
\[
 M'\to \mathbb D(L[m+1])\to T\to M'[1].
\]
Since the
homology of $T$ has finite length,
$\HH^i_{\mm}(T)=\HH_{-i}(T)$. Applying derived local cohomology and
using (D3$'$), for $i\le m$ one obtains
\[
 \HH^i_{\mm}(M')\simeq \HH^{i-1}_{\mm}(T)=\HH_{1-i}(T).
\]
In particular each $\HH^i_{\mm}(M')$ with $i\le m$ has finite length.
These vanish for $i=0,1$, and
$\sum_{j=2}^m\length \HH^j_{\mm}(M')\le h_B\sum_{i=0}^m\tau_i$. Now use
\eqref{hd:free-ranks}. Finally, Step~5 of the proof of
Proposition~\ref{transfer:replacement} applies verbatim to
$M'=\Hom_A(L,A)$ and shows that $M'$ is reflexive.
\end{proof}
For $m=2$ this is the module $M'$ of
Proposition~\ref{transfer:replacement}, and \eqref{hd:replacement-bound}
is weaker than \eqref{transfer:replacement-bound}.
By Lemma~\ref{hd:punctured-freeness}(b), $M'$ is maximal
Cohen--Macaulay on the punctured spectrum.

\subsection{Section rings and surface sections}
Let $X$ be smooth, connected and projective of dimension $m\ge2$
over $\CC$, and let $H$ be ample and globally generated. Set
\[
 R(X,H)=\bigoplus_{t\ge0}\HH^0(X,\OO_X(tH)),\qquad
 \mm=R(X,H)_+,\qquad d=H^m,
\]
let $B(X,H)$ be the completion of $R(X,H)_{\mm}$; its maximal ideal
is again denoted $\mm$. The same notation is used for the pairs $(X',H')$ and
$(C,H_C)$ below.

\begin{proposition}[The higher-dimensional cone model]
\label{hd:cone-model}
The ring $R(X,H)$ is a finite-type normal graded domain of
dimension $m+1$. Its completed vertex $B(X,H)$ is a normal
complete local domain of dimension $m+1$ with residue field
$\CC$, regular outside its closed point. If $x_0,\ldots,x_m$ are
sections of $\OO_X(H)$ without common zero, then $R(X,H)$ is
finite over $S=\CC[x_0,\ldots,x_m]$, the ideal $(\mathbf x)R$
is $\mm$-primary, and
\[
 B(X,H)=R(X,H)\otimes_S A,\qquad
 A=\CC[[x_0,\ldots,x_m]],\qquad \rk_A B(X,H)=d.
\]
\end{proposition}
\begin{proof}
The proof of Proposition~\ref{loc:model} applies in every dimension;
compare \cite[Proposition~10.4]{An3}. The sections define a morphism
$f\colon X\to\mathbf P^m$ with $f^*\OO(1)=\OO_X(H)$. It is finite
because $\OO_X(H)$ is ample, and its degree is $H^m=d$.
The full section ring is finite over the polynomial ring $S$ and
is normal, since it is the full section ring of an ample line
bundle on a normal variety. The radical of $(\mathbf x)R$ is
$R_+$, and completion is therefore $R\otimes_S A$ of rank $d$.
The punctured cone is the complement of the zero section in the
total space of $\OO_X(-H)$, hence is smooth. Normality and
regularity away from the vertex persist on completion: $R_{\mm}$ is
excellent, so its completion homomorphism is regular, the completion
is normal, and regularity ascends along this homomorphism
\cite[Tags 0AH2, 0C23 and 031E]{St}. The completion is local and
normal, hence a domain, and its dimension and residue field are
$m+1$ and $\CC$.
\end{proof}
In particular this $B$ satisfies the hypotheses of the preceding
subsection. By \cite[Theorem~A4.1]{Eis}, applied to $R$ as a module
over $S$,
$\HH^{j+1}_{\mm}(R)\simeq\bigoplus_{t\in\mathbf Z}\HH^j(X,\OO_X(tH))$ for
$j\ge1$. For $1\le j\le m-1$ the summands vanish for $t\gg0$ by Serre
vanishing and for $t\ll0$ by Serre duality, so these modules have
finite length; completion does not change them. Hence
\[
 h_B=\sum_{j=1}^{m-1}\sum_{t\in\mathbf Z}
 h^j(X,\OO_X(tH)).
\]

\begin{lemma}[Hyperplane sections]\label{hd:hyperplane}
Let $m\ge3$, and let $l\in \HH^0(X,\OO_X(H))$ have a smooth
connected zero divisor $X'$. Put $H'=H|_{X'}$. Restriction gives
an injective graded map $R(X,H)/lR(X,H)\to R(X',H')$ with
finite-dimensional cokernel, and a finite injective local map
\[
 B(X,H)/lB(X,H)\longrightarrow B(X',H')
\]
with finite-length cokernel. The latter induces an isomorphism of
punctured spectra.
\end{lemma}
\begin{proof}
Write $R=R(X,H)$ and $R'=R(X',H')$. For $t\ge0$, the cohomology
sequences of
\[
 0\to\OO_X((t-1)H)\xrightarrow{l}\OO_X(tH)\to\OO_{X'}(tH')\to0
\]
give
\[
 0\longrightarrow R(-1)\xrightarrow{l}R\longrightarrow R'
 \longrightarrow\Gamma\longrightarrow0,\qquad
 \Gamma\subset\bigoplus_{t\ge0}\HH^1(X,\OO_X((t-1)H)).
\]
Serre vanishing makes $\Gamma$ finite-dimensional. Thus $R'$ is finite
over $R$, and $R'/\mm R'$ is finite-dimensional, so $\mm R'$
is primary to $R'_+$. Hence $B(X',H')=R'\otimes_R B(X,H)$.
Tensoring the sequence with the flat completion, and using
$\Gamma\otimes_R B(X,H)=\Gamma$, proves the local assertion. Finally, an
injective finite local map with finite-length cokernel becomes an
isomorphism after localizing at every nonmaximal prime, and the
maximal ideal of its target is the only prime over the maximal
ideal of its source. Hence the punctured spectra are isomorphic.
\end{proof}

By Bertini \cite[Corollary~III.10.9]{Ha}, general
$l_1,\ldots,l_{m-2}\in \HH^0(X,\OO_X(H))$ cut out a chain of smooth closed subschemes
\[
 X=X_0\supset X_1\supset\cdots\supset X_{m-2}=C,\qquad
 X_i=X_{i-1}\cap\{l_i=0\},\qquad \dim X_i=m-i.
\]
They are connected, hence irreducible since they are smooth: on a
connected $X_{i-1}$ of dimension at
least two, Kodaira vanishing and Serre duality give
$\HH^1(X_{i-1},\OO_{X_{i-1}}(-H))=0$; the divisor exact sequence
then gives $h^0(\OO_{X_i})=1$. Call $C$ a \emph{surface section}
and put $H_C=H|_C$. If $m=2$, take $C=X$.

\begin{proposition}[Restriction to a surface section]
\label{hd:restriction}
If $N$ is a finite $B(X,H)$-module with $\HH^j_{\mm}(N)$ of finite
length for every $j\le m$, then there is a finite reflexive
$B(C,H_C)$-module $V$ such that
\begin{equation}\label{hd:restriction-bound}
 \rk V=\rk N,\qquad
 \length \HH^2_{\mm}(V)\le
 \sum_{j=0}^{m-2}\binom{m-2}{j}\length \HH^{2+j}_{\mm}(N).
\end{equation}
\end{proposition}
\begin{proof}
Put $B_i=B(X_i,H|_{X_i})$, normal complete local domains regular
away from their closed points by Proposition~\ref{hd:cone-model}.
Write $\HH^a(-)$ for local cohomology at the closed point. By
Lemma~\ref{hd:hyperplane}, $B_i$ is finite and local over
$B_{i-1}/l_iB_{i-1}$. Local cohomology and finite lengths therefore
agree over $B_{i-1}$, $B_{i-1}/l_iB_{i-1}$ and $B_i$, all with
residue field $\CC$.

We construct finite $B_i$-modules $F_i$ of rank $r=\rk N$, with
$\HH^0(F_i)=0$, locally free on the punctured spectrum, and satisfying
\begin{equation}\label{hd:inductive-restriction}
 \length \HH^a(F_i)\le\sum_{j=0}^i\binom{i}{j}
 \length \HH^{a+j}(N)\qquad(a\ge1).
\end{equation}
The right side is finite when $a+i\le m$. Take
$F_0=N/\HH^0(N)$; its positive local cohomology agrees with that of
$N$, and Lemma~\ref{hd:punctured-freeness}(b) gives the punctured
local freeness.

Suppose $F=F_{i-1}$ has been constructed. Its only possible
associated prime is zero: the maximal ideal is excluded by
$\HH^0(F)=0$, and other nonzero primes are excluded by local
freeness over a domain. The image $l$ of $l_i$ is nonzero and
therefore is $F$-regular. The module $\overline F=F/lF$ is locally
free of rank $r$ on the punctured spectrum of $B_{i-1}/lB_{i-1}$,
and the local cohomology sequence gives
\[
 \length \HH^a(\overline F)\le
 \length \HH^a(F)+\length \HH^{a+1}(F)\qquad(a\ge0).
\]
Set $G=\overline F\otimes_{B_{i-1}}B_i$ and
$F_i=G/\HH^0(G)$. If $\Gamma_i$ is the finite-length cokernel of
$B_{i-1}/lB_{i-1}\to B_i$, then the kernel and cokernel of
$\overline F\to G$ have finite length: they are a quotient of
$\Tor_1^{B_{i-1}/lB_{i-1}}(\overline F,\Gamma_i)$ and the module
$\overline F\otimes\Gamma_i$, respectively. The map is an isomorphism
on punctured spectra. Thus $F_i$ is locally free there of rank
$r$. Let $I$ be the image of $\overline F\to G$. Then
$\HH^a(\overline F)\simeq \HH^a(I)$ for $a\ge1$,
$\HH^1(I)\to \HH^1(G)$ is surjective, and
$\HH^a(I)\simeq \HH^a(G)$ for $a\ge2$. Consequently
\[
 \length \HH^a(F_i)=\length \HH^a(G)
 \le\length \HH^a(\overline F)
 \le\length \HH^a(F_{i-1})+\length \HH^{a+1}(F_{i-1})
 \quad(a\ge1).
\]
Pascal's rule establishes \eqref{hd:inductive-restriction}.

The final module $Q=F_{m-2}$ is torsion-free over $B(C,H_C)$:
its only associated prime is zero when $r>0$, and it is zero
when $r=0$. Set $V=Q^{**}$. Since $Q$ is free away from the
closed point, $V/Q$ has finite length. Thus $\rk V=r$ and
$\HH^2(V)\simeq \HH^2(Q)$. Apply
\eqref{hd:inductive-restriction} with $i=m-2$, $a=2$.
\end{proof}

\subsection{The higher-dimensional criterion}
\begin{theorem}[Higher-dimensional criterion]\label{hd:criterion}
For a smooth connected complex projective variety $X$ of dimension
$m\ge2$ and an ample globally generated $H$, put
\[
 d=H^m,\qquad
 \sigma_{X,H}=\frac{(K_X^2-2c_2(X))\cdot H^{m-2}-(m-2)d}{3},
 \qquad \lambda_{X,H}=\frac{K_X\cdot H^{m-1}}{d}+m-2.
\]
Suppose an integer $q\ge1$ satisfies
\begin{equation}\label{hd:condition}
 \bigl(q^2d<\sigma_{X,H}\ \text{and}\ \lambda_{X,H}<5q\bigr)
 \quad\text{or}\quad
 \bigl(q^2d+16<\sigma_{X,H}\ \text{and}\ \lambda_{X,H}\le5q\bigr).
\end{equation}
Let $c_q=q(\sigma_{X,H}-q^2d)/8$ if $\lambda_{X,H}<5q$, and let
$c_q=q((\sigma_{X,H}-q^2d)/8-2)$ if $\lambda_{X,H}=5q$. Write
$B=B(X,H)$.
\begin{enumerate}[label=\textup{(\alph*)}]
\item Every finite $B$-module $N$ with finite-length
$\HH^j_{\mm}(N)$ for $j\le m$ satisfies
\begin{equation}\label{hd:weighted-gap}
 \sum_{j=0}^{m-2}\binom{m-2}{j}\length \HH^{2+j}_{\mm}(N)
 \ge c_q\rk_B N.
\end{equation}
\item For each system of parameters $\mathbf x$ and every
nonzero finite $B$-module $M$,
\[
 \kappa_{\mathbf x}(M)\ge\eta_{\mathbf x}\nu_B(M),\qquad
 \eta_{\mathbf x}=\frac{c_q}
 {(m+1)(1+b)^m\theta_m e h_B+c_q},\qquad
 \theta_m=\binom{m-2}{\lfloor(m-2)/2\rfloor}.
\]
Here $e$, $b$ and $h_B$ are the invariants of $B$ and $\mathbf x$
defined at the beginning of Section~\ref{hd:section}; $e=d$
for parameters given by $m+1$ sections of $H$ without common zero.
\item Neither $B(X,H)$ nor $R(X,H)_{\mm}$ admits an ordinary lim
Cohen--Macaulay sequence or a small Cohen--Macaulay module.
The completed ring is a normal local domain of dimension $m+1$
with residue field $\CC$ and an isolated singularity.
\end{enumerate}
\end{theorem}
\begin{proof}
Let $C$ be a surface section. Its normal bundle is
$\OO_C(H_C)^{\oplus(m-2)}$, so adjunction and the Chern class
identity $c(T_C)=c(T_X)|_C(1+H_C)^{-(m-2)}$ give
\[
 K_C=(K_X+(m-2)H)|_C,\qquad
 c_2(C)=\left(c_2(X)+(m-2)K_X\cdot H+
 \binom{m-1}{2}H^2\right)\Big|_C.
\]
Hence $H_C^2=d$, $K_C\cdot H_C/d=\lambda_{X,H}$, and
\[
 \begin{split}
 3\sigma(C)&=K_C^2-2c_2(C)\\
 &=(K_X^2-2c_2(X))\cdot H^{m-2}
   +\bigl((m-2)^2-(m-1)(m-2)\bigr)d
 =3\sigma_{X,H}.
 \end{split}
\]
Thus \eqref{hd:condition} is condition \eqref{main:local-condition}
of Theorem~\ref{main:criteria} for $(C,H_C)$, and $c_q$ is the
constant of Theorem~\ref{main:local} for $B(C,H_C)$.
For (a), use the reflexive module $V$ of
Proposition~\ref{hd:restriction} and Theorem~\ref{main:local}:
$\length \HH^2_{\mm}(V)\ge c_q\rk V=c_q\rk N$.
Then \eqref{hd:restriction-bound} proves \eqref{hd:weighted-gap}.

For (b), apply Proposition~\ref{hd:replacement} to obtain $M'$.
Part (a) and \eqref{hd:replacement-bound} yield
\[
 c_q\rk M\le c_q\rk M'\le
 \theta_m\sum_{j=2}^m\length \HH^j_{\mm}(M')
 \le(m+1)(1+b)^m\theta_mh_B\kappa_{\mathbf x}(M).
\]
In a minimal $A$-free resolution,
\[
 \beta_0^A(M)-\beta_1^A(M)
 =e\rk_B M-\sum_{i\ge2}(-1)^i\beta_i^A(M).
\]
The last sum is the nonnegative rank of the second syzygy. Thus
\[
 \nu_B(M)\le\beta_0^A(M)\le e\rk_B M+\beta_1^A(M)
 \le\left(\frac{(m+1)(1+b)^m\theta_m e h_B}{c_q}+1\right)
 \kappa_{\mathbf x}(M),
\]
which is the desired estimate.

For (c), part (b) excludes an ordinary lim Cohen--Macaulay
sequence over $B$;
completion preserves finite Koszul homology and numbers of
generators, so it excludes one before completion as well. A small
Cohen--Macaulay module would give a constant sequence.
Proposition~\ref{hd:cone-model} gives regularity away from the
vertex. The completed ring cannot be Cohen--Macaulay, since it
would itself be a small Cohen--Macaulay module. Hence its closed
point is an isolated singular point.
\end{proof}

For $m=2$, Theorem~\ref{hd:criterion} recovers the local assertion
of Theorem~\ref{main:criteria}, with a smaller constant than
$\eta_{\mathbf x}$ there. The strict numerical hypothesis can also be
written
\[
 (K_X^2-2c_2(X))\cdot H^{m-2}>(m-2+3q^2)H^m,
 \qquad K_X\cdot H^{m-1}<(5q-m+2)H^m.
\]
By \cite[Theorem~10.3]{An3}, absence of small Cohen--Macaulay
modules already follows from
$(K_X^2-2c_2(X))\cdot H^{m-2}>(m+4)H^m$, equivalently
$\sigma_{X,H}>2d$, without a condition on $\lambda_{X,H}$.

\subsection{Products with projective spaces}
Let $Y$ be a smooth connected complex projective surface and $H_Y$
an ample globally generated divisor, with invariants
$d_Y=H_Y^2$, $\sigma_Y=\sigma(Y)$ and
$\lambda_Y=K_Y\cdot H_Y/d_Y$. For $m\ge2$, take
\[
 X=Y\times\mathbf P^{m-2},\qquad
 H=\operatorname{pr}_1^*H_Y+
 \operatorname{pr}_2^*\OO_{\mathbf P^{m-2}}(1).
\]
For $m=2$ this is $(Y,H_Y)$. The divisor $H$ is globally generated,
as a sum of pullbacks of globally generated divisors, and ample,
since a multiple of it is very ample through the Segre embedding.

\begin{corollary}[Products]\label{hd:products}
For this pair,
\begin{equation}\label{hd:product-invariants}
 d=\binom m2d_Y,\qquad
 \sigma_{X,H}=\sigma_Y-\binom{m-1}{2}d_Y,\qquad
 \lambda_{X,H}=\frac{2\lambda_Y+m-2}{m}.
\end{equation}
If an integer $q\ge1$ satisfies
\begin{equation}\label{hd:product-condition}
 \sigma_Y>\left(\binom{m-1}{2}+q^2\binom m2\right)d_Y,
 \qquad 2\lambda_Y+m-2<5mq,
\end{equation}
then $B(X,H)$ is a normal complete local domain of dimension
$m+1$ with residue field $\CC$ and an isolated singularity, and
neither it nor $R(X,H)_{\mm}$ admits an ordinary lim Cohen--Macaulay
sequence or a small Cohen--Macaulay module. The weighted
inequality \eqref{hd:weighted-gap} holds with
$c_q=q(\sigma_{X,H}-q^2d)/8$.
\end{corollary}
\begin{proof}
Let $\gamma$ be the pulled-back hyperplane class from
$\mathbf P^{m-2}$, and also write $H_Y,K_Y,c_2(Y)$ for pullbacks.
Then $H=H_Y+\gamma$, $\gamma^{m-1}=0$, and intersection with
$\gamma^{m-2}$ preserves the degree of a codimension-two class
on $Y$. Thus $H^m=\binom m2H_Y^2\gamma^{m-2}$.
From $c(T_X)=c(T_Y)(1+\gamma)^{m-1}$ one gets
\[
 K_X=K_Y-(m-1)\gamma,\qquad
 c_2(X)=c_2(Y)-(m-1)K_Y\cdot\gamma+
 \binom{m-1}{2}\gamma^2,
\]
and hence
\[
 K_X^2-2c_2(X)=K_Y^2-2c_2(Y)+(m-1)\gamma^2,
\]
\[
 (K_X^2-2c_2(X))\cdot H^{m-2}
 =3\sigma_Y+(m-1)\binom{m-2}{2}d_Y.
\]
Use
$(m-1)\binom{m-2}{2}-(m-2)\binom m2
=-3\binom{m-1}{2}$ to obtain the expression for $\sigma_{X,H}$;
compare \cite[Corollary~10.9]{An3}. Finally,
\[
 K_X\cdot H^{m-1}=(m-1)K_Y\cdot H_Y
 -(m-1)\binom{m-1}{2}d_Y,
\]
so
\[
 \lambda_{X,H}=\frac{2\lambda_Y}{m}
 -\frac{(m-1)(m-2)}m+m-2
 =\frac{2\lambda_Y+m-2}{m}.
\]
For $m=2$ the same identities hold with $\gamma=0$ and the
binomial coefficients $\binom{0}{2}=\binom{1}{2}=0$.
Condition~\eqref{hd:product-condition} is the strict case of
\eqref{hd:condition}, so Theorem~\ref{hd:criterion} applies.
\end{proof}

\begin{corollary}[Dual Hesse products]\label{hd:dual-hesse}
For $m\ge2$ let $Y_{32m}$ be the dual Hesse surface $Y_n$ of
Section~\ref{ex:families} with $n=32m$, and put
\[
  Z_m=Y_{32m}\times\mathbf P^{m-2},\qquad
 H=\operatorname{pr}_1^*(6A_{32m}-E_{32m})+
 \operatorname{pr}_2^*\OO_{\mathbf P^{m-2}}(1).
\]
Condition~\eqref{hd:product-condition} holds with $q=13$.
Consequently $B(Z_m,H)$ is a normal complete local domain of
dimension $m+1$, residue field $\CC$, and an isolated singularity.
Neither it nor $R(Z_m,H)_{\mm}$ admits an ordinary lim
Cohen--Macaulay sequence or a small Cohen--Macaulay module.
The weighted inequality \eqref{hd:weighted-gap} holds with
\[
 c_{13}=\frac{13}{8}(8m^2+2064m-37)(32m)^6.
\]
With $D=m+1$, this proves Theorem~\ref{main:dimensions}(a).
\end{corollary}
\begin{proof}
Table~\ref{ex:family-table} gives
$d_Y=24n^6$, $\sigma_Y=(2n^2-13)n^6$ and
$\lambda_Y=n-5/4$. With $n=32m$ and $q=13$,
\[
 \lambda_{X,H}=\frac{64m-5/2+m-2}{m}
 =65-\frac9{2m}<65=5q,
\]
while
\[
 \begin{split}
 \sigma_{X,H}-q^2d
 &=\bigl(2n^2-13-12(m-1)(m-2)-2028m(m-1)\bigr)n^6\\
 &=(8m^2+2064m-37)(32m)^6>0.
 \end{split}
\]
The divisor $H$ is ample and globally generated, as noted before
Corollary~\ref{hd:products}. Apply Corollary~\ref{hd:products} and
Theorem~\ref{hd:criterion}.
\end{proof}
For $n=32m$, $q=13$ is the least admissible value, since
$\lambda_{X,H}/5>12$. Among choices $n=km$ with one $q$ for all $m$,
this choice has the least $q$ and the least $k$. Indeed $\lambda_{X,H}=2k+1-9/(2m)$ is
less than $5q$ for all $m$ only if $2k+1\le5q$, and the coefficient
$2k^2-12q^2-12$ of $m^2$ in $(\sigma_{X,H}-q^2d)/n^6$ must be
nonnegative, hence positive, since $6(q^2+1)$ is never a square. The
least $q$ allowing both is $13$, which forces $k=32$; no $k<32$ works
for any $q$. For a single $m$, much
smaller $n$ suffice (Example~\ref{hd:dimension-four}).

\begin{example}[Dimension four]\label{hd:dimension-four}
Take $m=3$ and the dual Hesse surface $Y_8$ with
$H_Y=6A_8-E_8$. Its invariants are
$d_Y=24\cdot8^6$, $\sigma_Y=115\cdot8^6$ and
$\lambda_Y=27/4$. For $X=Y_8\times\mathbf P^1$,
\eqref{hd:product-invariants} gives
\[
 d=72\cdot8^6,\qquad \sigma_{X,H}=91\cdot8^6,
 \qquad \lambda_{X,H}=29/6<5.
\]
Thus $q=1$ applies. The four-dimensional completed cone has an
isolated singularity and no ordinary lim Cohen--Macaulay sequence.
Every finite module $N$ with $\HH^j_{\mm}(N)$ of finite length for
$j\le3$ satisfies
\[
 \length \HH^2_{\mm}(N)+\length \HH^3_{\mm}(N)
 \ge19\cdot8^5\rk N.
\]
Here $\sigma_{X,H}<2d$, so \cite[Theorem~10.3]{An3} does not give
the small Cohen--Macaulay obstruction; it follows from
Theorem~\ref{hd:criterion}.

The endpoint case of \eqref{hd:condition} gives an example of
smaller degree. Take the quadrangle
surface
$Y_9$ with $H_Y=3A_9-E_9$. Then
$d_Y=5\cdot9^3$, $\sigma_Y=(71/3)9^3$, $\lambda_Y=7$.
For $Y_9\times\mathbf P^1$, the surface-section invariants are
\[
 d=10935,\qquad \sigma_{X,H}=13608,\qquad \lambda_{X,H}=5.
\]
Since $d+16<\sigma_{X,H}$, the endpoint of
Theorem~\ref{hd:criterion} applies with $q=1$ and constant
$(\sigma_{X,H}-d)/8-2=2657/8$. Thus this four-dimensional completed cone
also has an isolated singularity, no ordinary lim Cohen--Macaulay
sequence and no small Cohen--Macaulay module. Again $\sigma_{X,H}<2d$, so
\cite[Theorem~10.3]{An3} does not apply.
\end{example}

\begin{remark}[Scope]\label{hd:scope}
\begin{enumerate}
\item The strict case of \eqref{hd:condition} is equivalent to
$\lambda_{X,H}/5<q<\sqrt{\sigma_{X,H}/d}$. The second alternative adds only
the endpoint $q=\lambda_{X,H}/5$, when this is an integer, and then
requires $q^2d+16<\sigma_{X,H}$. An admissible $q$ gives $d<\sigma_{X,H}$, so
Lemma~\ref{ex:range}, applied to a surface section, gives
$\lambda_{X,H}>0$; squaring $\lambda_{X,H}/5\le q<\sqrt{\sigma_{X,H}/d}$ shows
that a necessary condition is
$\lambda_{X,H}^2d/\sigma_{X,H}<25$. If this ratio stays below a fixed
constant less than $25$ and $\lambda_{X,H}$ is large, the strict
interval contains an integer. Fix $m$ and let $n\to\infty$ in the products of the three standard
families with $\mathbf P^{m-2}$; by \eqref{hd:product-invariants},
\[
 \lim_{n\to\infty}\frac{\lambda_{X,H}^2d}{\sigma_{X,H}}
 =\begin{cases}
 24(m-1)/m,&\text{dual Hesse},\\
 216(m-1)/(7m),&\text{Hesse},\\
 30(m-1)/m,&\text{quadrangle}.
 \end{cases}
\]
The first is less than $25$ for every $m$. The other two are less
than $25$ only for $m\le5$. Thus for large $n$ the Hesse and
quadrangle products satisfy the criterion when $m\le5$, that is, in
dimension at most six, and fail it when $m\ge6$. At the borderline quadrangle value $m=6$,
the ratio is exactly $25n^2/(n^2-160)>25$, although its limit is
$25$.
\item This section makes no claim for weakly lim Cohen--Macaulay
sequences. Proposition~\ref{hd:replacement} uses all positive Betti
numbers, not their alternating sum. For $m\ge3$, depth at least
two no longer forces projective dimension at most one over $A$, so
the depth-two Euler gap \cite[Proposition~2.2]{AnII}, which uses
$\operatorname{pd}_AM\le1$, and the arguments of \cite{AnII} built on
it do not extend directly.
\end{enumerate}
\end{remark}

\section{Arithmetic transfer and formal variables}\label{arith:section}

The arithmetic construction follows \cite[Remark~3.3]{J}, where an
example defined over $\mathbf Q$ yields, for every prime $p$, a
four-dimensional example with residue field $\mathbf F_p$. We make
three changes. The pair may be defined over a number field, as our
construction of the Hesse families requires
(Section~\ref{ex:families}). The residue field is the algebraically
closed field $\overline{\mathbf F}_p$, as in \cite[Conjectures~5.2 and~13.2]{Ho}
and \cite[Corollary~6.3]{BHM}. Formal variables give every dimension
at least four (Proposition~\ref{arith:mixed-variables}). The transfer
also applies to the weak condition (Lemma~\ref{arith:flat-transfer};
see Remark~\ref{arith:scope}).

Given a pair $(Y_F,H_F)$ over a number field $F$, we construct a
finitely generated graded model $\mathcal R$ of $R_F\otimes_FK$ over a
complete discrete valuation ring $V$ with residue field
$\overline{\mathbf F}_p$ and fraction field $K$
(Lemmas~\ref{arith:coefficients} and~\ref{arith:model}), and complete
it at the vertex of its special fibre. By
Lemma~\ref{arith:flat-transfer}, an ordinary lim Cohen--Macaulay
sequence over the resulting ring $T_p$ would give one over the
completed generic vertex $\widehat R_K$, and flat base change along an
embedding $\theta\colon K\hookrightarrow\CC$ then gives one over the
complex cone $\widehat R_{\CC}$, which is excluded by the hypothesis
of Theorem~\ref{arith:transfer}.

\subsection{A flat transfer lemma}
For a regular local ring $\Lambda$ with residue field $\Bbbk$ and a finite
$\Lambda$-module $M$, write
$\beta_i^\Lambda(M)=\dim_{\Bbbk}\Tor_i^\Lambda(\Bbbk,M)$.
The following lemma controls positive Koszul homology and the first
partial Euler characteristic under a flat change of the regular
subring.

\begin{lemma}[Flat base change of the regular subring]\label{arith:flat-transfer}
Let $\Lambda'$ be a regular local ring and $T'$ a local domain which is
a finite $\Lambda'$-algebra of rank $e>0$, with the same residue field.
Let $\Lambda'\to\Lambda''$ be a flat ring homomorphism to a regular
local ring $\Lambda''$, not necessarily a local homomorphism. Suppose
that $T''=T'\otimes_{\Lambda'}\Lambda''$ is a local domain with the
same residue field as $\Lambda''$. For a finite $T'$-module $M$, put
$N=M\otimes_{\Lambda'}\Lambda''$. Then
\begin{equation}\label{arith:betti-transfer}
 \rk_{T''}N=\rk_{T'}M,\qquad
 \beta_i^{\Lambda''}(N)\le\beta_i^{\Lambda'}(M)\quad(i\ge0).
\end{equation}
Consequently, an ordinary lim Cohen--Macaulay sequence over $T'$ gives
one over $T''$, and the same implication holds for weakly lim
Cohen--Macaulay sequences.
\end{lemma}
\begin{proof}
The map $\Lambda'\to T'$ is injective, since $T'$ has positive rank
over the domain $\Lambda'$. By flatness,
$\Lambda''\to T''$ is also injective and finite. Let $G_\bullet\to M$
be a minimal free resolution over $\Lambda'$. It is finite, and
$\rk G_i=\beta_i^{\Lambda'}(M)$. Flat base change gives a free
resolution of $N$ over $\Lambda''$, not necessarily minimal, proving
the Betti inequalities. Taking alternating sums of ranks gives
\[
 \rk_{\Lambda''}N=\sum_i(-1)^i\rk G_i=\rk_{\Lambda'}M.
\]
Applying this first to $M=T'$ shows that $\rk_{\Lambda''}T''=e$.
Since $T'$ and $T''$ are domains,
$\rk_{\Lambda'}M=e\rk_{T'}M$ and
$\rk_{\Lambda''}N=e\rk_{T''}N$, proving the rank equality.

Let $\mathbf y'$ and $\mathbf y''$ be regular systems of parameters
of $\Lambda'$ and $\Lambda''$; they are also systems of parameters of
$T'$ and $T''$. The agreement of residue fields gives
\[
 \length \HH_i(\mathbf y';M)=\beta_i^{\Lambda'}(M),\qquad
 \chi_1(\mathbf y';M)=\beta_0^{\Lambda'}(M)-e\rk_{T'}M,
\]
and the analogous identities for $N$. In particular,
\begin{equation}\label{arith:generator-rank}
 \nu_{T'}(M)\le\beta_0^{\Lambda'}(M)
 =e\rk_{T'}M+\chi_1(\mathbf y';M),\qquad
 0\le\chi_1(\mathbf y';M)\le\sum_{i\ge1}\beta_i^{\Lambda'}(M).
\end{equation}
For an ordinary lim Cohen--Macaulay sequence $(M_j)$, set
$r_j=\rk_{T'}M_j$ and
$\varepsilon_j=\sum_{i\ge1}\beta_i^{\Lambda'}(M_j)$.
Parameter independence over $T'$ \cite[Theorem~4.6]{BHM} gives
$\varepsilon_j=o(\nu_{T'}(M_j))$. Since
$\nu_{T'}(M_j)\le er_j+\varepsilon_j$, the ranks are positive
eventually and $\varepsilon_j=o(r_j)$. The transferred modules $N_j$
have positive rank $r_j$, hence full dimension and
$\nu_{T''}(N_j)\ge r_j$. By \eqref{arith:betti-transfer}, all their
positive Koszul homology lengths are at most $\varepsilon_j$. By
parameter independence over $T''$, they form an ordinary lim
Cohen--Macaulay sequence after discarding finitely many initial
terms.

For a weakly lim Cohen--Macaulay sequence, parameter independence
\cite[Proposition~2.6]{Ma} gives
$\chi_1(\mathbf y';M_j)=o(\nu_{T'}(M_j))$.
Equation~\eqref{arith:generator-rank} again implies $r_j>0$ eventually
and $\chi_1(\mathbf y';M_j)=o(r_j)$. Moreover,
\[
 0\le\chi_1(\mathbf y'';N_j)
 =\beta_0^{\Lambda''}(N_j)-er_j
 \le\beta_0^{\Lambda'}(M_j)-er_j=o(r_j).
\]
Thus the transferred sequence is weakly lim Cohen--Macaulay.
\end{proof}

\subsection{The arithmetic model}
\begin{theorem}[Arithmetic transfer]\label{arith:transfer}
Let $F\subset\CC$ be a number field, $Y_F$ a smooth geometrically
integral projective surface over $F$, and $H_F$ an ample globally
generated divisor. Put
\[
 R_F=\bigoplus_{t\ge0}\HH^0(Y_F,\OO_{Y_F}(tH_F)),\qquad
 R_{\CC}=R_F\otimes_F\CC,
\]
and let $\widehat R_{\CC}$ be the completion of $R_{\CC}$ at its
irrelevant ideal. If $\widehat R_{\CC}$ admits no ordinary lim
Cohen--Macaulay sequence, then for every prime $p$ there is a
four-dimensional complete Noetherian normal local domain $T_p$ of
mixed characteristic $(0,p)$, with residue field
$\overline{\mathbf F}_p$, admitting no such sequence.
The ring $T_p$ is finite over $V[[s_0,s_1,s_2]]$, where $V$ is a
complete discrete valuation ring with residue field
$\overline{\mathbf F}_p$.
If $\widehat R_{\CC}$ admits no weakly lim Cohen--Macaulay sequence,
then neither does this same ring $T_p$.
\end{theorem}

\begin{lemma}[A coefficient ring and a complex embedding]\label{arith:coefficients}
For every prime $p$ there are a complete discrete valuation ring $V$
of mixed characteristic $(0,p)$ with residue field
$\Bbbk=\overline{\mathbf F}_p$, a field embedding
$F\hookrightarrow K=\operatorname{Frac}(V)$, and a field embedding
$\theta\colon K\hookrightarrow\CC$ extending the given inclusion
$F\subset\CC$.
\end{lemma}
\begin{proof}
Let $V_0=W(\Bbbk)$ and $K_0=\operatorname{Frac}(V_0)$. Choose a place
$v$ of $F$ above $p$ and embed the finite extension $F_v/\mathbf Q_p$
into an algebraic closure of $K_0$. Put $K=K_0F_v$, and let $V$ be
the integral closure of $V_0$ in $K$. This is a finite, possibly
ramified extension of complete discrete valuation rings. Its residue
field is a finite extension of the algebraically closed field $\Bbbk$,
and therefore equals $\Bbbk$. We have $F\subset F_v\subset K$.

The field $\Bbbk$ is countable, so $W(\Bbbk)$ and $K$ have cardinality at
most $2^{\aleph_0}$. A transcendence basis of $K/F$ can therefore be
mapped injectively to an algebraically independent subset of $\CC/F$,
whose transcendence degree is $2^{\aleph_0}$. Extend the given
inclusion of $F$ first to the resulting purely transcendental
extension and then to its algebraic extension $K$, using algebraic
closedness of $\CC$. The embedding $\theta$ need not be continuous
for the valuation topology.
\end{proof}

\begin{lemma}[A finitely generated graded integral model]\label{arith:model}
Let $V,K$ be as above, with uniformizer $\varpi$, and put
$R_K=R_F\otimes_FK$. There is a finitely generated normal graded
$V$-algebra $\mathcal R\subset R_K$ such that
\[
 \mathcal R_0=V,\qquad \mathcal R[1/\varpi]=R_K,\qquad
 \mathcal R_t\text{ is finite free over }V\text{ for every }t.
\]
There are homogeneous elements $s_0,s_1,s_2\in\mathcal R_+$ of
positive degree such that $\mathcal R$ is finite over the polynomial
subring $\Lambda_0=V[s_0,s_1,s_2]$. The positive-degree ideal of
$\mathcal R/(s_0,s_1,s_2)\mathcal R$ is nilpotent.
\end{lemma}
\begin{proof}
Choose finitely many homogeneous generators $w_i$ of positive degree
for $R_F$ over $F$, and let $\mathcal R'\subset R_K$ be the
$V$-subalgebra that they generate.
The ring $R_K$ is the full section ring of a smooth geometrically
integral projective surface, hence is a normal domain, and
$\mathcal R'[1/\varpi]=R_K$. The two rings have the same fraction
field. Since $\mathcal R'$ is of finite type over the excellent ring
$V$, its normalization $\mathcal R$ is finite over $\mathcal R'$.
Normality of $R_K$ places $\mathcal R$ inside $R_K$, and normalization
commutes with localization, giving $\mathcal R[1/\varpi]=R_K$.
Hence $\mathcal R$ is the integral closure of $\mathcal R'$ in the
graded ring $R_K$, and it is graded \cite[Tag~077G]{St}.

If $c\in\mathcal R_0\subset K$, apply the augmentation $R_K\to K$
to an equation of integral dependence over $\mathcal R'$.
It makes $c$ integral over $V$, so $c\in V$. Thus $\mathcal R_0=V$.
Finite homogeneous generation implies that every $\mathcal R_t$ is
finite over $V$; it is torsion-free, being a submodule of $(R_K)_t$,
and therefore free.

Write $\overline{\mathcal R}=\mathcal R/\varpi\mathcal R$.
The equality
\[
 \dim_{\Bbbk}\overline{\mathcal R}_t
 =\rk_V\mathcal R_t=\dim_K(R_K)_t
\]
shows that $\dim\overline{\mathcal R}=\dim R_K=3$. Indeed, choose
$N\ge1$ such that the Veronese subrings $\overline{\mathcal R}^{(N)}$
and $R_K^{(N)}$ are generated by their components of degree $N$. Both
rings are integral over these subrings, which have the same Hilbert
functions; after regrading, their dimension is one more than the
degree of the Hilbert polynomial \cite[Theorem~4.1.3]{BrHe}. Choose a homogeneous system of parameters
$\bar s_0,\bar s_1,\bar s_2$ of $\overline{\mathcal R}$ and homogeneous
lifts $s_i\in\mathcal R$. For
$\Xi=\mathcal R/(s_0,s_1,s_2)\mathcal R$, the quotient
$\Xi/\varpi\Xi$ is finite-dimensional over $\Bbbk$, and each $\Xi_t$ is
finite over $V$. Thus $\Xi_t/\varpi\Xi_t=0$ for $t\gg0$, so
$\Xi_t=0$ by Nakayama's lemma. It follows that $\Xi$ is finite over
$V$ and that $\Xi_+$ is nilpotent.

Lift a finite homogeneous set of $V$-generators of $\Xi$ to
$\mathcal R$. Induction on degree, subtracting a linear combination
of these lifts and writing the remainder as $\sum_i s_ic_i$, shows
that they generate $\mathcal R$ over $V[s_0,s_1,s_2]$.
After inverting $\varpi$, the three-dimensional domain $R_K$ is
finite over $K[s_0,s_1,s_2]$. Consequently the $s_i$ are algebraically
independent over $K$, and hence over $V$.
\end{proof}

\begin{proof}[Proof of Theorem~\ref{arith:transfer}]
Fix $p$ and the preceding objects, and put
\[
 \mathfrak m_{\mathcal R}=(\varpi,\mathcal R_+),\qquad
 \Lambda=V[[s_0,s_1,s_2]],\qquad
 T=\widehat{\mathcal R_{\mathfrak m_{\mathcal R}}}.
\]
The ring $\mathcal R/(\varpi,s_0,s_1,s_2)\mathcal R$ is a
finite-dimensional graded $\Bbbk$-algebra with degree-zero part $\Bbbk$ and
nilpotent positive-degree ideal. Hence $\mathfrak m_{\mathcal R}$ is
the unique prime over the maximal ideal
$\mathfrak n_0=(\varpi,s_0,s_1,s_2)$ of $\Lambda_0$. Since
$\mathcal R$ is finite over $\Lambda_0$, its $\mathfrak n_0$-adic
completion is
\begin{equation}\label{arith:finite-completion}
 T=\mathcal R\otimes_{\Lambda_0}\Lambda.
\end{equation}
Thus $T$ is finite over the four-dimensional regular local ring
$\Lambda$, and $\Lambda\to T$ is injective by flatness of completion.
It is local with residue field $\Bbbk$, and is a normal domain because
$\mathcal R_{\mathfrak m_{\mathcal R}}$ is excellent and normal
\cite[Tag~0C23]{St}. It has dimension four and mixed characteristic
$(0,p)$.

Let $\mathfrak s=(s_0,s_1,s_2)\Lambda$ and $P=\mathcal R_+T$.
Equation~\eqref{arith:finite-completion} gives $T/P=V$, and
the nilpotence in Lemma~\ref{arith:model} gives
$\sqrt{\mathfrak sT}=P$. Hence $P$ is the unique prime over
$\mathfrak s$. Since $\Lambda/\mathfrak s^n=V[\mathbf s]/(\mathbf s)^n$,
localizing at $\mathfrak s$ inverts every nonzero element of $V$;
hence
\[
 \widehat{\Lambda_{\mathfrak s}}=K[[s_0,s_1,s_2]]=:\Lambda_K.
\]
Because $P$ is the only prime of $T$ over $\mathfrak s$, we have
$T\otimes_\Lambda\Lambda_{\mathfrak s}=T_P$, and on $T_P$ the
$\mathfrak s$-adic and $P$-adic topologies agree; hence
\begin{equation}\label{arith:generic-vertex}
 \widehat{T_P}
 =T\otimes_\Lambda\Lambda_K
 =\mathcal R\otimes_{\Lambda_0}\Lambda_K
 =R_K\otimes_{K[s_0,s_1,s_2]}\Lambda_K
 =\widehat R_K.
\end{equation}
The last equality holds because the irrelevant ideal of $R_K$ is the
only prime over $(s_0,s_1,s_2)$. The completed generic vertex $\widehat R_K$
is a complete normal local domain of dimension three with residue
field $K$; the map $T\to\widehat R_K$ is localization followed by
completion, and is flat and injective.

By Lemma~\ref{arith:flat-transfer}, an ordinary lim Cohen--Macaulay
sequence over $T$ would give one over $\widehat R_K$, and the same is true of
a weakly lim Cohen--Macaulay sequence.

Choose $\theta\colon K\hookrightarrow\CC$ as in
Lemma~\ref{arith:coefficients}. The induced local homomorphism
\[
 \Lambda_K=K[[s_0,s_1,s_2]]\longrightarrow
 \Lambda_{\CC}=\CC[[s_0,s_1,s_2]]
\]
is flat: both rings are regular of dimension three and the closed
fibre is the field $\CC$ \cite[Tag~00R4]{St}.
Because $\theta$ agrees with the given embedding on $F$,
\begin{equation}\label{arith:complex-fibre}
 \widehat R_K\otimes_{\Lambda_K}\Lambda_{\CC}
 =R_K\otimes_{K[s_0,s_1,s_2]}\Lambda_{\CC}
 =R_{\CC}\otimes_{\CC[s_0,s_1,s_2]}\Lambda_{\CC}
 =\widehat R_{\CC}.
\end{equation}
Here $s_i\mapsto s_i\otimes1\in R_K\otimes_{K,\theta}\CC=R_{\CC}$, and
completing
$R_{\CC}$ at $(s_0,s_1,s_2)$ is completing at its irrelevant ideal.
By Proposition~\ref{loc:model}, $\widehat R_{\CC}$ is a local domain with residue field $\CC$,
and Lemma~\ref{arith:flat-transfer}, applied to
$\Lambda_K\to\Lambda_{\CC}$ with $T'=\widehat R_K$ and
$T''=\widehat R_{\CC}$, carries either kind of sequence over
$\widehat R_K$ to the same kind over $\widehat R_{\CC}$; compare
\cite[Lemma~2.7]{Ma} and \cite[Proposition~4.12(c)]{BHM}. The
respective nonexistence hypotheses give the required contradictions.
Take $T_p=T$.
\end{proof}

\begin{corollary}[Arithmetic arrangement examples]\label{arith:examples}
For every pair covered by the local assertions of
Corollaries~\ref{ex:standard-cones} and~\ref{ex:other-polarizations},
and every prime $p$, Theorem~\ref{arith:transfer} produces a
four-dimensional complete normal local domain of mixed characteristic
$(0,p)$ with residue field $\overline{\mathbf F}_p$ and no ordinary
lim Cohen--Macaulay sequence.
\end{corollary}
\begin{proof}
By Section~\ref{ex:families}, $Y_n$, $A_n$ and $E_n$, and hence the
polarizations of Corollaries~\ref{ex:standard-cones}
and~\ref{ex:other-polarizations}, are defined over
$F=\mathbf Q(\zeta_3)$ for the Hesse and dual Hesse families and over
$F=\mathbf Q$ for the complete quadrangle. The individual singular
points and exceptional components need not be defined over $F$. Since
$F$ is perfect, forming the reduced singular locus, its blow-up and
the reduced exceptional divisor commutes with the extension
$F\subset\CC$. Smoothness, geometric integrality, ampleness and global
generation descend from $\CC$ to $F$.
Corollaries~\ref{ex:standard-cones} and~\ref{ex:other-polarizations}
supply the nonexistence hypothesis of Theorem~\ref{arith:transfer}.
No prime is excluded, including $3$ and the primes dividing $n$:
the special fibre of $\mathcal R$ may be bad or nonreduced, and the
proof requires only the normality of the total model $\mathcal R$
and its prescribed generic fibre.
\end{proof}

\begin{remark}[Scope and singular locus]\label{arith:scope}
Beyond nonexistence over $\widehat R_{\CC}$, the proof uses only
that $R_F\otimes_FK$ is a finitely generated normal graded domain of
dimension three with degree-zero part $K$. Generation in degree one
and degree-one parameters are unnecessary. Under the hypothesis of
Theorem~\ref{arith:transfer}, $T$ does not have an isolated
singularity. Indeed $T/P=V$, and $\widehat{T_P}=\widehat R_K$ is not
Cohen--Macaulay: otherwise its constant sequence would remain lim
Cohen--Macaulay after the complex base change. Completion preserves
the Cohen--Macaulay property, so $T_P$ is not Cohen--Macaulay.
If $T_{P'}$ were Cohen--Macaulay for $P'\supset P$, its localization
$T_P$ would be so as well. Thus the non-Cohen--Macaulay locus contains
the one-dimensional closed subset $\Spec(T/P)$ of $\Spec T$.

The weak assertions of Theorem~\ref{arith:transfer} and
Proposition~\ref{arith:mixed-variables} are conditional: weak
nonexistence over $\widehat R_{\CC}$ is not known for any of our
examples, and the partial results of \cite{AnII} do not supply it.
\end{remark}

\subsection{Formal variables in mixed characteristic}
\begin{proposition}[Formal variables]\label{arith:mixed-variables}
Let $T=T_p$ be the ring constructed in the proof of
Theorem~\ref{arith:transfer}, let $g\ge0$ and
$\mathbf u=(u_1,\ldots,u_g)$. The ring $T[[\mathbf u]]$
is a complete
Noetherian normal local domain of dimension $4+g$ and mixed
characteristic $(0,p)$, with residue field
$\overline{\mathbf F}_p$, finite over $V[[s_0,s_1,s_2,\mathbf u]]$.
If $\widehat R_{\CC}$ admits no ordinary lim Cohen--Macaulay sequence,
then neither does $T[[\mathbf u]]$. If $\widehat R_{\CC}$ admits no weakly lim
Cohen--Macaulay sequence, then neither does $T[[\mathbf u]]$.
\end{proposition}
\begin{proof}
Put $\Lambda'=\Lambda[[\mathbf u]]$ and
\[
 T'=T[[\mathbf u]]
 =T\otimes_\Lambda\Lambda'
 =\mathcal R\otimes_{\Lambda_0}\Lambda'.
\]
These equalities use finiteness over the indicated regular subrings.
The normal domain $\mathcal R[\mathbf u]$ is of finite type over $V$
and finite over $\Lambda_0[\mathbf u]$. Its only prime over
$(\varpi,s_0,s_1,s_2,\mathbf u)$ is
$(\mathfrak m_{\mathcal R},\mathbf u)$. Its completion there is
$T'$, a complete normal local domain by \cite[Tag~0C23]{St}, of
dimension $4+g$ and residue field $\Bbbk$.

Let $V'=V[[\mathbf u]]$, $L=\operatorname{Frac}(V')$ and
$\mathfrak s'=(s_0,s_1,s_2)\Lambda'$. Localizing the truncated
polynomial rings over $V'$ gives
\[
 \widehat{\Lambda'_{\mathfrak s'}}=L[[s_0,s_1,s_2]].
\]
The field $L$ contains $K\supset F$ and has cardinality at most
$2^{\aleph_0}$. As in Lemma~\ref{arith:coefficients}, there is an
embedding $\theta':L\hookrightarrow\CC$ extending the given
inclusion of $F$. The induced map
$L[[s_0,s_1,s_2]]\to\Lambda_{\CC}$ is flat by
\cite[Tag~00R4]{St}; hence so is the composite
\[
 \Lambda'\longrightarrow\Lambda'_{\mathfrak s'}
 \longrightarrow L[[s_0,s_1,s_2]]\longrightarrow\Lambda_{\CC}.
\]
The image of $\varpi$ is a unit, and $V\to\CC$ factors through
$\theta'|_K$. As in \eqref{arith:complex-fibre},
\[
 T'\otimes_{\Lambda'}\Lambda_{\CC}
 =\mathcal R\otimes_{\Lambda_0}\Lambda_{\CC}
 =R_K\otimes_{K[s_0,s_1,s_2]}\Lambda_{\CC}
 =\widehat R_{\CC}.
\]
This is a local domain with residue field $\CC$.
Lemma~\ref{arith:flat-transfer} carries any hypothetical sequence,
ordinary or weak, to one of the same kind over the given complex
cone. The corresponding nonexistence hypothesis proves the assertion.
\end{proof}

The argument uses flat transfer, which treats ordinary and weak
sequences alike. For ordinary sequences one cannot instead reduce
modulo the formal variables, since these need not form a regular
sequence on the modules of a hypothetical sequence.

Under the hypothesis on $\widehat R_{\CC}$, the ring
$T[[\mathbf u]]$ does not have an isolated singularity. Indeed,
$Q=(P,\mathbf u)T[[\mathbf u]]$ has quotient $V$.
The flat local map $T_P\to T[[\mathbf u]]_Q$ has regular closed
fibre, a localization of $V[[\mathbf u]]\otimes_VK$, and its target
cannot be Cohen--Macaulay because that property descends along a
flat local map.

Applying Proposition~\ref{arith:mixed-variables} with $g=D-4$ to
the rings of Corollary~\ref{arith:examples} proves
Theorem~\ref{main:dimensions}(b).

\section{Scope and open questions}\label{sec:questions}

\subsection{Consequences of the proved obstructions}
The inequalities in Theorem~\ref{main:criteria} are uniform in the
module, for each fixed ring and parameter system. They exclude
ordinary lim Cohen--Macaulay sequences and hence small
Cohen--Macaulay modules, strongly lim Cohen--Macaulay sequences and
lim Ulrich sequences \cite[Corollary~10.10 and Section~4.2]{BHM};
the graded criterion excludes the graded versions of these
objects.

The graded criterion requires $d<\sigma$. For arbitrary modules on
the completed cone, the additional slope conditions in
\eqref{main:local-condition} remain necessary for the argument; see
the first of the remaining questions below.

Weakly lim Cohen--Macaulay sequences are not excluded here. The
companion paper \cite{AnII} excludes weak sequences of modules of
depth at least two on the completed cones satisfying
\eqref{main:local-condition}, but its partial results do not cover all
torsion-free modules of depth one.

Big Cohen--Macaulay algebras exist over every Noetherian local ring,
by Hochster--Huneke in equal characteristic and by Andr\'e \cite{And}
in mixed characteristic (see \cite[Section~13]{Ho}), and a lim
Cohen--Macaulay sequence over a complete local ring yields a big
Cohen--Macaulay module \cite[Theorem~9.4]{BHM}; the rings of
Theorem~\ref{main:dimensions} show that the converse fails.

Theorem~\ref{main:dimensions} disproves the conjecture emphasized in
\cite[Section~1]{Ho}, as well as \cite[Conjectures~5.2
and~13.2]{Ho} and \cite[Conjecture~1.1]{IS}. Consequently, the
hypothesis of \cite[Corollary~6.3]{BHM}, existence for all complete
local domains of dimension at most $D$ with perfect residue field,
fails for every $D\ge3$, and for every $D\ge4$ even when restricted
to mixed characteristic. Serre's conjecture itself is
unaffected.

By \cite[Theorem~3.17 and Corollary~3.18]{IS}, a complete Noetherian
local domain $T$ of mixed characteristic with residue field
$\overline{\mathbf F}_p$ has a lim Cohen--Macaulay sequence if it
admits a perfectoid tower arising from $(T,I_0)$ for some ideal $I_0$
\cite[Definition~2.3]{IS}, or a $\delta$-ring structure with $T/pT$
reduced. The rings of Theorem~\ref{main:dimensions}(b) therefore admit
neither.

\subsection{Characteristic and dimension}
The surface inequalities use characteristic zero through Bogomolov's
inequality. Their positive-characteristic counterparts cannot hold in
this generality. Frobenius pushforwards of a line bundle have bounded
$h^1_*$ and unbounded rank \cite[Remark~2.19(3)]{AnII}, in
agreement with the existence of lim Cohen--Macaulay sequences over
$F$-finite rings \cite[Theorem~5.4]{BHM}. The replacement and Koszul-gap
arguments of Section~\ref{transfer:section} and the depth-two gap of
\cite{AnII} work over any coefficient field once the cohomological gap
is known; the trace and Artin-approximation arguments of \cite{AnII} are
written in characteristic zero.

Dimension three is the smallest possible dimension for a complete
normal local domain without ordinary lim Cohen--Macaulay sequences,
since normal local rings of dimension at most two are Cohen--Macaulay.
The complex examples have isolated singularities in every dimension
at least three. The mixed-characteristic examples have dimension at
least four and a positive-dimensional singular locus. Their
construction is a transfer from the vertex of a characteristic-zero
generic fibre, not an obstruction over their positive-characteristic special fibre.

\subsection{Remaining questions}
\begin{enumerate}[label=\textup{(\arabic*)},leftmargin=*]
\item \emph{Removing the slope restriction.}
Does $d<\sigma$ alone imply
$\length \HH^2_{\mathfrak m}(M)\ge(\sigma-d)\rk M/8$ for all finite
modules over $\widehat R$? Theorem~\ref{main:graded} gives this for
graded modules. The two omitted standard cones are concrete local
test cases; both have $2d<\sigma$ and no small Cohen--Macaulay
modules, but whether they admit lim Cohen--Macaulay sequences is
not determined here.

\item \emph{Optimal constants.}
What is the infimum of $h^1_*(E)/\rk E$ over coherent sheaves $E$
of positive rank on a fixed pair with $d<\sigma$,
and how does it compare with the infimum for completed-local modules?
Is the correction $-2$ in Theorem~\ref{main:local}(b) necessary? The stronger
vector-bundle inequality in \cite{AC} is relevant to the graded
question but does not by itself remove the local persistence issue.

\item \emph{Higher-dimensional local inequalities.}
Can one obtain a local cohomology obstruction directly on
$B(X,H)$, without restriction to a surface section?
The small Cohen--Macaulay obstruction in \cite[Section~10]{An3}
requires no surface-slope condition, whereas the ordinary lim
criterion here does. The products of Hesse and quadrangle surfaces
outside that criterion and the weak behavior of the isolated
higher-dimensional examples remain open.

\item \emph{Mixed characteristic.}
Are there complete normal local domains of mixed characteristic
without ordinary lim Cohen--Macaulay sequences which have dimension
three or an isolated singularity? For the constructed models,
nonexistence of weak sequences on the corresponding complex cone
would transfer to mixed characteristic and to their formal-variable
extensions.
\end{enumerate}

\appendix
\section{The graded setting}\label{app:graded}

This appendix proves the graded forms of Lemma~\ref{transfer:base-ring}
and Propositions~\ref{transfer:replacement} and~\ref{transfer:gap},
which give the graded assertion of Theorem~\ref{main:criteria}
(Section~\ref{ex:graded-subsection}). The arguments are those of
Section~\ref{transfer:section}; we write out the points where the
graded category needs care.

\subsection*{Setting}
Let $B=\bigoplus_{t\ge0}B_t$ be a finitely generated graded algebra
over the field $\Bbbk=B_0$ such that $B$ is a normal domain of
dimension three, and put $\mm=B_+$. Let $x_0,x_1,x_2\in B_1$ be such that $B$ is
finite over $A=\Bbbk[x_0,x_1,x_2]$, and put $\nn=A_+$. Since
$\dim B=3$, the $x_i$ are algebraically independent, so $A$ is a
standard graded polynomial ring. The elements of positive degree of
the finite-dimensional graded algebra $B/\nn B$ are nilpotent, so
$\sqrt{\nn B}=\mm$. Hence $\mm$ is the only prime of $B$ over $\nn$,
$B\otimes_AA_{\nn}=B_{\mm}$, and local cohomology with supports in
$\nn$ and in $\mm$ coincide.

All modules are $\ZZ$-graded, and homomorphisms are homogeneous of
degree zero unless stated otherwise; $M(s)_t=M_{s+t}$. For graded
modules, $\Hom(M,N)=\bigoplus_s\Hom_0(M,N(s))$ denotes the module of
homogeneous homomorphisms of all degrees, where $\Hom_0$ consists of
those of degree zero, and $\Ext^i$ are its derived functors; for
finitely generated $M$ these are the usual modules
\cite[Section~1.5]{BrHe}. A graded module of finite length is
finite-dimensional over $\Bbbk$, and its length is its dimension.
Complexes are indexed homologically as in
Section~\ref{transfer:duality}. Shifts $M\mapsto M(s)$ change neither
ranks, nor minimal numbers of generators, nor dimensions over $\Bbbk$.

Define $e$, $a_i$ and $C_{\mathbf x}$ as in \eqref{transfer:constants},
and put $h=\dim_{\Bbbk}\HH^2_{\mm}(B)$. For a finitely generated graded
$B$-module $M$, let $\nu_B(M)=\dim_{\Bbbk}M/\mm M$ be its minimal number
of homogeneous generators, and define $\beta_i^A$ and
$\kappa_{\mathbf x}$ with dimensions over $\Bbbk$ in place of lengths.
The Koszul complex $K(\mathbf x;A)$ is a graded free resolution of
$\Bbbk$, so \eqref{transfer:koszul-betti} holds in this form.

\begin{lemma}[The graded ring $B$ over $A$]\label{app:base-ring}
One has $\operatorname{pd}_AB\le1$, so $a_i=0$ for $i\ge2$ and
$a_0-a_1=e$. Moreover $\Ext_A^1(B,A)$ has dimension $h$ over $\Bbbk$,
$a_1\le h$, and $C_{\mathbf x}\le1+e+2h+(e+h)^2$.
\end{lemma}
\begin{proof}
A minimal graded free resolution of $B$ over $A$ localizes to a
minimal free resolution of $B_{\mm}=B\otimes_AA_{\nn}$ over the
regular local ring $A_{\nn}$. The normal local domain $B_{\mm}$ of
dimension three has depth at least two, also as an $A_{\nn}$-module,
because $\nn B_{\mm}$ is primary to its maximal ideal.
Auslander--Buchsbaum gives $\operatorname{pd}_AB\le1$, and ranks give
$a_0-a_1=e$.

The module $E=\Ext_A^1(B,A)$ is finitely generated and graded, and the
minimal primes of its support are homogeneous. A homogeneous prime
$\mathfrak p\ne\nn$ of $A$ has height at most two, and the argument of
Lemma~\ref{transfer:base-ring} shows that $B_{\mathfrak p}$ is free
over $A_{\mathfrak p}$, so $E_{\mathfrak p}=0$. Hence $E$ has finite
length. Graded local duality \cite[Theorem~3.6.19]{BrHe}, for $A$ with
graded canonical module $A(-3)$, identifies the graded $\Bbbk$-dual of
$\HH^2_{\nn}(B)=\HH^2_{\mm}(B)$ with $\Ext_A^1(B,A(-3))=E(-3)$, so
$\dim_{\Bbbk}E=h$. Dualizing the minimal resolution presents $E$ as the
cokernel of a homogeneous map with entries in $\nn$, so $E$ needs
exactly $a_1$ homogeneous generators. Thus $a_1\le h$ and
$a_0=e+a_1\le e+h$.
\end{proof}

\subsection*{Graded duality}
The polynomial ring $A$ has graded injective dimension three
\cite[Section~3.6]{BrHe}. Let $I^\bullet$ be a resolution of $A(-3)$
by graded-injective $A$-modules, placed in homological degrees
$0,-1,-2,-3$. For a complex $X$ of graded $B$-modules set
\[
 \omega^\bullet=\Hom_A(B,I^\bullet)[3],\qquad
 \mathbb D(X)=\Hom_B(X,\omega^\bullet)=\Hom_A(X,I^\bullet)[3],
\]
the second identification being adjunction. By adjunction the terms
of $\omega^\bullet$ are graded-injective $B$-modules, so $\mathbb D$
preserves quasi-isomorphisms and takes short exact sequences of
complexes to exact triangles in the derived category of graded
$B$-modules. If $X$ has bounded, finitely generated homology, then
over $A$ it is quasi-isomorphic to a bounded complex $G$ of finitely
generated graded free modules, by Hilbert's syzygy theorem, and
$\mathbb D(X)\simeq\Hom_A(G,A(-3))[3]$ over $A$. Hence $\mathbb D(X)$
has bounded, finitely generated homology, and the biduality map
$X\to\mathbb D(\mathbb D(X))$ is a quasi-isomorphism: under adjunction
it is the biduality map of $X$ with respect to $I^\bullet$ over $A$,
which for $X=G$ is that of $G$ with respect to $A(-3)$. The graded
forms of (D1)--(D3) are as follows.
\begin{enumerate}
\item[(D1$^{\mathrm{gr}}$)] For a finitely generated graded module $M$,
$\HH_i(\mathbb D(M))=\Ext_A^{3-i}(M,A(-3))$, vanishing outside
$i=0,1,2,3$. In particular $\omega^\bullet=\mathbb D(B)$ has homology
only in degrees two and three, and $\dim_{\Bbbk}\HH_2(\omega^\bullet)=h$
by Lemma~\ref{app:base-ring}.
\item[(D2$^{\mathrm{gr}}$)] For such $M$,
$\beta_i^A(\mathbb D(M))=\beta_{3-i}^A(M)$.
\item[(D3$^{\mathrm{gr}}$)] For a finitely generated graded module $L$,
$\mathbb D(L[3])\simeq\RHom_A(L,A(-3))$ and
$\HH^i_{\mm}(\mathbb D(L[3]))=0$ for $i\ne3$.
\end{enumerate}
Here $\HH^i_{\mm}$ is computed with the \v{C}ech complex on
$x_0,x_1,x_2$. Assertion (D1$^{\mathrm{gr}}$) is the definition of
$\Ext$. For (D2$^{\mathrm{gr}}$), if $G\to M$ is the minimal graded free
resolution over $A$, then $\mathbb D(M)\simeq\Hom_A(G,A(-3))[3]$ is a
complex of finitely generated graded free $A$-modules whose
differentials have entries in $\nn$ and whose term in degree $i$ has
rank $\beta_{3-i}^A(M)$; tensoring with $\Bbbk$ gives the claim. For
(D3$^{\mathrm{gr}}$), one has
$\mathbb D(L[3])=\Hom_A(L,I^\bullet)$, and if $G\to L$ is a finite
graded free resolution, then $\mathbb D(L[3])\simeq Q=\Hom_A(G,A(-3))$.
The \v{C}ech complex of a finitely generated graded free $A$-module has
cohomology only in degree three. Hence $\HH^i_{\mm}(Q)$ is the homology
at $p=3-i$ of the complex $p\mapsto\HH^3_{\nn}(Q_p)=\Hom_A(G_{-p},E_A)$,
where $E_A=\HH^3_{\nn}(A)(-3)$ is the graded injective hull of $\Bbbk$
\cite[Section~3.6]{BrHe}. As $\Hom_A(-,E_A)$ is exact on graded
modules, this complex has homology only at $p=0$.

\begin{proposition}[Graded reflexive replacement]
\label{app:graded-replacement}
For every finitely generated graded $B$-module $M$ there are a finitely
generated graded reflexive $B$-module $M'$ and a homogeneous
homomorphism $M'\to M$ of degree zero with finite-dimensional cokernel,
such that $\HH^0_{\mm}(M')=\HH^1_{\mm}(M')=0$, $\rk_BM'\ge\rk_BM$ and
\[
 \dim_{\Bbbk}\HH^2_{\mm}(M')
 \le h\left[\beta_1^A(M)+a_0\beta_2^A(M)+(a_1+a_0^2)\beta_3^A(M)\right]
 \le hC_{\mathbf x}\kappa_{\mathbf x}(M).
\]
\end{proposition}
\begin{proof}
We follow the five steps of the proof of
Proposition~\ref{transfer:replacement}.

\emph{Step 1.} By (D1$^{\mathrm{gr}}$), $\mathbb D(M)$ has finitely
generated homology in degrees $0,\ldots,3$. Using the graded Nakayama
lemma, choose degree by degree a quasi-isomorphism $P\to\mathbb D(M)$
from a complex of finitely generated graded free $B$-modules with
$P_i=0$ for $i<0$ and all differentials having entries in $\mm$, and
write $P_i=\bigoplus_{l=1}^{\tau_i}B(-c_{il})$.

\emph{Step 2.} Since $\beta_j^A(P_i[i])=\tau_ia_{j-i}$ and
(D2$^{\mathrm{gr}}$) holds, Step~2 of the proof of
Proposition~\ref{transfer:replacement} applies word for word and gives
\eqref{transfer:tau-two}.

\emph{Step 3.} As there, $P_{\ge3}\simeq L[3]$ with
$L=\coker(P_4\to P_3)$, now graded. Applying $\mathbb D$ to
$0\to P_{<3}\to P\to P_{\ge3}\to0$ and using
$M\simeq\mathbb D(\mathbb D(M))\simeq\mathbb D(P)$ gives the triangle
\eqref{transfer:replacement-triangle} of complexes of graded modules.
Here $\mathbb D(P_i[i])=\bigoplus_l\omega^\bullet(c_{il})[-i]$ has
homology only in degrees $2-i$ and $3-i$, and the argument of Step~3
shows that $\HH_0(\mathbb D(P_{<3}))$ is a quotient of
$\bigoplus_l\HH_2(\omega^\bullet)(c_{2l})$, of dimension at most
$h\tau_2$. Put
\[
 M'=\HH_0(\mathbb D(L[3]))=\Hom_A(L,A(-3)),\qquad
 Z=\HH_{-1}(\mathbb D(L[3]))=\Ext_A^1(L,A(-3)).
\]
The homology sequence of the triangle is an exact sequence of graded
modules $M'\to M\to\HH_0(\mathbb D(P_{<3}))\to Z\to0$, and, as in
Step~3 there, $\mathbb D(L[3])$ has homology only in degrees $0$ and
$-1$. Thus $M'\to M$ has finite-dimensional cokernel, so
$\rk_BM'\ge\rk_BM$, and $\dim_{\Bbbk}Z\le h\tau_2$.

\emph{Step 4.} Apply $\HH^*_{\mm}$ to the triangle
$M'\to\mathbb D(L[3])\to Z[-1]\to M'[1]$. By (D3$^{\mathrm{gr}}$) the
middle term has local cohomology only in degree three, and $Z[-1]$
has local cohomology $Z$ in degree one and zero otherwise, since $Z$
is finite-dimensional. Hence
$\HH^0_{\mm}(M')=\HH^1_{\mm}(M')=0$ and $\HH^2_{\mm}(M')\simeq Z$ as
graded modules. With $\dim_{\Bbbk}Z\le h\tau_2$,
\eqref{transfer:tau-two} and the graded form of
\eqref{transfer:koszul-betti}, this gives the bound.

\emph{Step 5.} Step~5 of the proof of
Proposition~\ref{transfer:replacement} uses only that $A$ is a regular
domain and $B$ a normal domain finite over $A$. It applies to
$M'=\Hom_A(L,A(-3))$ and shows that $M'$ is reflexive.
\end{proof}

\begin{proposition}[Graded Koszul gap]\label{app:graded-gap}
Suppose that a constant $c>0$ satisfies
$\dim_{\Bbbk}\HH^2_{\mm}(N)\ge c\rk_BN$ for every finitely generated
graded reflexive $B$-module $N$. Then $h\ge c$, and every nonzero
finitely generated graded $B$-module $M$ satisfies
\[
 \kappa_{\mathbf x}(M)
 \ge\frac{c}{ehC_{\mathbf x}+c}\,\beta_0^A(M)
 \ge\frac{c}{ehC_{\mathbf x}+c}\,\nu_B(M).
\]
\end{proposition}
\begin{proof}
The proof of Proposition~\ref{transfer:gap} applies, with
Proposition~\ref{app:graded-replacement} in place of
Proposition~\ref{transfer:replacement}: the module $B$ is graded and
reflexive of rank one, and the minimal graded free resolution of $M$
over $A$ has length at most three.
\end{proof}

\addtocontents{toc}{\protect\setcounter{tocdepth}{0}}
\section*{Statements and declarations}
\noindent\emph{Use of generative AI.}
ChatGPT and Claude were used for exploratory computations, literature
search and editorial assistance. The authors are solely responsible for the mathematical arguments and
the final text.

\addtocontents{toc}{\protect\setcounter{tocdepth}{1}}

\end{document}